# SUPER SPECIAL CODES USING SUPER MATRICES

**W. B. Vasantha Kandasamy**
**Florentin Smarandache**
**K. Ilanthenral**

2010

# SUPER SPECIAL CODES USING SUPER MATRICES

# CONTENTS

# PREFACE

The new classes of super special codes are constructed in this book using the specially constructed super special vector spaces. These codes mainly use the super matrices. These codes can be realized as a special type of concatenated codes. This book has four chapters.

In chapter one basic properties of codes and super matrices are given. A new type of super special vector space is constructed in chapter two of this book. Three new classes of super special codes namely, super special row code, super special column code and super special codes are introduced in chapter three. Applications of these codes are given in the final chapter.

These codes will be useful in cryptography, when ARQ protocols are impossible or very costly, in scientific experiments where stage by stage recording of the results are needed, can be used in bulk transmission of information and in medical fields.

The reader should be familiar with both coding theory and in super linear algebras and super matrices.

Our thanks are due to Dr. K. Kandasamy for proof-reading this book. We also acknowledge our gratitude to Kama and Meena for their help with corrections and layout.

W.B.VASANTHA KANDASAMY
FLORENTIN SMARANDACHE
K.ILANTHENRAL

**Chapter One**

# INTRODUCTION TO SUPERMATRICES AND LINEAR CODES

This chapter has two sections. In section we one introduce the basic properties about supermatrices which are essential to build super special codes. Section two gives a brief introduction to algebraic coding theory and the basic properties related with linear codes.

## 1.1 Introduction to Supermatrices

The general rectangular or square array of numbers such as

$$A = \begin{bmatrix} 2 & 3 & 1 & 4 \\ -5 & 0 & 7 & -8 \end{bmatrix},\ B = \begin{bmatrix} 1 & 2 & 3 \\ -4 & 5 & 6 \\ 7 & -8 & 11 \end{bmatrix},$$

$$C = [3, 1, 0, -1, -2] \text{ and } D = \begin{bmatrix} -7/2 \\ 0 \\ \sqrt{2} \\ 5 \\ -41 \end{bmatrix}$$

are known as matrices.

We shall call them as simple matrices [10]. By a simple matrix we mean a matrix each of whose elements are just an ordinary number or a letter that stands for a number. In other words, the elements of a simple matrix are scalars or scalar quantities.

A supermatrix on the other hand is one whose elements are themselves matrices with elements that can be either scalars or other matrices. In general the kind of supermatrices we shall deal with in this book, the matrix elements which have any scalar for their elements. Suppose we have the four matrices;

$$a_{11} = \begin{bmatrix} 2 & -4 \\ 0 & 1 \end{bmatrix},\ a_{12} = \begin{bmatrix} 0 & 40 \\ 21 & -12 \end{bmatrix}$$

$$a_{21} = \begin{bmatrix} 3 & -1 \\ 5 & 7 \\ -2 & 9 \end{bmatrix} \text{ and } a_{22} = \begin{bmatrix} 4 & 12 \\ -17 & 6 \\ 3 & 11 \end{bmatrix}.$$

One can observe the change in notation $a_{ij}$ denotes a matrix and not a scalar of a matrix $(1 \leq i, j \leq 2)$.

Let

$$a = \begin{bmatrix} a_{11} & a_{12} \\ a_{21} & a_{22} \end{bmatrix};$$

we can write out the matrix a in terms of the original matrix elements i.e.,

$$a = \left[\begin{array}{cc|cc} 2 & -4 & 0 & 40 \\ 0 & 1 & 21 & -12 \\ \hline 3 & -1 & 4 & 12 \\ 5 & 7 & -17 & 6 \\ -2 & 9 & 3 & 11 \end{array}\right].$$

Here the elements are divided vertically and horizontally by thin lines. If the lines were not used the matrix a would be read as a simple matrix.

Thus far we have referred to the elements in a supermatrix as matrices as elements. It is perhaps more usual to call the elements of a supermatrix as submatrices. We speak of the submatrices within a supermatrix. Now we proceed on to define the order of a supermatrix.

The order of a supermatrix is defined in the same way as that of a simple matrix. The height of a supermatrix is the number of rows of submatrices in it. The width of a supermatrix is the number of columns of submatrices in it.

All submatrices with in a given row must have the same number of rows. Likewise all submatrices with in a given column must have the same number of columns.

A diagrammatic representation is given by the following figure:

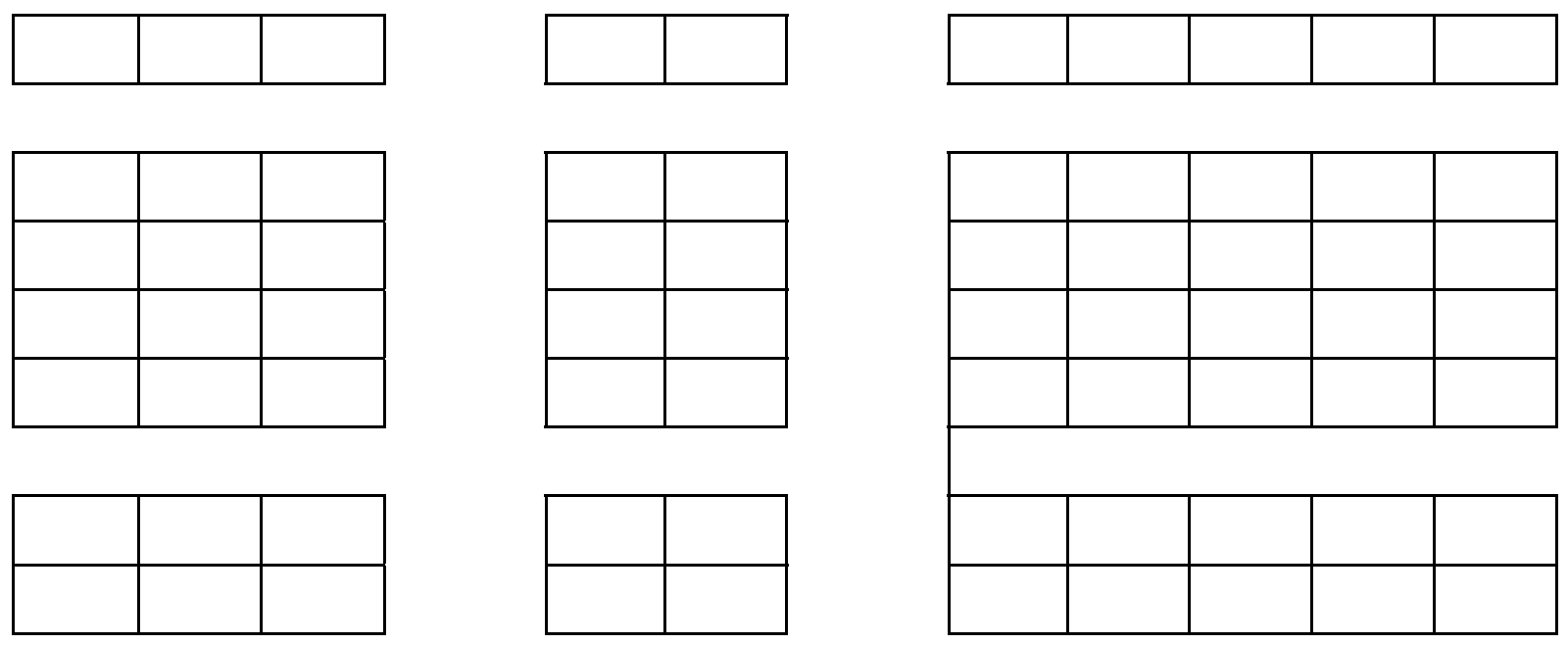

In the first row of rectangles we have one row of a square for each rectangle; in the second row of rectangles we have four rows of squares for each rectangle and in the third row of rectangles we have two rows of squares for each rectangle. Similarly for the first column of rectangles three columns of squares for each rectangle. For the second column of rectangles we have two column of squares for each rectangle, and for the third column of rectangles we have five columns of squares for each rectangle.

Thus we have for this supermatrix 3 rows and 3 columns.

One thing should now be clear from the definition of a supermatrix. The super order of a supermatrix tells us nothing about the simple order of the matrix from which it was obtained by partitioning. Furthermore, the order of supermatrix tells us nothing about the orders of the submatrices within that supermatrix.

Now we illustrate the number of rows and columns of a supermatrix.

***Example 1.1.1:*** Let

$$a = \left[\begin{array}{c|cccc} 3 & 3 & 0 & 1 & 4 \\ -1 & 2 & 1 & -1 & 6 \\ \hline 0 & 3 & 4 & 5 & 6 \\ 1 & 7 & 8 & -9 & 0 \\ 2 & 1 & 2 & 3 & -4 \end{array}\right].$$

a is a supermatrix with two rows and two columns.

Now we proceed on to define the notion of partitioned matrices. It is always possible to construct a supermatrix from any simple matrix that is not a scalar quantity.

The supermatrix can be constructed from a simple matrix this process of constructing supermatrix is called the partitioning.

A simple matrix can be partitioned by dividing or separating the matrix between certain specified rows, or the procedure may be reversed. The division may be made first between rows and then between columns.

We illustrate this by a simple example.

***Example 1.1.2:*** Let

$$A = \begin{bmatrix} 3 & 0 & 1 & 1 & 2 & 0 \\ 1 & 0 & 0 & 3 & 5 & 2 \\ 5 & -1 & 6 & 7 & 8 & 4 \\ 0 & 9 & 1 & 2 & 0 & -1 \\ 2 & 5 & 2 & 3 & 4 & 6 \\ 1 & 6 & 1 & 2 & 3 & 9 \end{bmatrix}$$

is a 6 × 6 simple matrix with real numbers as elements.

$$A_1 = \left[\begin{array}{cc|cccc} 3 & 0 & 1 & 1 & 2 & 0 \\ 1 & 0 & 0 & 3 & 5 & 2 \\ 5 & -1 & 6 & 7 & 8 & 4 \\ 0 & 9 & 1 & 2 & 0 & -1 \\ 2 & 5 & 2 & 3 & 4 & 6 \\ 1 & 6 & 1 & 2 & 3 & 9 \end{array}\right].$$

Now let us draw a thin line between the 2nd and 3rd columns.

This gives us the matrix $A_1$. Actually $A_1$ may be regarded as a supermatrix with two matrix elements forming one row and two columns.

Now consider

$$A_2 = \left[\begin{array}{cccccc} 3 & 0 & 1 & 1 & 2 & 0 \\ 1 & 0 & 0 & 3 & 5 & 2 \\ 5 & -1 & 6 & 7 & 8 & 4 \\ 0 & 9 & 1 & 2 & 0 & -1 \\ \hline 2 & 5 & 2 & 3 & 4 & 6 \\ 1 & 6 & 1 & 2 & 3 & 9 \end{array}\right].$$

Draw a thin line between the rows 4 and 5 which gives us the new matrix $A_2$. $A_2$ is a supermatrix with two rows and one column.

Now consider the matrix

$$A_3 = \left[\begin{array}{cc|cccc} 3 & 0 & 1 & 1 & 2 & 0 \\ 1 & 0 & 0 & 3 & 5 & 2 \\ 5 & -1 & 6 & 7 & 8 & 4 \\ 0 & 9 & 1 & 2 & 0 & -1 \\ \hline 2 & 5 & 2 & 3 & 4 & 6 \\ 1 & 6 & 1 & 2 & 3 & 9 \end{array}\right],$$

$A_3$ is now a second order supermatrix with two rows and two columns. We can simply write $A_3$ as

$$\begin{bmatrix} a_{11} & a_{12} \\ a_{21} & a_{22} \end{bmatrix}$$

where

$$a_{11} = \begin{bmatrix} 3 & 0 \\ 1 & 0 \\ 5 & -1 \\ 0 & 9 \end{bmatrix},$$

$$a_{12} = \begin{bmatrix} 1 & 1 & 2 & 0 \\ 0 & 3 & 5 & 2 \\ 6 & 7 & 8 & 4 \\ 1 & 2 & 0 & -1 \end{bmatrix},$$

$$a_{21} = \begin{bmatrix} 2 & 5 \\ 1 & 6 \end{bmatrix} \text{ and } a_{22} = \begin{bmatrix} 2 & 3 & 4 & 6 \\ 1 & 2 & 3 & 9 \end{bmatrix}.$$

The elements now are the submatrices defined as $a_{11}$, $a_{12}$, $a_{21}$ and $a_{22}$ and therefore $A_3$ is in terms of letters.

According to the methods we have illustrated a simple matrix can be partitioned to obtain a supermatrix in any way that happens to suit our purposes.

The natural order of a supermatrix is usually determined by the natural order of the corresponding simple matrix. Further more we are not usually concerned with natural order of the submatrices within a supermatrix.

Now we proceed on to recall the notion of symmetric partition, for more information about these concepts please refer [10]. By a symmetric partitioning of a matrix we mean that the rows and columns are partitioned in exactly the same way. If the matrix is partitioned between the first and second column and between the third and fourth column, then to be symmetrically partitioning, it must also be partitioned between the first and second rows and third and fourth rows. According to this rule of symmetric partitioning only square simple matrix can be symmetrically partitioned. We give an example of a symmetrically partitioned matrix $a_s$,

***Example 1.1.3:*** Let

$$a_s = \left[\begin{array}{c|cc|c} 2 & 3 & 4 & 1 \\ \hline 5 & 6 & 9 & 2 \\ 0 & 6 & 1 & 9 \\ \hline 5 & 1 & 1 & 5 \end{array}\right].$$

Here we see that the matrix has been partitioned between columns one and two and three and four. It has also been partitioned between rows one and two and rows three and four.

Now we just recall from [10] the method of symmetric partitioning of a symmetric simple matrix.

***Example 1.1.4:*** Let us take a fourth order symmetric matrix and partition it between the second and third rows and also between the second and third columns.

$$a = \left[\begin{array}{cc|cc} 4 & 3 & 2 & 7 \\ 3 & 6 & 1 & 4 \\ \hline 2 & 1 & 5 & 2 \\ 7 & 4 & 2 & 7 \end{array}\right].$$

We can represent this matrix as a supermatrix with letter elements.

$$a_{11} = \begin{bmatrix} 4 & 3 \\ 3 & 6 \end{bmatrix}, a_{12} = \begin{bmatrix} 2 & 7 \\ 1 & 4 \end{bmatrix}$$

$$a_{21} = \begin{bmatrix} 2 & 1 \\ 7 & 4 \end{bmatrix} \text{ and } a_{22} = \begin{bmatrix} 5 & 2 \\ 2 & 7 \end{bmatrix},$$

so that

$$a = \begin{bmatrix} a_{11} & a_{12} \\ a_{21} & a_{22} \end{bmatrix}.$$

The diagonal elements of the supermatrix a are $a_{11}$ and $a_{22}$. We also observe the matrices $a_{11}$ and $a_{22}$ are also symmetric matrices.

The non diagonal elements of this supermatrix a are the matrices $a_{12}$ and $a_{21}$. Clearly $a_{21}$ is the transpose of $a_{12}$.

The simple rule about the matrix element of a symmetrically partitioned symmetric simple matrix are (1) The diagonal submatrices of the supermatrix are all symmetric matrices. (2) The matrix elements below the diagonal are the transposes of the corresponding elements above the diagonal.

The forth order supermatrix obtained from a symmetric partitioning of a symmetric simple matrix a is as follows:

$$a = \begin{bmatrix} a_{11} & a_{12} & a_{13} & a_{14} \\ a'_{12} & a_{22} & a_{23} & a_{24} \\ a'_{13} & a'_{23} & a_{33} & a_{34} \\ a'_{14} & a'_{24} & a'_{34} & a_{44} \end{bmatrix}.$$

How to express that a symmetric matrix has been symmetrically partitioned (i) $a_{11}$ and $a^t_{11}$ are equal. (ii) $a^t_{ij}$ $(i \neq j)$; $a^t_{ij} = a_{ji}$ and $a^t_{ji} = a_{ij}$. Thus the general expression for a symmetrically partitioned symmetric matrix;

$$a = \begin{bmatrix} a_{11} & a_{12} & \cdots & a_{1n} \\ a'_{12} & a_{22} & \cdots & a_{2n} \\ \vdots & \vdots & & \vdots \\ a'_{1n} & a'_{2n} & \cdots & a_{nn} \end{bmatrix}.$$

If we want to indicate a symmetrically partitioned simple diagonal matrix we would write

$$D = \begin{bmatrix} D_1 & 0 & \ldots & 0 \\ 0' & D_2 & \ldots & 0 \\ & & & \\ 0' & 0' & \ldots & D_n \end{bmatrix}$$

0' only represents the order is reversed or transformed. We denote $a_{ij}^{t}$ = a'$_{ij}$ just the ' means the transpose.

D will be referred to as the super diagonal matrix. The identity matrix

$$I = \begin{bmatrix} I_s & 0 & 0 \\ 0 & I_t & 0 \\ 0 & 0 & I_r \end{bmatrix}$$

s, t and r denote the number of rows and columns of the first second and third identity matrices respectively (zeros denote matrices with zero as all entries).

***Example 1.1.5:*** We just illustrate a general super diagonal matrix d;

$$d = \left[\begin{array}{ccc|cc} 3 & 1 & 2 & 0 & 0 \\ 5 & 6 & 0 & 0 & 0 \\ \hline 0 & 0 & 0 & 2 & 5 \\ 0 & 0 & 0 & -1 & 3 \\ 0 & 0 & 0 & 9 & 10 \end{array}\right]$$

$$\text{i.e., } d = \begin{bmatrix} m_1 & 0 \\ 0 & m_2 \end{bmatrix}.$$

An example of a super diagonal matrix with vector elements is given, which can be useful in experimental designs.

***Example 1.1.6:*** Let

$$\left[\begin{array}{c|c|c|c} 1 & 0 & 0 & 0 \\ 1 & 0 & 0 & 0 \\ 1 & 0 & 0 & 0 \\ \hline 0 & 1 & 0 & 0 \\ 0 & 1 & 0 & 0 \\ \hline 0 & 0 & 1 & 0 \\ 0 & 0 & 1 & 0 \\ 0 & 0 & 1 & 0 \\ \hline 0 & 0 & 1 & 0 \\ \hline 0 & 0 & 0 & 1 \\ 0 & 0 & 0 & 1 \\ 0 & 0 & 0 & 1 \\ 0 & 0 & 0 & 1 \end{array}\right].$$

Here the diagonal elements are only column unit vectors. In case of supermatrix [10] has defined the notion of partial triangular matrix as a supermatrix.

***Example 1.1.7:*** Let

$$u = \left[\begin{array}{ccc|cc} 2 & 1 & 1 & 3 & 2 \\ 0 & 5 & 2 & 1 & 1 \\ 0 & 0 & 1 & 0 & 2 \end{array}\right]$$

u is a partial upper triangular supermatrix.

***Example 1.1.8:*** Let

$$L = \left[\begin{array}{ccccc} 5 & 0 & 0 & 0 & 0 \\ 7 & 2 & 0 & 0 & 0 \\ 1 & 2 & 3 & 0 & 0 \\ 4 & 5 & 6 & 7 & 0 \\ 1 & 2 & 5 & 2 & 6 \\ \hline 1 & 2 & 3 & 4 & 5 \\ 0 & 1 & 0 & 1 & 0 \end{array}\right];$$

L is partial upper triangular matrix partitioned as a supermatrix.

Thus T = $\left[\begin{array}{c} T \\ \hline a' \end{array}\right]$ where T is the lower triangular submatrix, with

$$T = \begin{bmatrix} 5 & 0 & 0 & 0 & 0 \\ 7 & 2 & 0 & 0 & 0 \\ 1 & 2 & 3 & 0 & 0 \\ 4 & 5 & 6 & 7 & 0 \\ 1 & 2 & 5 & 2 & 6 \end{bmatrix} \text{ and } a' = \begin{bmatrix} 1 & 2 & 3 & 4 & 5 \\ 0 & 1 & 0 & 1 & 0 \end{bmatrix}.$$

We proceed on to define the notion of supervectors i.e., Type I column supervector. A simple vector is a vector each of whose elements is a scalar. It is nice to see the number of different types of supervectors given by [10].

***Example 1.1.9:*** Let

$$v = \left[\begin{array}{c} 1 \\ 3 \\ 4 \\ \hline 5 \\ 7 \end{array}\right].$$

This is a type I i.e., type one column supervector.

$$v = \begin{bmatrix} v_1 \\ v_2 \\ \vdots \\ v_n \end{bmatrix}$$

where each $v_i$ is a column subvectors of the column vector v.

Type I row supervector is given by the following example.

***Example 1.1.10:*** $v^1$ = [2 3 1 | 5 7 8 4] is a type I row supervector. i.e., v' = [$v'_1$, $v'_2$, …, $v'_n$] where each $v'_i$ is a row subvector; $1 \leq i \leq n$.

Next we recall the definition of type II supervectors.

Type II column supervectors.

**DEFINITION 1.1.1:** *Let*

$$a = \begin{bmatrix} a_{11} & a_{12} & \dots & a_{1m} \\ a_{21} & a_{22} & \dots & a_{2m} \\ \dots & \dots & \dots & \dots \\ a_{n1} & a_{n2} & \dots & a_{nm} \end{bmatrix}$$

$$\begin{aligned} a_1^1 &= [a_{11} \dots a_{1m}] \\ a_2^1 &= [a_{21} \dots a_{2m}] \\ &\dots \\ a_n^1 &= [a_{n1} \dots a_{nm}] \end{aligned}$$

*i.e.,*

$$a = \begin{bmatrix} a_1^1 \\ a_2^1 \\ \vdots \\ a_n^1 \end{bmatrix}_m$$

*is defined to be the type II column supervector.*
*Similarly if*

$$a^1 = \begin{bmatrix} a_{11} \\ a_{21} \\ \vdots \\ a_{n1} \end{bmatrix}, \; a^2 = \begin{bmatrix} a_{12} \\ a_{22} \\ \vdots \\ a_{n2} \end{bmatrix}, \dots, \; a^m = \begin{bmatrix} a_{1m} \\ a_{2m} \\ \vdots \\ a_{nm} \end{bmatrix}.$$

*Hence now* $a = [a^1\ a^2 \dots a^m]_n$ *is defined to be the type II row supervector.*

*Clearly*

$$a = \begin{bmatrix} a_1^1 \\ a_2^1 \\ \vdots \\ a_n^1 \end{bmatrix}_m = [a^1\ a^2 \ldots a^m]_n$$

*the equality of supermatrices.*

***Example 1.1.11:*** Let

$$A = \begin{bmatrix} 3 & 6 & 0 & 4 & 5 \\ 2 & 1 & 6 & 3 & 0 \\ 1 & 1 & 1 & 2 & 1 \\ 0 & 1 & 0 & 1 & 0 \\ 2 & 0 & 1 & 2 & 1 \end{bmatrix}$$

be a simple matrix. Let a and b the supermatrix made from A.

$$a = \left[\begin{array}{ccc|cc} 3 & 6 & 0 & 4 & 5 \\ 2 & 1 & 6 & 3 & 0 \\ 1 & 1 & 1 & 2 & 1 \\ \hline 0 & 1 & 0 & 1 & 0 \\ 2 & 0 & 1 & 2 & 1 \end{array}\right]$$

where

$$a_{11} = \begin{bmatrix} 3 & 6 & 0 \\ 2 & 1 & 6 \\ 1 & 1 & 1 \end{bmatrix}, \; a_{12} = \begin{bmatrix} 4 & 5 \\ 3 & 0 \\ 2 & 1 \end{bmatrix},$$

$$a_{21} = \begin{bmatrix} 0 & 1 & 0 \\ 2 & 0 & 1 \end{bmatrix} \text{ and } a_{22} = \begin{bmatrix} 1 & 0 \\ 2 & 1 \end{bmatrix}.$$

i.e., $$a = \begin{bmatrix} a_{11} & a_{12} \\ a_{21} & a_{22} \end{bmatrix}.$$

$$b = \left[\begin{array}{cccc|c} 3 & 6 & 0 & 4 & 5 \\ 2 & 1 & 6 & 3 & 0 \\ 1 & 1 & 1 & 2 & 1 \\ 0 & 1 & 0 & 1 & 0 \\ \hline 2 & 0 & 1 & 2 & 1 \end{array}\right] = \begin{bmatrix} b_{11} & b_{12} \\ b_{21} & b_{22} \end{bmatrix}$$

where

$$b_{11} = \begin{bmatrix} 3 & 6 & 0 & 4 \\ 2 & 1 & 6 & 3 \\ 1 & 1 & 1 & 2 \\ 0 & 1 & 0 & 1 \end{bmatrix}, \; b_{12} = \begin{bmatrix} 5 \\ 0 \\ 1 \\ 0 \end{bmatrix},$$

$b_{21}$ = [2 0 1 2 ] and $b_{22}$ = [1].

$$a = \left[\begin{array}{ccc|cc} 3 & 6 & 0 & 4 & 5 \\ 2 & 1 & 6 & 3 & 0 \\ 1 & 1 & 1 & 2 & 1 \\ \hline 0 & 1 & 0 & 1 & 0 \\ 2 & 0 & 1 & 2 & 1 \end{array}\right]$$

and

$$b = \left[\begin{array}{cccc|c} 3 & 6 & 0 & 4 & 5 \\ 2 & 1 & 6 & 3 & 0 \\ 1 & 1 & 1 & 2 & 1 \\ 0 & 1 & 0 & 1 & 0 \\ \hline 2 & 0 & 1 & 2 & 1 \end{array}\right].$$

We see that the corresponding scalar elements for matrix a and matrix b are identical. Thus two supermatrices are equal if and only if their corresponding simple forms are equal.

Now we give examples of type III supervector for more refer [10].

***Example 1.1.12:***

$$a = \left[\begin{array}{ccc|cc} 3 & 2 & 1 & 7 & 8 \\ 0 & 2 & 1 & 6 & 9 \\ 0 & 0 & 5 & 1 & 2 \end{array}\right] = [T' \mid a']$$

and

$$b = \left[\begin{array}{ccc} 2 & 0 & 0 \\ 9 & 4 & 0 \\ 8 & 3 & 6 \\ \hline 5 & 2 & 9 \\ 4 & 7 & 3 \end{array}\right] = \left[\begin{array}{c} T \\ \hline b' \end{array}\right]$$

are type III supervectors.

One interesting and common example of a type III supervector is a prediction data matrix having both predictor and criterion attributes.

The next interesting notion about supermatrix is its transpose. First we illustrate this by an example before we give the general case.

***Example 1.1.13:*** Let

$$a = \left[\begin{array}{ccc|cc} 2 & 1 & 3 & 5 & 6 \\ 0 & 2 & 0 & 1 & 1 \\ 1 & 1 & 1 & 0 & 2 \\ \hline 2 & 2 & 0 & 1 & 1 \\ 5 & 6 & 1 & 0 & 1 \\ \hline 2 & 0 & 0 & 0 & 4 \\ 1 & 0 & 1 & 1 & 5 \end{array}\right]$$

$$= \begin{bmatrix} a_{11} & a_{12} \\ a_{21} & a_{22} \\ a_{31} & a_{32} \end{bmatrix}$$

where

$$a_{11} = \begin{bmatrix} 2 & 1 & 3 \\ 0 & 2 & 0 \\ 1 & 1 & 1 \end{bmatrix}, a_{12} = \begin{bmatrix} 5 & 6 \\ 1 & 1 \\ 0 & 2 \end{bmatrix},$$

$$a_{21} = \begin{bmatrix} 2 & 2 & 0 \\ 5 & 6 & 1 \end{bmatrix}, a_{22} = \begin{bmatrix} 1 & 1 \\ 0 & 1 \end{bmatrix},$$

$$a_{31} = \begin{bmatrix} 2 & 0 & 0 \\ 1 & 0 & 1 \end{bmatrix} \text{ and } a_{32} = \begin{bmatrix} 0 & 4 \\ 1 & 5 \end{bmatrix}.$$

The transpose of a

$$a^t = a' = \left[\begin{array}{ccc|cc|cc} 2 & 0 & 1 & 2 & 5 & 2 & 1 \\ 1 & 2 & 1 & 2 & 6 & 0 & 0 \\ 3 & 0 & 1 & 0 & 1 & 0 & 1 \\ \hline 5 & 1 & 0 & 1 & 0 & 0 & 1 \\ 6 & 1 & 2 & 1 & 1 & 4 & 5 \end{array}\right].$$

Let us consider the transposes of $a_{11}$, $a_{12}$, $a_{21}$, $a_{22}$, $a_{31}$ and $a_{32}$.

$$a'_{11} = a_{11}^t = \begin{bmatrix} 2 & 0 & 1 \\ 1 & 2 & 1 \\ 3 & 0 & 1 \end{bmatrix}$$

$$a'_{12} = a_{12}^t = \begin{bmatrix} 5 & 1 & 0 \\ 6 & 1 & 2 \end{bmatrix}$$

$$a'_{21} = a_{21}^t = \begin{bmatrix} 2 & 5 \\ 2 & 6 \\ 0 & 1 \end{bmatrix}$$

$$a'_{31} = a^t_{31} = \begin{bmatrix} 2 & 1 \\ 0 & 0 \\ 0 & 1 \end{bmatrix}$$

$$a'_{22} = a^t_{22} = \begin{bmatrix} 1 & 0 \\ 1 & 1 \end{bmatrix}$$

$$a'_{32} = a^t_{32} = \begin{bmatrix} 0 & 1 \\ 4 & 5 \end{bmatrix}.$$

$$a' = \begin{bmatrix} a'_{11} & a'_{21} & a'_{31} \\ a'_{12} & a'_{22} & a'_{32} \end{bmatrix}.$$

Now we describe the general case. Let

$$a = \begin{bmatrix} a_{11} & a_{12} & \cdots & a_{1m} \\ a_{21} & a_{22} & \cdots & a_{2m} \\ \vdots & \vdots & & \vdots \\ a_{n1} & a_{n2} & \cdots & a_{nm} \end{bmatrix}$$

be a $n \times m$ supermatrix. The transpose of the supermatrix a denoted by

$$a' = \begin{bmatrix} a'_{11} & a'_{21} & \cdots & a'_{n1} \\ a'_{12} & a'_{22} & \cdots & a'_{n2} \\ \vdots & \vdots & & \vdots \\ a'_{1m} & a'_{2m} & \cdots & a'_{nm} \end{bmatrix}.$$

a' is a m by n supermatrix obtained by taking the transpose of each element i.e., the submatrices of a.

Now we will find the transpose of a symmetrically partitioned symmetric simple matrix. Let a be the symmetrically partitioned symmetric simple matrix.

Let a be a m × m symmetric supermatrix i.e.,

$$a = \begin{bmatrix} a_{11} & a_{21} & \cdots & a_{m1} \\ a_{12} & a_{22} & \cdots & a_{m2} \\ \vdots & \vdots & & \vdots \\ a_{1m} & a_{2m} & \cdots & a_{mm} \end{bmatrix}$$

the transpose of the supermatrix is given by a'

$$a' = \begin{bmatrix} a'_{11} & (a'_{12})' & \cdots & (a'_{1m})' \\ a'_{12} & a'_{22} & \cdots & (a'_{2m})' \\ \vdots & \vdots & & \vdots \\ a'_{1m} & a'_{2m} & \cdots & a'_{mm} \end{bmatrix}$$

The diagonal matrix $a_{11}$ are symmetric matrices so are unaltered by transposition. Hence

$$a'_{11} = a_{11}, a'_{22} = a_{22}, \ldots, a'_{mm} = a_{mm}.$$

Recall also the transpose of a transpose is the original matrix. Therefore

$$(a'_{12})' = a_{12}, (a'_{13})' = a_{13}, \ldots, (a'_{ij})' = a_{ij}.$$

Thus the transpose of supermatrix constructed by symmetrically partitioned symmetric simple matrix a of a' is given by

$$a' = \begin{bmatrix} a_{11} & a_{12} & \cdots & a_{1m} \\ a'_{21} & a_{22} & \cdots & a_{2m} \\ \vdots & \vdots & & \vdots \\ a'_{1m} & a'_{2m} & \cdots & a_{mm} \end{bmatrix}.$$

Thus a = a'.
Similarly transpose of a symmetrically partitioned diagonal matrix is simply the original diagonal supermatrix itself;

i.e., if

$$D = \begin{bmatrix} d_1 & & & \\ & d_2 & & \\ & & \ddots & \\ & & & d_n \end{bmatrix}$$

$$D' = \begin{bmatrix} d'_1 & & & \\ & d'_2 & & \\ & & \ddots & \\ & & & d'_n \end{bmatrix}$$

$d'_1 = d_1$, $d'_2 = d_2$ etc. Thus D = D'.

Now we see the transpose of a type I supervector.

***Example 1.1.14:*** Let

$$V = \left[\begin{array}{c} 3 \\ 1 \\ 2 \\ \hline 4 \\ 5 \\ 7 \\ \hline 5 \\ 1 \end{array}\right]$$

The transpose of V denoted by V' or $V^t$ is

$$V' = [3\ 1\ 2\ |\ 4\ 5\ 7\ |\ 5\ 1].$$

If

$$V = \begin{bmatrix} v_1 \\ v_2 \\ v_3 \end{bmatrix}$$

where

$$v_1 = \begin{bmatrix} 3 \\ 1 \\ 2 \end{bmatrix}, v_2 = \begin{bmatrix} 4 \\ 5 \\ 7 \end{bmatrix} \text{ and } v_3 = \begin{bmatrix} 5 \\ 1 \end{bmatrix}$$

$$V' = [v'_1\ v'_2\ v'_3].$$

Thus if

$$V = \begin{bmatrix} v_1 \\ v_2 \\ \vdots \\ v_n \end{bmatrix}$$

then

$$V' = [v'_1\ v'_2 \ldots v'_n].$$

***Example 1.1.15:*** Let

$$t = \left[\begin{array}{cccc|cc} 3 & 0 & 1 & 1 & 5 & 2 \\ 4 & 2 & 0 & 1 & 3 & 5 \\ 1 & 0 & 1 & 0 & 1 & 6 \end{array}\right]$$

= [T | a ]. The transpose of t

$$\text{i.e., } t' = \left[\begin{array}{ccc} 3 & 4 & 1 \\ 0 & 2 & 0 \\ 1 & 0 & 1 \\ 1 & 1 & 0 \\ \hline 5 & 3 & 1 \\ 2 & 5 & 6 \end{array}\right] = \left[\begin{array}{c} T' \\ \hline a' \end{array}\right].$$

## 1.2 Introduction of Linear Codes and their Properties

In this section we just recall the definition of linear code and enumerate a few important properties about them. We begin by describing a simple model of a communication transmission system given by the figure 1.2.

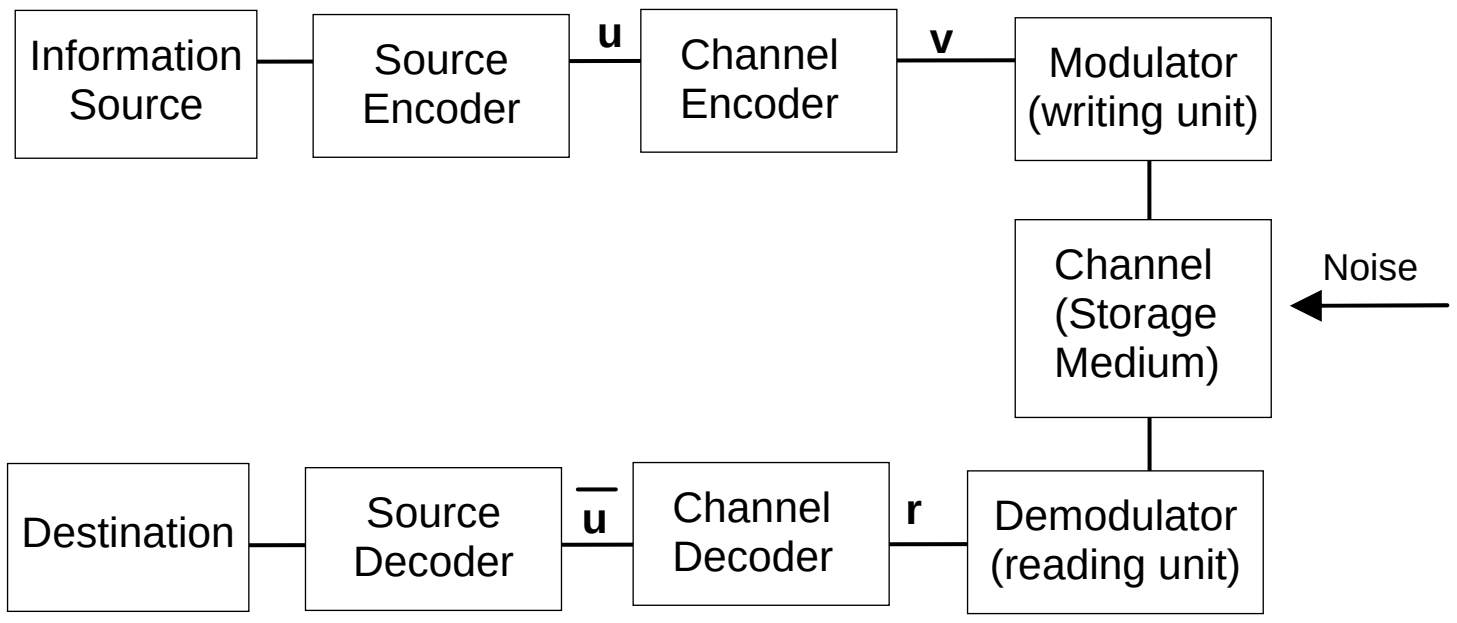

**Figure 1.2: General Coding System**

Messages go through the system starting from the source (sender). We shall only consider senders with a finite number of discrete signals (eg. Telegraph) in contrast to continuous sources (eg. Radio). In most systems the signals emanating from the source cannot be transmitted directly by the channel. For instance, a binary channel cannot transmit words in the usual Latin alphabet. Therefore an encoder performs the important task of data reduction and suitably transforms the message into usable form. Accordingly one distinguishes between source encoding the channel encoding. The former reduces the message to its essential(recognizable) parts, the latter adds redundant information to enable detection and correction of possible errors in the transmission. Similarly on the receiving end one distinguishes between channel decoding and source decoding, which invert the corresponding channel and source encoding besides detecting and correcting errors.

One of the main aims of coding theory is to design methods for transmitting messages error free cheap and as fast as possible. There is of course the possibility of repeating the message. However this is time consuming, inefficient and crude. We also note that the possibility of errors increases with an increase in the length of messages. We want to find efficient algebraic methods (codes) to improve the reliability of the transmission of messages. There are many types of algebraic codes; here we give a few of them.

Throughout this book we assume that only finite fields represent the underlying alphabet for coding. Coding consists of transforming a block of k message symbols $a_1, a_2, \ldots, a_k$; $a_i \in F_q$ into a code word $x = x_1\ x_2 \ldots x_n$; $x_i \in F_q$, where $n \geq k$. Here the first $k_i$ symbols are the message symbols i.e., $x_i = a_i$; $1 \leq i \leq k$; the remaining $n - k$ elements $x_{k+1}, x_{k+2}, \ldots, x_n$ are check symbols or control symbols. Code words will be written in one of the forms x; $x_1, x_2, \ldots, x_n$ or $(x_1\ x_2 \ldots x_n)$ or $x_1\ x_2 \ldots x_n$. The check symbols can be obtained from the message symbols in such a way that the code words x satisfy a system of linear equations; $Hx^T = (0)$ where H is the given $(n - k) \times n$ matrix with elements in $F_q = Z_p^n$ $(q = p^n)$. A standard form for H is $(A, I_{n-k})$ with $n - k \times k$ matrix and $I_{n-k}$, the $n - k \times n - k$ identity matrix.

We illustrate this by the following example:

***Example 1.2.1:*** Let us consider $Z_2 = \{0, 1\}$. Take n = 7, k = 3. The message $a_1\ a_2\ a_3$ is encoded as the code word $x = a_1\ a_2\ a_3\ x_4\ x_5\ x_6\ x_7$. Here the check symbols $x_4\ x_5\ x_6\ x_7$ are such that for this given matrix

$$H = \begin{bmatrix} 0 & 1 & 0 & 1 & 0 & 0 & 0 \\ 1 & 0 & 1 & 0 & 1 & 0 & 0 \\ 0 & 0 & 1 & 0 & 0 & 1 & 0 \\ 0 & 0 & 1 & 0 & 0 & 0 & 1 \end{bmatrix} = (A; I_4);$$

we have $Hx^T = (0)$ where $x = a_1\ a_2\ a_3\ x_4\ x_5\ x_6\ x_7$.

$$a_2 + x_4 = 0;\ a_1 + a_3 + x_5 = 0;\ a_3 + x_6 = 0;\ a_3 + x_7 = 0.$$

Thus the check symbols $x_4\ x_5\ x_6\ x_7$ are determined by $a_1\ a_2\ a_3$. The equation $Hx^T = (0)$ are also called check equations. If the message a = 1 0 0 then, $x_4 = 0$, $x_5 = 1$, $x_6 = 0$ and $x_7 = 0$. The code word x is 1 0 0 0 1 0 0. If the message a = 1 1 0 then $x_4$ =1, $x_5$ = 1, $x_6$ = 1 = $x_7$. Thus the code word x = 1 1 0 1 1 0 0.

We will have altogether $2^3$ code words given by

| | |
|---|---|
| 0 0 0 0 0 0 0 | 1 1 0 1 1 0 0 |
| 1 0 0 0 1 0 0 | 1 0 1 0 0 1 1 |
| 0 1 0 1 1 0 0 | 0 1 1 1 1 1 1 |
| 0 0 1 0 1 1 1 | 1 1 1 1 0 1 1 |

**DEFINITION 1.2.1:** *Let H be an n – k × n matrix with elements in $Z_q$. The set of all n-dimensional vectors satisfying $Hx^T = (0)$ over $Z_q$ is called a linear code(block code) C over $Z_q$ of block length n. The matrix H is called the parity check matrix of the code C. C is also called a linear(n, k) code.*

*If H is of the form(A, $I_{n-k}$) then the k-symbols of the code word x is called massage(or information) symbols and the last n – k symbols in x are the check symbols. C is then also called a systematic linear(n, k) code. If q = 2, then C is a binary code. k/n is called transmission (or information) rate.*

*The set C of solutions of x of $Hx^T = (0)$. i.e., the solution space of this system of equations, forms a subspace of this system of equations, forms a subspace of $Z_q^n$ of dimension k. Since the code words form an additive group, C is also called a group code. C can also be regarded as the null space of the matrix H.*

***Example 1.2.2:*** (Repetition Code) If each codeword of a code consists of only one message symbol $a_1 \in Z_2$ and (n – 1) check symbols $x_2 = x_3 = \ldots = x_n$ are all equal to $a_1$ ($a_1$ is repeated n – 1 times) then we obtain a binary (n, 1) code with parity check matrix

$$H = \begin{bmatrix} 1 & 1 & 0 & 0 & \ldots & 1 \\ 0 & 0 & 1 & 0 & \ldots & 0 \\ 0 & 0 & 0 & 1 & \ldots & 0 \\ \vdots & \vdots & \vdots & \vdots & & \vdots \\ 1 & 0 & 0 & 0 & \ldots & 1 \end{bmatrix}.$$

There are only two code words in this code namely 0 0 … 0 and 1 1 …1.

If is often impracticable, impossible or too expensive to send the original message more than once. Especially in the transmission of information from satellite or other spacecraft, it is impossible to repeat such messages owing to severe time limitations. One such cases is the photograph from spacecraft as it is moving it may not be in a position to retrace its path. In such cases it is impossible to send the original message more than once. In repetition codes we can of course also consider code words with more than one message symbol.

***Example 1.2.3:*** (Parity-Check Code): This is a binary (n, n – 1) code with parity-check matrix to be H = (1 1 … 1). Each code word has one check symbol and all code words are given by all binary vectors of length n with an even number of ones. Thus if sum of the ones of a code word which is received is odd then atleast one error must have occurred in the transmission.

Such codes find its use in banking. The last digit of the account number usually is a control digit.

**DEFINITION 1.2.2:** *The matrix $G = (I_k, -A^T)$ is called a canonical generator matrix (or canonical basic matrix or encoding matrix) of a linear (n, k) code with parity check matrix $H = (A, I_{n-k})$. In this case we have $GH^T = (0)$.*

***Example 1.2.4:*** Let

$$G = \begin{bmatrix} 1 & 0 & 0 & 0 & 1 & 0 & 0 \\ 0 & 1 & 0 & 1 & 0 & 0 & 0 \\ 0 & 0 & 1 & 0 & 1 & 1 & 1 \end{bmatrix}$$

be the canonical generator matrix of the code given in example 1.2.1. The $2^3$ code words x of the binary code can be obtained from x = aG with $a = a_1 a_2 a_3$, $a_i \in Z_2$, $1 \le i \le 3$. We have the set of $a = a_1 a_2 a_3$ which correspond to the message symbols which is as follows:

[0 0 0], [1 0 0], [0 1 0], [0 0 1],
[1 1 0], [1 0 1], [0 1 1] and [1 1 1].

$$x = \begin{bmatrix} 0 & 0 & 0 \end{bmatrix} \begin{bmatrix} 1 & 0 & 0 & 0 & 1 & 0 & 0 \\ 0 & 1 & 0 & 1 & 0 & 0 & 0 \\ 0 & 0 & 1 & 0 & 1 & 1 & 1 \end{bmatrix}$$
$$= \begin{bmatrix} 0 & 0 & 0 & 0 & 0 & 0 & 0 \end{bmatrix}$$

$$x = \begin{bmatrix} 1 & 0 & 0 \end{bmatrix} \begin{bmatrix} 1 & 0 & 0 & 0 & 1 & 0 & 0 \\ 0 & 1 & 0 & 1 & 0 & 0 & 0 \\ 0 & 0 & 1 & 0 & 1 & 1 & 1 \end{bmatrix}$$
$$= \begin{bmatrix} 1 & 0 & 0 & 0 & 1 & 0 & 0 \end{bmatrix}$$

$$x = \begin{bmatrix} 0 & 1 & 0 \end{bmatrix} \begin{bmatrix} 1 & 0 & 0 & 0 & 1 & 0 & 0 \\ 0 & 1 & 0 & 1 & 0 & 0 & 0 \\ 0 & 0 & 1 & 0 & 1 & 1 & 1 \end{bmatrix}$$
$$= \begin{bmatrix} 0 & 1 & 0 & 1 & 0 & 0 & 0 \end{bmatrix}$$

$$x = \begin{bmatrix} 0 & 0 & 1 \end{bmatrix} \begin{bmatrix} 1 & 0 & 0 & 0 & 1 & 0 & 0 \\ 0 & 1 & 0 & 1 & 0 & 0 & 0 \\ 0 & 0 & 1 & 0 & 1 & 1 & 1 \end{bmatrix}$$
$$= \begin{bmatrix} 0 & 0 & 1 & 0 & 1 & 1 & 1 \end{bmatrix}$$

$$x = \begin{bmatrix} 1 & 1 & 0 \end{bmatrix} \begin{bmatrix} 1 & 0 & 0 & 0 & 1 & 0 & 0 \\ 0 & 1 & 0 & 1 & 0 & 0 & 0 \\ 0 & 0 & 1 & 0 & 1 & 1 & 1 \end{bmatrix}$$
$$= \begin{bmatrix} 1 & 1 & 0 & 1 & 1 & 0 & 0 \end{bmatrix}$$

$$x = \begin{bmatrix} 1 & 0 & 1 \end{bmatrix} \begin{bmatrix} 1 & 0 & 0 & 0 & 1 & 0 & 0 \\ 0 & 1 & 0 & 1 & 0 & 0 & 0 \\ 0 & 0 & 1 & 0 & 1 & 1 & 1 \end{bmatrix}$$

$$= \begin{bmatrix}1 & 0 & 1 & 0 & 0 & 1 & 1\end{bmatrix}$$

$$x = \begin{bmatrix}0 & 1 & 1\end{bmatrix}\begin{bmatrix}1 & 0 & 0 & 0 & 1 & 0 & 0\\0 & 1 & 0 & 1 & 0 & 0 & 0\\0 & 0 & 1 & 0 & 1 & 1 & 1\end{bmatrix}$$

$$= \begin{bmatrix}0 & 1 & 1 & 1 & 1 & 1 & 1\end{bmatrix}$$

$$x = \begin{bmatrix}1 & 1 & 1\end{bmatrix}\begin{bmatrix}1 & 0 & 0 & 0 & 1 & 0 & 0\\0 & 1 & 0 & 1 & 0 & 0 & 0\\0 & 0 & 1 & 0 & 1 & 1 & 1\end{bmatrix}$$

$$= \begin{bmatrix}1 & 1 & 1 & 1 & 0 & 1 & 1\end{bmatrix}.$$

The set of codes words generated by this G are

(0 0 0 0 0 0 0), (1 0 0 0 1 0 0), (0 1 0 1 0 0 0), (0 0 1 0 1 1 1), (1 1 0 1 1 0 0), (1 0 1 0 0 1 1), (0 1 1 1 1 1 1) and (1 1 1 1 0 1 1).

The corresponding parity check matrix H obtained from this G is given by

$$H = \begin{bmatrix}0 & 1 & 0 & 1 & 0 & 0 & 0\\1 & 0 & 1 & 0 & 1 & 0 & 0\\0 & 0 & 1 & 0 & 0 & 1 & 0\\0 & 0 & 1 & 0 & 0 & 0 & 1\end{bmatrix}.$$

Now

$$GH^T = \begin{bmatrix}1 & 0 & 0 & 0 & 1 & 0 & 0\\0 & 1 & 0 & 1 & 0 & 0 & 0\\0 & 0 & 1 & 0 & 1 & 1 & 1\end{bmatrix}\begin{bmatrix}0 & 1 & 0 & 0\\1 & 0 & 0 & 0\\0 & 1 & 1 & 1\\1 & 0 & 0 & 0\\0 & 1 & 0 & 0\\0 & 0 & 1 & 0\\0 & 0 & 0 & 1\end{bmatrix}$$

$$= \begin{bmatrix} 0 & 0 & 0 & 0 \\ 0 & 0 & 0 & 0 \\ 0 & 0 & 0 & 0 \end{bmatrix}.$$

We recall just the definition of Hamming distance and Hamming weight between two vectors. This notion is applied to codes to find errors between the sent message and the received message. As finding error in the received message happens to be one of the difficult problems more so is the correction of errors and retrieving the correct message from the received message.

**DEFINITION 1.2.3:** *The Hamming distance d(x, y) between two vectors* $x = x_1\, x_2 \ldots x_n$ *and* $y = y_1\, y_2 \ldots y_n$ *in* $F_q^n$ *is the number of coordinates in which x and y differ. The Hamming weight* $\omega(x)$ *of a vector* $x = x_1\, x_2 \ldots x_n$ *in* $F_q^n$ *is the number of non zero co ordinates in* $x_i$. *In short* $\omega(x) = d(x, 0)$.

We just illustrate this by a simple example.
Suppose x = [1 0 1 1 1 1 0] and y ∈ [0 1 1 1 1 0 1 ] belong to $F_2^7$ then D(x, y) = (x ~ y) = (1 0 1 1 1 1 0) ~ (0 1 1 1 1 0 1) = (1~0, 0~1, 1~1, 1~1, 1~1, 1~0, 0~1) = (1 1 0 0 0 1 1) = 4. Now Hamming weight ω of x is ω(x) = d(x, 0) = 5 and ω(y) = d(y, 0) = 5.

**DEFINITION 1.2.4:** *Let C be any linear code then the minimum distance* $d_{min}$ *of a linear code C is given as*

$$d_{min} = \min_{\substack{u,v \in C \\ u \neq v}} d(u,v).$$

*For linear codes we have*

$$d(u, v) = d(u - v, 0) = \omega(u - v).$$

Thus it is easily seen minimum distance of C is equal to the least weight of all non zero code words. A general code C of length n with k message symbols is denoted by C(n, k) or by a binary (n, k) code. Thus a parity check code is a binary (n, n – 1) code and a repetition code is a binary (n, 1) code.

If $H = (A, I_{n-k})$ be a parity check matrix in the standard form then $G = (I_k, -A^T)$ is the canonical generator matrix of the linear (n, k) code.

The check equations $(A, I_{n-k})\, x^T = (0)$ yield

$$\begin{bmatrix} x_{k+1} \\ x_{k+2} \\ \vdots \\ x_n \end{bmatrix} = -A \begin{bmatrix} x_1 \\ x_2 \\ \vdots \\ x_k \end{bmatrix} = -A \begin{bmatrix} a_1 \\ a_2 \\ \vdots \\ a_k \end{bmatrix}.$$

Thus we obtain

$$\begin{bmatrix} x_1 \\ x_2 \\ \vdots \\ x_n \end{bmatrix} = \begin{bmatrix} I_k \\ -A \end{bmatrix} \begin{bmatrix} a_1 \\ a_2 \\ \vdots \\ a_k \end{bmatrix}.$$

We transpose and denote this equation as

$$(x_1\ x_2 \ldots x_n) = (a_1\ a_2 \ldots a_k)\,(I_k, -A^7)$$
$$= (a_1\ a_2 \ldots a_k)\, G.$$

We have just seen that minimum distance

$$d_{min} = \min_{\substack{u,v \in C \\ u \neq v}} d(u,v).$$

If d is the minimum distance of a linear code C then the linear code of length n, dimension k and minimum distance d is called an (n, k, d) code.

Now having sent a message or vector x and if y is the received message or vector a simple decoding rule is to find the code word closest to x with respect to Hamming distance, i.e., one chooses an error vector e with the least weight. The decoding method is called “nearest neighbour decoding” and amounts to comparing y with all $q^k$ code words and choosing the closest among them. The nearest neighbour decoding is the maximum likelihood decoding if the probability p for correct transmission is > ½.

Obviously before, this procedure is impossible for large k but with the advent of computers one can easily run a program in few seconds and arrive at the result.
We recall the definition of sphere of radius r. The set $S_r(x) = \{y \in F_q^n / d(x, y) \leq r\}$ is called the sphere of radius r about $x \in F_q^n$.

In decoding we distinguish between the detection and the correction of error. We can say a code can correct t errors and can detect t + s, s ≥ 0 errors, if the structure of the code makes it possible to correct up to t errors and to detect t + j, 0 < j ≤ s errors which occurred during transmission over a channel.

A mathematical criteria for this, given in the linear code is ; A linear code C with minimum distance $d_{min}$ can correct upto t errors and can detect t + j, 0 < j ≤ s, errors if and only if zt + s ≤ $d_{min}$ or equivalently we can say "A linear code C with minimum distance d can correct t errors if and only if $t = \left[\frac{(d-1)}{2}\right]$. The real problem of coding theory is not merely to minimize errors but to do so without reducing the transmission rate unnecessarily. Errors can be corrected by lengthening the code blocks, but this reduces the number of message symbols that can be sent per second. To maximize the transmission rate we want code blocks which are numerous enough to encode a given message alphabet, but at the same time no longer than is necessary to achieve a given Hamming distance. One of the main problems of coding theory is "Given block length n and Hamming distance d, find the maximum number, A(n, d) of binary blocks of length n which are at distances ≥ d from each other".

Let $u = (u_1, u_2, \ldots, u_n)$ and $v = (v_1, v_2, \ldots, v_n)$ be vectors in $F_q^n$ and let $u.v = u_1v_1 + u_2v_2 + \ldots + u_nv_n$ denote the dot product of u and v over $F_q^n$. If u.v = 0 then u and v are called orthogonal.

**DEFINITION 1.2.5:** *Let C be a linear (n, k) code over $F_q$. The dual(or orthogonal)code $C^{\perp} = \{u \mid u.v = 0$ for all $v \in C\}$, $u \in F_q^n$. If C is a k-dimensional subspace of the n-dimensional* *vector space $F_q^n$ the orthogonal complement is of dimension n – k and an (n, n – k) code. It can be shown that if the code C has a generator matrix G and parity check matrix H then $C^\perp$ has generator matrix H and parity check matrix G.*

Orthogonality of two codes can be expressed by $GH^T = (0)$.

***Example 1.1 .5:*** Let us consider the parity check matrix H of a (7, 3) code where

$$H = \begin{bmatrix} 1 & 0 & 0 & 1 & 0 & 0 & 0 \\ 0 & 0 & 1 & 0 & 1 & 0 & 0 \\ 1 & 1 & 0 & 0 & 0 & 1 & 0 \\ 1 & 0 & 1 & 0 & 0 & 0 & 1 \end{bmatrix}.$$

The code words got using H are as follows

$$\begin{matrix} 0\ 0\ 0\ 0\ 0\ 0\ 0 & \\ 1\ 0\ 0\ 1\ 0\ 1\ 1 & 0\ 1\ 1\ 0\ 1\ 1\ 1 \\ 0\ 1\ 0\ 0\ 0\ 1\ 0 & 1\ 0\ 1\ 1\ 1\ 1\ 0 \\ 0\ 0\ 1\ 0\ 1\ 0\ 1 & 1\ 1\ 1\ 1\ 1\ 0\ 0 \\ 1\ 1\ 0\ 1\ 0\ 0\ 1 & \end{matrix}.$$

Now for the orthogonal code the parity check matrix H of the code happens to be generator matrix,

$$G = \begin{bmatrix} 1 & 0 & 0 & 1 & 0 & 0 & 0 \\ 0 & 0 & 1 & 0 & 1 & 0 & 0 \\ 1 & 1 & 0 & 0 & 0 & 1 & 0 \\ 1 & 0 & 1 & 0 & 0 & 0 & 1 \end{bmatrix}.$$

$$x = \begin{bmatrix} 0 \\ 0 \\ 0 \\ 0 \end{bmatrix}^T \begin{bmatrix} 1 & 0 & 0 & 1 & 0 & 0 & 0 \\ 0 & 0 & 1 & 0 & 1 & 0 & 0 \\ 1 & 1 & 0 & 0 & 0 & 1 & 0 \\ 1 & 0 & 1 & 0 & 0 & 0 & 1 \end{bmatrix} = [0\ 0\ 0\ 0\ 0\ 0\ 0].$$

$$x = \begin{bmatrix} 1 \\ 0 \\ 0 \\ 0 \end{bmatrix}^{T} \begin{bmatrix} 1 & 0 & 0 & 1 & 0 & 0 & 0 \\ 0 & 0 & 1 & 0 & 1 & 0 & 0 \\ 1 & 1 & 0 & 0 & 0 & 1 & 0 \\ 1 & 0 & 1 & 0 & 0 & 0 & 1 \end{bmatrix} = [1\ 0\ 0\ 1\ 0\ 0\ 0].$$

$$x = \begin{bmatrix} 0 & 1 & 0 & 0 \end{bmatrix} \begin{bmatrix} 1 & 0 & 0 & 1 & 0 & 0 & 0 \\ 0 & 0 & 1 & 0 & 1 & 0 & 0 \\ 1 & 1 & 0 & 0 & 0 & 1 & 0 \\ 1 & 0 & 1 & 0 & 0 & 0 & 1 \end{bmatrix} = [0\ 0\ 1\ 0\ 1\ 0\ 0]$$

$$x = \begin{bmatrix} 0 & 0 & 1 & 0 \end{bmatrix} \begin{bmatrix} 1 & 0 & 0 & 1 & 0 & 0 & 0 \\ 0 & 0 & 1 & 0 & 1 & 0 & 0 \\ 1 & 1 & 0 & 0 & 0 & 1 & 0 \\ 1 & 0 & 1 & 0 & 0 & 0 & 1 \end{bmatrix} = [1\ 1\ 0\ 0\ 0\ 1\ 0]$$

$$x = \begin{bmatrix} 0 & 0 & 0 & 1 \end{bmatrix} \begin{bmatrix} 1 & 0 & 0 & 1 & 0 & 0 & 0 \\ 0 & 0 & 1 & 0 & 1 & 0 & 0 \\ 1 & 1 & 0 & 0 & 0 & 1 & 0 \\ 1 & 0 & 1 & 0 & 0 & 0 & 1 \end{bmatrix} = [1\ 0\ 1\ 0\ 0\ 0\ 1]$$

$$x = \begin{bmatrix} 1 & 1 & 0 & 0 \end{bmatrix} \begin{bmatrix} 1 & 0 & 0 & 1 & 0 & 0 & 0 \\ 0 & 0 & 1 & 0 & 1 & 0 & 0 \\ 1 & 1 & 0 & 0 & 0 & 1 & 0 \\ 1 & 0 & 1 & 0 & 0 & 0 & 1 \end{bmatrix} = [1\ 0\ 1\ 1\ 1\ 0\ 0]$$

$$x = \begin{bmatrix} 1 & 0 & 1 & 0 \end{bmatrix} \begin{bmatrix} 1 & 0 & 0 & 1 & 0 & 0 & 0 \\ 0 & 0 & 1 & 0 & 1 & 0 & 0 \\ 1 & 1 & 0 & 0 & 0 & 1 & 0 \\ 1 & 0 & 1 & 0 & 0 & 0 & 1 \end{bmatrix} = [0\ 1\ 0\ 1\ 0\ 1\ 0]$$

$$x = \begin{bmatrix} 1 & 0 & 0 & 1 \end{bmatrix} \begin{bmatrix} 1 & 0 & 0 & 1 & 0 & 0 & 0 \\ 0 & 0 & 1 & 0 & 1 & 0 & 0 \\ 1 & 1 & 0 & 0 & 0 & 1 & 0 \\ 1 & 0 & 1 & 0 & 0 & 0 & 1 \end{bmatrix} = [0\ 0\ 1\ 1\ 0\ 0\ 1]$$

$$x = \begin{bmatrix} 0 & 1 & 1 & 0 \end{bmatrix} \begin{bmatrix} 1 & 0 & 0 & 1 & 0 & 0 & 0 \\ 0 & 0 & 1 & 0 & 1 & 0 & 0 \\ 1 & 1 & 0 & 0 & 0 & 1 & 0 \\ 1 & 0 & 1 & 0 & 0 & 0 & 1 \end{bmatrix} = [1\ 1\ 1\ 0\ 1\ 1\ 0]$$

$$x = \begin{bmatrix} 0 & 1 & 0 & 1 \end{bmatrix} \begin{bmatrix} 1 & 0 & 0 & 1 & 0 & 0 & 0 \\ 0 & 0 & 1 & 0 & 1 & 0 & 0 \\ 1 & 1 & 0 & 0 & 0 & 1 & 0 \\ 1 & 0 & 1 & 0 & 0 & 0 & 1 \end{bmatrix} = [1\ 0\ 0\ 0\ 1\ 0\ 1]$$

$$x = \begin{bmatrix} 0 & 0 & 1 & 1 \end{bmatrix} \begin{bmatrix} 1 & 0 & 0 & 1 & 0 & 0 & 0 \\ 0 & 0 & 1 & 0 & 1 & 0 & 0 \\ 1 & 1 & 0 & 0 & 0 & 1 & 0 \\ 1 & 0 & 1 & 0 & 0 & 0 & 1 \end{bmatrix} = [0\ 1\ 1\ 0\ 0\ 1\ 1]$$

$$x = \begin{bmatrix} 1 & 1 & 1 & 0 \end{bmatrix} \begin{bmatrix} 1 & 0 & 0 & 1 & 0 & 0 & 0 \\ 0 & 0 & 1 & 0 & 1 & 0 & 0 \\ 1 & 1 & 0 & 0 & 0 & 1 & 0 \\ 1 & 0 & 1 & 0 & 0 & 0 & 1 \end{bmatrix} = [0\ 1\ 1\ 1\ 1\ 1\ 0]$$

$$x = \begin{bmatrix} 1 & 1 & 0 & 1 \end{bmatrix} \begin{bmatrix} 1 & 0 & 0 & 1 & 0 & 0 & 0 \\ 0 & 0 & 1 & 0 & 1 & 0 & 0 \\ 1 & 1 & 0 & 0 & 0 & 1 & 0 \\ 1 & 0 & 1 & 0 & 0 & 0 & 1 \end{bmatrix} = [0\ 0\ 0\ 1\ 1\ 0\ 1]$$

$$x = \begin{bmatrix} 1 & 0 & 1 & 1 \end{bmatrix} \begin{bmatrix} 1 & 0 & 0 & 1 & 0 & 0 & 0 \\ 0 & 0 & 1 & 0 & 1 & 0 & 0 \\ 1 & 1 & 0 & 0 & 0 & 1 & 0 \\ 1 & 0 & 1 & 0 & 0 & 0 & 1 \end{bmatrix} = [1\ 1\ 1\ 1\ 0\ 1\ 1]$$

$$x = \begin{bmatrix} 0 & 1 & 1 & 1 \end{bmatrix} \begin{bmatrix} 1 & 0 & 0 & 1 & 0 & 0 & 0 \\ 0 & 0 & 1 & 0 & 1 & 0 & 0 \\ 1 & 1 & 0 & 0 & 0 & 1 & 0 \\ 1 & 0 & 1 & 0 & 0 & 0 & 1 \end{bmatrix} = [0\ 1\ 0\ 0\ 1\ 1\ 1]$$

$$x = \begin{bmatrix} 1 & 1 & 1 & 1 \end{bmatrix} \begin{bmatrix} 1 & 0 & 0 & 1 & 0 & 0 & 0 \\ 0 & 0 & 1 & 0 & 1 & 0 & 0 \\ 1 & 1 & 0 & 0 & 0 & 1 & 0 \\ 1 & 0 & 1 & 0 & 0 & 0 & 1 \end{bmatrix} = [1\ 1\ 0\ 1\ 1\ 1\ 1]\ .$$

The code words of C(7, 4) i.e., the orthogonal code of C(7, 3) are

{(0 0 0 0 0 0 0), (1 0 0 1 0 0 0), (0 0 1 0 1 0 0), (1 1 0 0 0 1 0), (1 0 1 0 0 0 1), (1 0 1 1 1 0 0), (0 1 0 1 0 1 0), (0 0 1 1 0 0 1), (1 1 1 0 1 1 0), (1 0 0 0 1 0 1), (0 1 1 0 0 1 1), (0 1 1 1 1 1 0), (0 0 0 1 1 0 1), (1 1 1 1 0 1 1), (0 1 0 0 1 1 1), (1 1 0 1 1 1 1)}

Thus we have found the orthogonal code for the given code. Now we just recall the definition of the cosets of a code C.

**DEFINITION 1.2.6:** *For* $a \in F_q^n$ *we have* $a + C = \{a + x \,/\, x \in C\}$*. Clearly each coset contains* $q^k$ *vectors. There is a partition of* $F_q^n$ *of the form* $F_q^n = C \cup \{a^{(1)} + C\} \cup \{a^{(2)} + C\} \cup \ldots \cup \{a^t + C\}$ *for* $t = q^{n-k} - 1$*. If y is a received vector then y must be an element of one of these cosets say* $a^i + C$*. If the code word* $x^{(1)}$ *has been transmitted then the error vector*

$$e = y - x^{(1)} \in a^{(i)} + C - x^{(1)} = a^{(i)} + C.$$

*Now we give the decoding rule which is as follows:*

*If a vector y is received then the possible error vectors e are the vectors in the coset containing y. The most likely error is the vector* $\overline{e}$ *with minimum weight in the coset of y. Thus y is decoded as* $\overline{x} = y - \overline{e}$ *. [18-21]*

*Now we show how to find the coset of y and describe the above method. The vector of minimum weight in a coset is called the coset leader.*

*If there are several such vectors then we arbitrarily choose one of them as coset leader. Let* $a^{(1)}, a^{(2)}, \ldots, a^{(t)}$ *be the coset leaders. We first establish the following table*

| $x^{(1)} = 0$ | $x^{(2)} = 0$ | $\ldots$ | $x^{(q^k)}$ | *code words in C* |
|---|---|---|---|---|
| $a^{(1)} + x^{(1)}$ | $a^{(1)} + x^{(2)}$ | $\ldots$ | $a^{(1)} + x^{(q^k)}$ | |
| $\vdots$ | $\vdots$ | | $\vdots$ | *other cosets* |
| $a^{(t)} + x^{(1)}$ | $a^{(t)} + x^{(2)}$ | $\ldots$ | $a^{(t)} + x^{(q^k)}$ | |
| *coset leaders* | | | | |

*If a vector y is received then we have to find y in the table. Let y =* $a^{(i)} + x^{(j)}$*; then the decoder decides that the error* $\overline{e}$ *is the coset leader* $a^{(i)}$*. Thus y is decoded as the code word* $\overline{x} = y - \overline{e} = x^{(j)}$ *. The code word* $\overline{x}$ *occurs as the first element in the column of y. The coset of y can be found by evaluating the so called syndrome.*

*Let H be parity check matrix of a linear (n, k) code. Then the vector* $S(y) = Hy^T$ *of length n–k is called syndrome of y. Clearly S(y) = (0) if and only if* $y \in C$*.*

$S(y^{(1)}) = S(y^{(2)})$ *if and only if* $y^{(1)} + C = y^{(2)} + C$ *.*

*We have the decoding algorithm as follows:*

*If* $y \in F_q^n$ *is a received vector find S(y), and the coset leader* $\overline{e}$ *with syndrome S(y). Then the most likely transmitted code word is* $\overline{x} = y - \overline{e}$ *we have* $d(\overline{x}, y) .= \min\{d(x, y)/x \in C\}$*.*

We illustrate this by the following example:

***Example 1.2.6:*** Let C be a (5, 3) code where the parity check matrix H is given by

$$H = \begin{bmatrix} 1 & 0 & 1 & 1 & 0 \\ 1 & 1 & 0 & 0 & 1 \end{bmatrix}$$

and

$$G = \begin{bmatrix} 1 & 0 & 0 & 1 & 1 \\ 0 & 1 & 0 & 0 & 1 \\ 0 & 0 & 1 & 1 & 0 \end{bmatrix}.$$

The code words of C are

{(0 0 0  0 0), (1 0 0 1 1), (0 1 0 0 1), (0 0 1 1 0), (1 1 0 1 0), (1 0 1 0 1), (0 1 1 1 1), (1 1 1 0 0)}.

The corresponding coset table is

| Message | 000 | 100 | 010 | 001 | 110 | 101 | 011 | 111 |
|---|---|---|---|---|---|---|---|---|
| code words | 00000 | 10011 | 01001 | 00110 | 11010 | 10101 | 01111 | 11100 |
| other cosets | 10000 | 00011 | 11001 | 10110 | 01010 | 00101 | 11111 | 01100 |
| | 01000 | 11011 | 00001 | 01110 | 10010 | 11101 | 00111 | 10100 |
| | 00100 | 10111 | 01101 | 00010 | 11110 | 10001 | 01011 | 11000 |

coset leaders

If y = (1  1  1  1  0) is received, then y is found in the coset with the coset leader (0 0 1 0 0)
y + (0  0  1  0  0) = (1  1  1  1  0) + (0  0  1  0  0 ) = (1  1  0  1  0) is the corresponding message.

Now with the advent of computers it is easy to find the real message or the sent word by using this decoding algorithm.

A binary code $C_m$ of length $n = 2^m - 1$, $m \geq 2$ with $m \times 2^m - 1$ parity check matrix H whose columns consists of all non zero binary vectors of length m is called a binary Hamming code.

We give example of them.

***Example 1.2.7:*** Let

$$H = \begin{bmatrix} 1 & 0 & 1 & 1 & 1 & 1 & 0 & 0 & 1 & 0 & 1 & 1 & 0 & 0 & 0 \\ 1 & 1 & 0 & 1 & 1 & 1 & 1 & 0 & 0 & 1 & 0 & 0 & 1 & 0 & 0 \\ 1 & 1 & 1 & 0 & 1 & 0 & 1 & 1 & 0 & 0 & 1 & 0 & 0 & 1 & 0 \\ 1 & 1 & 1 & 1 & 0 & 0 & 0 & 1 & 1 & 1 & 0 & 0 & 0 & 0 & 1 \end{bmatrix}$$

which gives a $C_4(15, 11, 4)$ Hamming code.

Cyclic codes are codes which have been studied extensively.

Let us consider the vector space $F_q^n$ over $F_q$. The mapping

$$Z: F_q^n \to F_q^n$$

where Z is a linear mapping called a "cyclic shift" if $Z(a_0, a_1, \ldots, a_{n-1}) = (a_{n-1}, a_0, \ldots, a_{n-2})$

$A = (F_q[x], +, ., .)$ is a linear algebra in a vector space over $F_q$. We define a subspace $V_n$ of this vector space by

$$\begin{aligned} V_n &= \{v \in F_q[x] \,/\, \text{degree } v < n\} \\ &= \{v_0 + v_1x + v_2x^2 + \ldots + v_{n-1}x^{n-1} \,/\, v_i \in F_q;\ 0 \leq i \leq n-1\}. \end{aligned}$$

We see that $V_n \cong F_q^n$ as both are vector spaces defined over the same field $F_q$. Let Γ be an isomorphism

$\Gamma(v_0, v_1, \ldots, v_{n-1}) \to \{v_0 + v_1x + v_2x^2 + \ldots + v_{n-1}x^{n-1}\}$.

w: $F_q^n$ ℧ $F_q[x] / x^n - 1$

i.e., $w(v_0, v_1, \ldots, v_{n-1}) = v_0 + v_1x + \ldots + v_{n-1}x^{n-1}$.

Now we proceed onto define the notion of a cyclic code.

**DEFINITION 1.2.7:** *A k-dimensional subspace C of $F_q^n$ is called a cyclic code if Z(v) ∈ C for all v ∈ C that is $v = v_0, v_1, \ldots, v_{n-1} \in C$ implies $(v_{n-1}, v_0, \ldots, v_{n-2}) \in C$ for $v \in F_q^n$.*

We just give an example of a cyclic code.

***Example 1.2.8:*** Let $C \subseteq F_2^7$ be defined by the generator matrix

$$G = \begin{bmatrix} 1 & 1 & 1 & 0 & 1 & 0 & 0 \\ 0 & 1 & 1 & 1 & 0 & 1 & 0 \\ 0 & 0 & 1 & 1 & 1 & 0 & 1 \end{bmatrix} = \begin{bmatrix} g^{(1)} \\ g^{(2)} \\ g^{(3)} \end{bmatrix}.$$

The code words generated by G are {(0 0 0 0 0 0 0), (1 1 1 0 1 0 0), (0 1 1 1 0 1 0), (0 0 1 1 1 0 1), (1 0 0 1 1 1 0), (1 1 0 1 0 0 1), (0 1 0 0 1 1 1), (1 0 1 0 0 1 1)}.

Clearly one can check the collection of all code words in C satisfies the rule if $(a_0 \ldots a_5) \in C$ then $(a_5\ a_0 \ldots a_4) \in C$ i.e., the codes are cyclic. Thus we get a cyclic code.

Now we see how the code words of the Hamming codes looks like.

***Example 1.2.9:*** Let

$$H = \begin{bmatrix} 1 & 0 & 0 & 1 & 1 & 0 & 1 \\ 0 & 1 & 0 & 1 & 0 & 1 & 1 \\ 0 & 0 & 1 & 0 & 1 & 1 & 1 \end{bmatrix}$$

be the parity check matrix of the Hamming (7, 4) code.

Now we can obtain the elements of a Hamming(7,4) code.

We proceed on to define parity check matrix of a cyclic code given by a polynomial matrix equation given by defining the generator polynomial and the parity check polynomial.

**DEFINITION 1.2.8:** *A linear code C in $V_n = \{v_0 + v_1x + \ldots + v_{n-1}x^{n-1} \mid v_i \in F_q, 0 \leq i \leq n-1\}$ is cyclic if and only if C is a principal ideal generated by $g \in C$.*

*The polynomial g in C can be assumed to be monic. Suppose in addition that $g / x^n - 1$ then g is uniquely determined and is called the generator polynomial of C. The elements of C are called code words, code polynomials or code vectors.*

*Let $g = g_0 + g_1x + \ldots + g_mx^m \in V_n$, $g / x^n - 1$ and $\deg g = m < n$. Let C be a linear (n, k) code, with $k = n - m$ defined by the generator matrix,*

$$G = \begin{bmatrix} g_0 & g_1 & \cdots & g_m & 0 & \cdots & 0 \\ 0 & g_0 & \cdots & g_{m-1} & g_m & \cdots & 0 \\ \vdots & \vdots & & & & & \\ 0 & 0 & & g_0 & g_1 & & g_m \end{bmatrix} = \begin{bmatrix} g \\ xg \\ \\ x^{k-1}g \end{bmatrix}.$$

*Then C is cyclic. The rows of G are linearly independent and rank G = k, the dimension of C.*

***Example 1.2.10:*** Let $g = x^3 + x^2 + 1$ be the generator polynomial having a generator matrix of the cyclic(7,4) code with generator matrix

$$G = \begin{bmatrix} 1 & 0 & 1 & 1 & 0 & 0 & 0 \\ 0 & 1 & 0 & 1 & 1 & 0 & 0 \\ 0 & 0 & 1 & 0 & 1 & 1 & 0 \\ 0 & 0 & 0 & 1 & 0 & 1 & 1 \end{bmatrix}.$$

The codes words associated with the generator matrix is

0000000, 1011000, 0101100, 0010110, 0001011, 1110100, 1001110, 1010011, 0111010, 0100111, 0011101, 1100010, 1111111, 1000101, 0110001, 1101001.

The parity check polynomial is defined to be

$$h = \frac{x^7 - 1}{g}$$

$$h = \frac{x^7 - 1}{x^3 + x^2 + 1} = x^4 + x^3 + x^2 + 1.$$

$$\text{If } \frac{x^n - 1}{g} = h_0 + h_1x + \ldots + h_kx^k.$$

the parity check matrix H related with the generator polynomial g is given by

$$H = \begin{bmatrix} 0 & \dots & 0 & h_k & \dots & h_1 & h_0 \\ 0 & \dots & h_k & h_{k-1} & & h_0 & 0 \\ \vdots & \vdots & \vdots & \vdots & & & \vdots \\ h_k & \dots & h_1 & h_0 & & \dots & 0 \end{bmatrix}.$$

For the generator polynomial g = $x^3 + x^2 + 1$ the parity check matrix

$$H = \begin{bmatrix} 0 & 0 & 1 & 1 & 1 & 0 & 1 \\ 0 & 1 & 1 & 1 & 0 & 1 & 0 \\ 1 & 1 & 1 & 0 & 1 & 0 & 0 \end{bmatrix}$$

where the parity check polynomial is given by $x^4 + x^3 + x^2 + 1 = \frac{x^7 - 1}{x^3 + x^2 + 1}$. It is left for the reader to verify that the parity check matrix gives the same set of cyclic codes.

We now proceed on to give yet another new method of decoding procedure using the method of best approximations.

We just recall this definition given by [9, 19-21]. We just give the basic concepts needed to define this notion. We know that $F_q^n$ is a finite dimensional vector space over $F_q$. If we take $Z_2$ = (0, 1) the finite field of characteristic two. $Z_2^5 = Z_2 \times Z_2 \times Z_2 \times Z_2 \times Z_2$ is a 5 dimensional vector space over $Z_2$. Infact {(1 0 0 0 0), (0 1 0 0 0), (0 0 1 0 0), (0 0 0 1 0), (0 0 0 0 1)} is a basis of $Z_2^5$. $Z_2^5$ has only $2^5 = 32$ elements in it. Let F be a field of real numbers and V a vector space over F. An inner product on V is a function which assigns to each ordered pair of vectors α, β in V a scalar ⟨α /β ⟩ in F in such a way that for all α, β, γ in V and for all scalars c in F.

(a) $\langle \alpha + \beta / \gamma \rangle = \langle \alpha/\gamma \rangle + \langle \beta/\gamma \rangle$

(b) $\langle c\alpha / \beta\rangle = c\langle\alpha/\beta\rangle$
(c) $\langle\beta/\alpha\rangle = \langle\alpha/\beta\rangle$
(d) $\langle\alpha/\alpha\rangle > 0$ if $\alpha \neq 0$.

On V there is an inner product which we call the standard inner product. Let $\alpha = (x_1, x_2, \ldots, x_n)$ and $\beta = (y_1, y_2, \ldots, y_n)$

$$\langle\alpha / \beta\rangle = \sum_i x_i y_i \ .$$

This is called as the standard inner product. $\langle\alpha/\alpha\rangle$ is defined as norm and it is denoted by $||\alpha||$. We have the Gram-Schmidt orthogonalization process which states that if V is a vector space endowed with an inner product and if $\beta_1, \beta_2, \ldots, \beta_n$ be any set of linearly independent vectors in V; then one may construct a set of orthogonal vectors $\alpha_1, \alpha_2, \ldots, \alpha_n$ in V such that for each $k = 1, 2, \ldots, n$ the set $\{\alpha_1, \ldots, \alpha_k\}$ is a basis for the subspace spanned by $\beta_1, \beta_2, \ldots, \beta_k$ where $\alpha_1 = \beta_1$.

$$\alpha_2 = \beta_2 - \frac{\langle\beta_1 / \alpha_1\rangle}{||\alpha_1||^2}\ \alpha_1$$

$$\alpha_3 = \beta_3 - \frac{\langle\beta_3 / \alpha_1\rangle}{||\alpha_1||^2}\ \alpha_1 - \frac{\langle\beta_3 / \alpha_2\rangle}{||\alpha_2||^2}\ \alpha_2$$

and so on.

Further it is left as an exercise for the reader to verify that if a vector $\beta$ is a linear combination of an orthogonal sequence of non-zero vectors $\alpha_1, \ldots, \alpha_m$, then $\beta$ is the particular linear combination, i.e.,

$$\beta = \sum_{k=1}^{m} \frac{\langle\beta / \alpha_k\rangle}{||\alpha_k||^2}\ \alpha_k \ .$$

In fact this property that will be made use of in the best approximations.

We just proceed on to give an example.

***Example 1.2.11:*** Let us consider the set of vectors $\beta_1 = (2, 0, 3)$, $\beta_2 = (-1, 0, 5)$ and $\beta_3 = (1, 9, 2)$ in the space $R^3$ equipped with the standard inner product.
Define $\alpha_1 = (2, 0, 3)$

$$\alpha_2 = (-1, 0, 5) - \frac{\langle(-1, 0, 5)/(2, 0, 3)\rangle}{13}(2, 0, 3)$$

$$= (-1, 0, 5) - \frac{13}{13}(2, 0, 3) = (-3, 0, 2)$$

$$\alpha_3 = (1,9,2) - \frac{\langle(-1, 9, 2)/(2, 0, 3)\rangle}{13}(2, 0, 3)$$

$$- \frac{\langle(1, 9, 2)/(-3, 0, 2)\rangle}{13}(-3,0,2)$$

$$= (1,9,2) - \frac{8}{13}(2,0,3) - \frac{1}{13}(-3,0,2)$$

$$= (1,9,2) - \left(\frac{16}{13}, 0, \frac{24}{13}\right) - \left(\frac{3}{13}, 0, \frac{2}{13}\right)$$

$$= (1,9,2) - \left\{\left(\frac{16-3}{13}, 0, \frac{24+2}{13}\right)\right\}$$

$$= (1, 9, 2) - (1, 0, 2)$$

$$= (0, 9, 0).$$

Clearly the set {(2, 0, 3), (–3, 0, 2), (0, 9, 0)} is an orthogonal set of vectors.

Now we proceed on to define the notion of a best approximation to a vector $\beta$ in V by vectors of a subspace W where $\beta \notin W$. Suppose W is a subspace of an inner product space V and let $\beta$ be an arbitrary vector in V. The problem is to find a best possible approximation to $\beta$ by vectors in W. This means we want to find a vector $\alpha$ for which $||\beta - \alpha||$ is as small as possible subject to the restriction that $\alpha$ should belong to W. To be precisely in mathematical terms: A best approximation to $\beta$ by vectors in W is a vector $\alpha$ in W such that $||\beta - \alpha|| \leq ||\beta - \gamma||$ for every vector $\gamma$ in W ; W a subspace of V.

By looking at this problem in $R^2$ or in $R^3$ one sees intuitively that a best approximation to $\beta$ by vectors in W ought to be a vector $\alpha$ in W such that $\beta - \alpha$ is perpendicular (orthogonal) to W and that there ought to be exactly one such $\alpha$. These intuitive ideas are correct for some finite dimensional subspaces, but not for all infinite dimensional subspaces.

We just enumerate some of the properties related with best approximation.

Let W be a subspace of an inner product space V and let $\beta$ be a vector in V.

(i) The vector $\alpha$ in W is a best approximation to $\beta$ by vectors in W if and only if $\beta - \alpha$ is orthogonal to every vector in W.

(ii) If a best approximation to $\beta$ by vectors in W exists, it is unique.

(iii) If W is finite-dimensional and $\{\alpha_1, \alpha_2, \ldots, \alpha_n\}$ is any orthonormal basis for W, then the vector $\alpha = \sum_k \frac{\langle \beta / \alpha_k \rangle}{\|\alpha_k\|^2} \alpha_k$, where $\alpha$ is the (unique)best approximation to $\beta$ by vectors in W.

Now this notion of best approximation for the first time is used in coding theory to find the best approximated sent code after receiving a message which is not in the set of codes used. Further we use for coding theory only finite fields $F_q$. i.e., $|F_q| < \infty$ . If C is a code of length n; C is a vector space over $F_q$ and C $\cong F_q^k \subseteq F_q^n$, k the number of message symbols in the code, i.e., C is a C(n, k) code. While defining the notion of inner product on vector spaces over finite fields we see all axiom of inner product defined over fields as reals or complex in general is not true. The main property which is not true is if $0 \neq x \in V$; the inner product of x with itself i.e., $\langle x / x \rangle = \langle x, x \rangle \neq 0$ if $x \neq 0$ is not true i.e., $\langle x / x \rangle = 0$ does not imply $x = 0$.

To overcome this problem we define for the first time the new notion of pseudo inner product in case of vector spaces defined over finite characteristic fields [9, 20-1].

**DEFINITION 1.2.9:** *Let V be a vector space over a finite field $F_p$ of characteristic p, p a prime. Then the pseudo inner product on*

*V is a map* $\langle , \rangle_p : V \times V \to F_p$ *satisfying the following conditions:*

1. $\langle x, x\rangle_p \geq 0$ *for all* $x \in V$ *and* $\langle x, x\rangle_p = 0$ *does not in general imply* $x = 0$.
2. $\langle x, y\rangle_p = \langle y, x\rangle_p$ *for all* $x, y \in V$.
3. $\langle x + y, z\rangle_p = \langle x, z\rangle_p + \langle y, z\rangle_p$ *for all* $x, y, z \in V$.
4. $\langle x, y + z\rangle_p = \langle x, y\rangle_p + \langle x, z\rangle_p$ *for all* $x, y, z \in V$.
5. $\langle \alpha.x, y\rangle_p = \alpha \langle x, y\rangle_p$ *and*
6. $\langle x, \beta.y\rangle_p = \beta\langle x, y\rangle_p$ *for all* $x, y, \in V$ *and* $\alpha, \beta \in F_p$.

*Let V be a vector space over a field* $F_p$ *of characteristic p, p is a prime; then V is said to be a pseudo inner product space if there is a pseudo inner product* $\langle , \rangle_p$ *defined on V. We denote the pseudo inner product space by* $(V, \langle , \rangle_p)$.

Now using this pseudo inner product space $(V, \langle , \rangle_p)$ we proceed on to define pseudo-best approximation.

**DEFINITION 1.2.10:** *Let V be a vector space defined over the finite field* $F_p$ *(or* $Z_p$*). Let W be a subspace of V. For* $\beta \in V$ *and for a set of basis* $\{\alpha_1, \ldots, \alpha_k\}$ *of the subspace W the pseudo best approximation to* $\beta$*, if it exists is given by* $\sum_{i=1}^{k} \langle \beta, \alpha_i \rangle_p \alpha_i$. *If* $\sum_{i=1}^{k} \langle \beta, \alpha_i \rangle_p \alpha_i = 0$*, then we say the pseudo best approximation does not exist for this set of basis* $\{\alpha_1, \alpha_2, \ldots, \alpha_k\}$*. In this case we choose another set of basis for W say* $\{\gamma_1, \gamma_2, \ldots, \gamma_k\}$ *and calculate* $\sum_{i=1}^{k} \langle \beta, \gamma_i \rangle_p \gamma_i$ *and* $\sum_{i=1}^{k} \langle \beta, \gamma_i \rangle_p \gamma_i$ *is called a pseudo best approximation to* $\beta$.

*Note:* We need to see the difference even in defining our pseudo best approximation with the definition of the best approximation. Secondly as we aim to use it in coding theory and most of our linear codes take only their values from the field of characteristic two we do not need $\langle x, x\rangle$ or the norm to

be divided by the pseudo inner product in the summation of finding the pseudo best approximation.

Now first we illustrate the pseudo inner product by an example.

***Example 1.2.12:*** Let V = $Z_2 \times Z_2 \times Z_2 \times Z_2$ be a vector space over $Z_2$. Define $\langle , \rangle_p$ to be the standard pseudo inner product on V; so if x = (1 0 1 1) and y = (1 1 1 1) are in V then the pseudo inner product of

$$\langle x, y\rangle_p = \langle(1\ 0\ 1\ 1\ ), (1\ 1\ 1\ 1)\rangle_p = 1 + 0 + 1 + 1 = 1.$$

Now consider

$$\langle x, x\rangle_p = \langle(1\ 0\ 1\ 1), (1\ 0\ 1\ 1)\rangle_p = 1 + 0 + 1 + 1 \neq 0$$

but

$$\langle y, y\rangle_p = \langle(1\ 1\ 1\ 1), (1\ 1\ 1\ 1)\rangle_p = 1 + 1 + 1 + 1 = 0.$$

We see clearly y ≠ 0, yet the pseudo inner product is zero.

Now having seen an example of the pseudo inner product we proceed on to illustrate by an example the notion of pseudo best approximation.

***Example 1.2.13:*** Let

$$V = Z_2^8 = \underbrace{Z_2 \times Z_2 \times \ldots \times Z_2}_{8 \text{ times}}$$

be a vector space over $Z_2$. Now
W = {0 0 0 0 0 0 0 0), (1 0 0 0 1 0 11), (0 1 0 0 1 1 0 0), (0 0 1 0 0 1 1 1), (0 0 0 1 1 1 0 1), (1 1 0 0 0 0 1 0), (0 1 1 0 1 1 1 0), (0 0 1 1 1 0 1 0), (0 1 0 1 0 1 0 0), (1 0 1 0 1 1 0 0), (1 0 0 1 0 1 1 0), (1 1 1 0 0 1 0 1), (0 1 1 1 0 0 1 1), (1 1 0 1 1 1 1 1), (1 0 1 1 0 0 0 1), (1 1 1 1 1 0 0 0)}

be a subspace of V. Choose a basis of W as B = {$\alpha_1, \alpha_2, \alpha_3, \alpha_4$} where

$$\alpha_1 = (0\ 1\ 0\ 0\ 1\ 0\ 0\ 1),$$
$$\alpha_2 = (1\ 1\ 0\ 0\ 0\ 0\ 1\ 0),$$
$$\alpha_3 = (1\ 1\ 1\ 0\ 0\ 1\ 0\ 1)$$

and

$$\alpha_4 = (1\ 1\ 1\ 1\ 1\ 0\ 0\ 0).$$

Suppose $\beta = (1\ 1\ 1\ 1\ 1\ 1\ 1\ 1)$ is a vector in V using pseudo best approximations find a vector in W close to $\beta$. This is given by $\alpha$ relative to the basis B of W where

$$\alpha = \sum_{k=-1}^{4} \langle \beta, \alpha_k \rangle_p \alpha_k$$

$= \langle(1\ 1\ 1\ 1\ 1\ 1\ 1\ 1), (0\ 1\ 0\ 0\ 1\ 0\ 0\ 1)\rangle_p\ \alpha_1 +$
$\langle(1\ 1\ 1\ 1\ 1\ 1\ 1\ 1), (1\ 1\ 0\ 0\ 0\ 0\ 1\ 0)\rangle_p\ \alpha_2 +$
$\langle(1\ 1\ 1\ 1\ 1\ 1\ 1\ 1), (1\ 1\ 1\ 0\ 0\ 1\ 0\ 1)\rangle_p\ \alpha_3 +$
$\langle(1\ 1\ 1\ 1\ 1\ 1\ 1\ 1), (1\ 1\ 1\ 1\ 1\ 0\ 0\ 0)\rangle_p\ \alpha_4.$

$= 1.\alpha_1 + 1.\alpha_2 + 1.\alpha_3 + 1.\alpha_4.$

$= (0\ 1\ 0\ 0\ 1\ 0\ 0\ 1) + (1\ 1\ 0\ 0\ 0\ 0\ 1\ 0) + (1\ 1\ 1\ 0\ 0\ 1\ 0\ 1) + (1\ 1\ 1\ 1\ 1\ 0\ 0\ 0)$

$= (1\ 0\ 0\ 1\ 0\ 1\ 1\ 0) \in W.$

Now having illustrated how the pseudo best approximation of a vector $\beta$ in V relative to a subspace W of V is determined, now we illustrate how the approximately the nearest code word is obtained.

***Example 1.2.14:*** Let C = C(4, 2) be a code obtained from the parity check matrix

$$H = \begin{bmatrix} 1 & 0 & 1 & 0 \\ 1 & 1 & 0 & 1 \end{bmatrix}.$$

The message symbols associated with the code C are {(0, 0), (1, 0), (1, 0), (1, 1)}. The code words associated with H are C = {(0 0 0 0), (1 0 1 1), (0 1 0 1), (1 1 1 0)}. The chosen basis for C is B = $\{\alpha_1, \alpha_2\}$ where $\alpha_1 = (0\ 1\ 0\ 1)$ and $\alpha_2 = (1\ 0\ 1\ 1)$. Suppose the received message is $\beta = (1\ 1\ 1\ 1)$, consider $H\beta^T = (0\ 1) \neq (0)$ so $\beta \notin C$. Let $\alpha$ be the pseudo best approximation to $\beta$ relative to the basis B given as

$$\alpha = \sum_{k=1}^{2} \langle \beta, \alpha_k \rangle_p \alpha_k = \langle (1\ 1\ 1\ 1), (0\ 1\ 0\ 1) \rangle_p \alpha_1 + \langle (1\ 1\ 1\ 1), (1\ 0\ 1\ 1) \rangle_p \alpha_2.$$

$$= (1\ 0\ 1\ 1)\ .$$

Thus the approximated code word is (1 0 1 1).

Now having seen some simple properties of codes we now proceed on to define super special vector spaces. This new algebraic structure basically makes used of supermatrices. For more super linear algebra please refer [10, 21].

**Chapter Two**

# SUPER SPECIAL VECTOR SPACES

In this chapter we for the first time define a new class of super special vector spaces. We describe mainly those properties essential for us to define super special codes and their properties like decoding, etc. Throughout this chapter V denotes a vector space over a field F. F may be a finite characteristic field or an infinite characteristic field.

**DEFINITION 2.1:** *Let $V_s = [V_1 \mid V_2 \mid \ldots \mid V_n]$ where each $V_i$ is a vector space of dimension m over a field F; i = 1, 2, …, n, then we call $V_s$ to be the super special finite dimensional vector space over F. Any element $v_s \in V_s$ would be of the form $\left[v_1^1 v_2^1 \ldots v_m^1 \mid v_1^2 v_2^2 \ldots v_m^2 \mid \ldots \mid v_1^n v_2^n \ldots v_m^n\right]$ where $v_i^t \in F$ ; $1 \le i \le m$ and $1 \le t \le n$ and $\left(v_1^t v_2^t \ldots v_m^t\right) \in V_t$. $V_t$ a vector space of dimension m over F i.e., $V_t \cong \underbrace{F \times F \times \ldots \times F}_{\text{m - times}}$ is true for t = 1, 2,*

*…, n. Thus any element $v_s \in V_s$ is a super row vector with entries from the field F. If $v_s$, $w_s \in V_s$ the sum $v_s + w_s$ is defined to be*

$$[\, v_1^1 + w_1^1 \;\; v_2^1 + w_2^1 \ldots v_m^1 + w_m^1 \;\Big|\; v_1^2 + w_1^2 \;\; v_2^2 + w_2^2 \ldots v_m^2 + w_m^2 \;\Big|\; \ldots \;\Big|\; v_1^n + w_1^n \;\; v_2^n + w_2^n \ldots v_m^n + w_m^n \,]$$

*where*

$$w_s = [\, w_1^1\, w_2^1 \ldots w_m^1 \;\Big|\; w_1^2 w_2^2 \ldots w_m^2 \;\Big|\; \ldots \;\Big|\; w_1^n w_2^n \ldots w_m^n \,] .$$

*Also*

$$av_s = [\, av_1^1 \;\; av_2^1 \ldots av_m^1 \;\Big|\; av_1^2 \;\; av_2^2 \ldots av_m^2 \;\Big|\; \ldots \;\Big|\; av_1^n \;\; av_2^n \ldots av_m^n \,],$$

$$a \in F.$$

We illustrate this by a simple example.

***Example 2.1:*** Let $V_s = [V_1 \mid V_2 \mid V_3]$ where each $V_i$ is a 3 dimensional vector space over Q, the field of rationals. Any element $v_s \in V_s$ would be of the form

$$v_s = \left[ a_1^1 \;\; a_2^1 \;\; a_3^1 \,|\, a_1^2 \;\; a_2^2 \;\; a_3^2 \,|\, a_1^3 \;\; a_2^3 \;\; a_3^3 \right]$$

is a super row vector; where $a_j^i \in Q$; $1 \le i, j \le 3$. Clearly $V_s$ is a super special finite dimensional super vector space over Q.

$B_s$ = {[0 0 1 | 0 0 1 | 0 0 1], [0 0 1 | 0 0 1 | 0 1 0],
[0 0 1 | 0 0 1 | 1 0 0], [0 0 1 | 0 1 0 | 0 0 1],
[0 0 1 | 0 1 0 | 0 1 0], [0 0 1 | 0 1 0 | 1 0 0],
[0 0 1 | 1 0 0 | 0 0 1], [0 0 1 | 1 0 0 | 1 1 0],
[0 0 1 | 1 0 0 | 1 0 0], [0 1 0 | 1 0 0 | 0 0 1],
[0 1 0 | 0 0 1 | 0 0 1], [0 1 0 | 0 0 1 | 0 1 0],
[0 1 0 | 0 1 0 | 0 1 0], [0 1 0 | 0 1 0 | 1 0 0],
[0 1 0 | 1 0 0 | 1 0 0], [0 1 0 | 1 0 0 | 0 1 0],
[1 0 0 | 1 0 0 | 0 1 0], [1 0 0 | 1 0 0 | 0 0 1],
[1 0 0 | 1 0 0 | 1 0 0], [1 0 0 | 0 1 0 | 0 0 1],
[1 0 0 | 0 1 0 | 0 1 0], [0 1 0 | 0 0 1 | 1 0 0],
[0 1 0 | 0 1 0 | 0 0 1], [1 0 0 | 0 1 0 | 1 0 0],
[1 0 0 | 0 0 1 | 1 0 0], [1 0 0 | 0 0 1 | 0 1 0]

and

[1 0 0 | 0 0 1 | 0 0 1]};

forms a super basis of $V_s$. Clearly this $B_s$ will generate $V_s$ over Q. It is easily verified $B_s \subset V_s$ and the elements of $B_s$ is a linearly independent set in $V_s$.

We have seen a super special finite dimensional vector space over Q. Now we proceed on to define the super special basis of a super special vector space over F.

**DEFINITION 2.2:** *Let $V_s = [V_1 \mid V_2 \mid \ldots \mid V_n]$ be a super special vector space where each $V_i$ is of dimension m over the field F. Let $B_s = \{ v_s^1, v_s^2, \ldots, v_s^t ; 0 < t \le m^n \}$ be the elements from $V_s$ i.e., each $v_s^t$ is a super row vector $1 \le p \le t$. We say $v_s^1, v_s^2, \ldots, v_s^t$ forms a super linearly independent set if, $\alpha_1, \ldots, \alpha_t$ are scalars in F such that*

$$\alpha_1 v_s^1 + \alpha_2 v_s^2 + \ldots + \alpha_t v_s^t = (0 \ldots 0 \mid 0 \ldots 0 \mid \ldots \mid 0 \ldots 0) \ldots \quad I$$

*then each $\alpha_i = 0$. If equation I is true for some non zero scalars $\alpha_1, \ldots, \alpha_t$ in F then we say $v_s^1, v_s^2, \ldots, v_s^t$ forms a super linearly dependent set.*

*If $B_s$ forms a super linearly independent set and if every element in $V_s$ can be expressed as a super linear combination from the set $B_s$ in a unique way then we call $B_s$ to be a special super basis of $V_s$ or super special basis of $V_s$. If the number of elements in $B_s$ is finite then we call $V_s$ to be a finite dimensional super special vector space otherwise an infinite dimensional super special vector space.*

***Example 2.2:*** Let $V_s = [V_1 \mid V_2 \mid V_3]$ where each $V_i$ is a two dimensional vector space over Q, i = 1, 2, 3; be the super special vector space i.e., $V_i \cong Q \times Q$; $1 \le i \le 3$.

We see

$B_s$ = {[1 0 | 1 1 | 0 1], [2 0 | 1 0 | 3 1], [3 0 | 2 1 | 3 2] $\subset V_s$

is a super linearly dependent set as [1 0 | 1 1 | 0 1] + [2 0 | 1 0 | 3 1] = [3 0 | 2 1 | 3 2].

Suppose we consider a subset

$C_s$ = {[1 0 | 1 0 | 1 0], [0 1 | 1 0 | 1 0],
[1 0 | 1 0 | 0 1], [1 0 | 0 1 | 0 1],

$$[1\ 0 \mid 0\ 1 \mid 1\ 0]\ [0\ 1 \mid 1\ 0 \mid 0\ 1],$$
$$[0\ 1 \mid 0\ 1 \mid 1\ 0],\ [0\ 1 \mid 0\ 1 \mid 0\ 1]\}$$

of $V_s$; clearly $C_s$ is a super linearly independent set as well as $C_s$ is a super basis of $V_s$.

It is a matter of routine to verify a super special vector space which is finite dimensional will have only the same number of elements in every super special basis of $V_s$.

Next we define a super special subvector space of a super special vector space $V_s = [V_1 \mid V_2 \mid \ldots \mid V_n]$.

**DEFINITION 2.3:** *Let $V_s = [V_1 \mid V_2 \mid \ldots \mid V_n]$ be a super special vector space; where each $V_i$ is of dimension m. A non empty subset $W_s = \{[W_1 \mid W_2 \mid \ldots \mid W_n]\}$ of $V_s$ is said to be a super special subspace of $V_s$ if each $W_i$ is a subspace of the vector space $V_i$ of dimension k, k < m for i = 1, 2, …, n.*

We give an example of a super special subspace of a super special vector space.

***Example 2.3:*** Let $V_s = [M_{2 \times 2} \mid Q \times Q \times Q \times Q] = [V_1 \mid V_2]$ be a super special vector space over the field Q. Clearly both $V_1$ and $V_2$ are vector spaces of dimension 4 over Q. Let $W_s = \{[W_1 \mid W_2]\}$ where

$$W_1 = \left\{ \begin{pmatrix} a & 0 \\ 0 & b \end{pmatrix} \middle|\ a, b \in Q \right\}$$

and

$$W_2 = Q \times Q\ \{0\} \times \{0\}$$

are subspaces of $V_1$ and $V_2$ respectively of dimension 2. $W_s$ is super special subspace of $V_s$.

Now we define the notion of super special mixed dimension vector space.

**DEFINITION 2.4:** *Let $V_s = [V_1 \mid V_2 \mid \ldots \mid V_n]$, where each $V_i$ is a finite dimensional vector space of dimension $m_i$ over a field F ;*

*$m_i \neq m_j$ for at least one $i \neq j$ for $i = 1, 2, \ldots, n$, $1 \leq j \leq n$. Then $V_s$ is defined to be a super special mixed dimension vector space.*

***Example 2.4:*** Let $V_s = [M_{2\times 2} \mid Q \times Q \times Q \times Q \mid V_3]$ where $V_3$ is Q[x] and Q[x] is the set of all polynomials of degree less than or equal to 3 with coefficients from Q. We see $M_{2\times 2} = V_1$ is a vector space of dimension 4 and $V_3$ is also a vector space of dimension 4.

Thus we see $M_{2\times 2} \cong Q \times Q \times Q \times Q$; $V_2 \cong Q \times Q \times Q \times Q$ and $V_3 \cong Q \times Q \times Q \times Q$. So any $v_s \in V_s$ will be of the form

$$\left[q_1^1\ q_2^1\ q_3^1\ q_4^1 \mid q_1^2\ q_2^2\ q_3^2\ q_4^2 \mid q_1^3\ q_2^3\ q_3^3\ q_4^3\right]$$

where $q_j^i \in Q$ with i = 1, 2, 3, 4.

Let us consider

$$B_s = \{[v_1^1\ v_1^2\ v_1^3], [v_1^1\ v_1^2\ v_2^3], [v_1^1\ v_1^2\ v_3^3], \ldots, [v_4^1\ v_4^2\ v_4^3]\}$$

where

$$v_1^1 = \begin{bmatrix} 1 & 0 \\ 0 & 0 \end{bmatrix}, v_2^1 = \begin{bmatrix} 0 & 1 \\ 0 & 0 \end{bmatrix}, v_3^1 = \begin{bmatrix} 0 & 0 \\ 0 & 1 \end{bmatrix}, v_4^1 = \begin{bmatrix} 0 & 0 \\ 0 & 1 \end{bmatrix};$$

$$v_1^2 = [0\ 0\ 0\ 1],\ v_2^2 = [0\ 0\ 1\ 0],\ v_3^2 = [0\ 1\ 0\ 0],\ v_4^2 = [1\ 0\ 0\ 0],$$
$$v_1^3 = 1,\ v_2^3 = x,\ v_3^3 = x^2 \text{ and } v_4^3 = x^3. \text{ Clearly } |B_s| = 64,$$

$B_s \subset V_s$ and $B_s$ is a special super basis of $V_s$ over Q.

We may have any other basis for $V_s$ but it is true as in case of vector spaces even in special super vector spaces the number of elements in each and every superbasis is only fixed and equal. Further as in case of usual vector spaces even in case of super special vector spaces $V_s$, every element of $V_s$ can be represented by a super row vector when the dimension of each $V_i$ is the same and would be represented by a super mixed row vector; in case the dimension of each of the vector spaces $V_i$'s are different.

Now we illustrate by an example the notion of super special mixed dimension vector space.

***Example 2.5:*** Let $V_s = [V_1 \mid V_2 \mid V_3 \mid V_4]$ be a special super mixed dimension vector space over Q, where $V_1 = Q \times Q \times Q$ is a vector space of dimension 3 over Q, $V_2 = M_{3 \times 2}$; the set of all 3 × 2 matrices with entries from Q is a vector space of dimension 6 over Q, $V_3 = Q[x]$ ; the set of all polynomials of degree less than or equal to four and $V_4 = P_{2 \times 2}$ the set of all 2 × 2 matrices with entries from Q. $V_4$ is a vector space of dimension 4 over Q. Thus $V_s$ is a super special mixed dimensional vector space over Q. Any element $v_s \in V_s$ is a super mixed row vector given by

$$\begin{aligned} v_s = \{[v_1^1\ v_2^1\ v_3^1 \mid v_1^2\ v_2^2\ v_3^2\ v_4^2\ v_5^2\ v_6^2 \mid \\ v_1^3\ v_2^3\ v_3^3\ v_4^3\ v_5^3 \mid v_1^4\ v_2^4\ v_3^4\ v_4^4]\} \end{aligned}$$

where $v_j^i \in V_i$ ; $1 \leq j \leq 3, 4, 5$ or 6 and $1 \leq i \leq 4$.

A super special subspace of $V_s$ can be either a super special subspace or it can also be a super special mixed dimension subspace.

Now it is important at this point to mention even a super special vector space can have a super special mixed dimension subspace also a super special mixed dimension vector space can have a super special vector subspace.

We illustrate this situation now by the following example.

***Example 2.6:*** Let $V_s = [V_1 \mid V_2 \mid V_3]$ where $V_1$ is the set of all 2 × 3 matrices with entries from Q, $V_2$ is the set of all 1 × 6 row vector with entries from Q and $V_3$ is the collection of all polynomials of degree less than or equal to five with coefficients from Q. All the three vector spaces $V_1$, $V_2$ and $V_3$ are of dimensions 6 over Q. $V_s$ is a super special vector space over Q.

Take $W_s = [W_1 \mid W_2 \mid W_3]$ to be a proper subset of $V_s$ where

$$W_1 = \left\{ \begin{bmatrix} a_1 & a_2 & a_3 \\ 0 & 0 & a_4 \end{bmatrix} \middle| \ a_1, a_2, a_3, a_4 \in Q \right\}.$$

$W_1$ is a subspace of $V_1$ of dimension 4.

$$W_2 = [[x_1\ 0\ x_2\ 0\ x_3\ 0] \mid x_1, x_2, x_3 \in Q\}$$
$$\cong Q \times \{0\} \times Q \times \{0\} \times Q \times \{0\},$$

$W_2$ is a subspace of $V_2$ of dimension 3. Let $W_3$ = {The collection of all polynomials of even degree (i.e., degree 2 and 4) with coefficients from Q}. $W_3$ is a subspace of $V_3$ of dimension 3 over Q.

Clearly $W_s$ is a special super subvector space of varying dimension or $W_s$ is a super special mixed dimension subspace of $V_s$ but $V_s$ is not a super special mixed dimension vector space, it is only a super special vector space over Q.
Consider $T_s = [T_1 \mid T_2 \mid T_3]$ a proper subset of $V_s$, where

$$T_1 = \left\{ \begin{bmatrix} a_1 & a_2 & a_3 \\ 0 & 0 & 0 \end{bmatrix} \middle| \, a_1, a_2, a_3 \in Q \right\}$$

a proper subspace of $V_1$ of dimension three.

$$T_2 = \{[a_1\ a_2\ a_3\ 0\ 0\ 0] \mid a_1, a_2, a_3 \in Q\}$$
$$\cong Q \times Q \times Q \times \{0\} \times \{0\} \times \{0\}$$

is a proper subspace of $V_2$ of dimension 3 over Q. Let $T_3$ = {all polynomials of degree 1 and 3 with coefficients over Q} = $\{[a_0 + a_1x + a_2x^3] \mid a_0, a_1, a_2 \in Q\}$. $T_3$ is a subspace of $V_3$ of dimension three over Q. Thus $T_s$ is a super special subvector space of $V_s$ over Q and the dimension of each $T_i$ is three; i = 1, 2, 3.

Next we proceed on to give an example of a super special mixed dimension vector space having subspaces which are super special mixed dimension subspace and super special subspace.

***Example 2.7:*** Let $V_s = [V_1 \mid V_2 \mid V_3]$ be a super special mixed dimension vector space over Q; where $V_1 = \{M_{3\times 3} = (m_{ij}) \mid m_{ij} \in Q, 1 \le i, j \le 3\}$, vector space of dimension 9 over Q. $V_2 = \{Q \times Q \times Q \times Q = (a, b, c, d) \mid a, b, c, d \in Q\}$ a vector space of dimension four over Q and $V_3 = \{M_{2\times 2} = (a_{ij}) \mid a_{ij} \in Q, 1 \leq i, j \leq 2\}$ a space of dimension four over Q.

Now let $W_s = [W_1 \mid W_2 \mid W_3]$ where $W_1$ is a diagonal matrix of the form

$$\left\{\begin{bmatrix} a & 0 & 0 \\ 0 & b & 0 \\ 0 & 0 & c \end{bmatrix} \middle| a,\ b,\ c \in Q\right\}$$

which is proper subspace of $V_1$ and of dimension 3,

$$W_2 = \{[a\ b\ c\ 0] \mid a, b, c, 0 \in Q\} \cong Q \times Q \times Q \times \{0\}$$

is a subspace of $V_2$ of dimension three over Q and

$$W_3 = \left\{\begin{bmatrix} a & b \\ 0 & c \end{bmatrix} \middle| a, b, c \in Q\right\}$$

is the proper subspace of dimension three of $V_3$. $W_s = [W_1 \mid W_2 \mid W_3]$ is a super special subspace of $V_s$. Clearly $W_s$ is not a super special mixed dimensional subspace of $V_s$.

Let $R_s = [R_1 \mid R_2 \mid R_3]$, a proper subset of $V_s$, where $R_1$ = {Set of all 3 × 3 upper triangular matrices} i.e.,

$$R_1 = \left\{\begin{bmatrix} a & b & c \\ 0 & d & e \\ 0 & 0 & f \end{bmatrix} \middle| a, b, c, d, e, f \in Q\right\}.$$

$R_1$ is a subspace of $V_1$ of dimension 6 over Q.

$$R_2 = \{[a\ 0\ b\ 0] \mid a, b \in Q\}; \cong Q \times \{0\} \times Q \times \{0\}.$$

$R_2$ is a subspace of $V_2$ of dimension 2 over Q.

Let

$$R_3 = \left\{ \begin{bmatrix} a & 0 \\ 0 & b \end{bmatrix} \middle| \ a, b \in Q \right\};$$

$R_3$ is a subspace of $V_3$ of dimension 2 over Q. We see $R_s$ is a super special mixed dimension subspace of $V_s$.

We mention here only those factors about the super special vector spaces which are essential for the study and introduction of super special codes.

Now we proceed on to define dot product of super special vector spaces.

**DEFINITION 2.5:** *Let $V_s$ be a real super special vector space over the field of reals F. A super special inner product or super inner product on $V_s$ is a function which assigns to each ordered pair of super row vectors $\alpha_s$, $\beta_s$ in $V_s$ a scalar in F in such a way that for $\alpha_s$, $\beta_s$, $\gamma_s$ in $V_s$ and for all scalars c in F we have*

*1. $(\alpha_s + \beta_s / \gamma_s) = (\alpha_s / \gamma_s) + (\beta_s / \gamma_s)$ where*

$$(\alpha_s / \gamma_s) = \left( \left[ v_1^1 v_2^1 \ldots v_{m_1}^1 \middle| v_1^2 v_2^2 \ldots v_{m_2}^2 \middle| \ldots \middle| v_1^n v_2^n \ldots v_{m_n}^n \right] \left[ w_1^1 w_2^1 \ldots w_{m_1}^1 \middle| w_1^2 w_2^2 \ldots w_{m_2}^2 \middle| \ldots \middle| w_1^n w_2^n \ldots w_{m_n}^n \right] \right)$$

$$= \left( v_1^1 w_1^1 + v_2^1 w_2^1 + \ldots + v_{m_1}^1 w_{m_1}^1 \right) + \left( v_1^2 w_1^2 + v_2^2 w_2^2 + \ldots + v_{m_2}^2 w_{m_2}^2 \right) + \ldots + \left( v_1^n w_1^n + v_2^n w_2^n + \ldots + v_{m_n}^n w_{m_n}^n \right)$$

$$= a_1 + a_2 + \ldots + a_n = c \in F.$$

*Here*

$$\alpha_s = \left[ v_1^1 v_2^1 \ldots v_{m_1}^1 \middle| v_1^2 v_2^2 \ldots v_{m_2}^2 \middle| \ldots \middle| v_1^n v_2^n \ldots v_{m_n}^n \right]$$

*and*

$$\gamma_s = \left[ w_1^1 w_2^1 \ldots w_{m_1}^1 \middle| w_1^2 w_2^2 \ldots w_{m_2}^2 \middle| \ldots \middle| w_1^n w_2^n \ldots w_{m_n}^n \right].$$

2. *$c(\alpha_s \mid \beta_s) = (c\alpha_s \mid \beta_s)$*
3. *$(\alpha_s \mid \beta_s) = (\beta_s \mid \alpha_s)$*
4. *$(\alpha_s \mid \alpha_s) > 0;\ \alpha_s \neq 0.$*

*We will say a super special vector space endowed with super special inner product as a super special inner product space. Suppose $V_s = [V_1 \mid V_2 \mid \ldots \mid V_n]$ is a super special mixed dimension vector space on which is endowed an inner product; we say $w_s, v_s \in V_s$ is orthogonal if $(v_s \mid w_s) = 0$.*

While defining super special codes. We may need the notion of orthogonality. As in case of usual vector spaces we see $(\alpha_s \mid \alpha_s) = \| \alpha_s \|^2$ where $\alpha_s \in V_s$ is defined as the super special norm. Also we see $V_s$ is an abelian group with respect to addition. If $W_s$ is a subspace of $V_s$ it is necessarily a subgroup of $V_s$.

Now we can define for any $x_s \in V_s$ the coset of $W_s$ by $x_s + W_s = \{x_s + w_s \mid w_s \in W_s\}$.

When we are carrying this inner product to super special vector spaces over finite fields condition 4 may not in general be true.

We may call a super special inner product in which $(\alpha_s \mid \alpha_s) = 0$ even if $\alpha_s \neq 0$ as a pseudo super special inner product, all the other conditions 1, 2 and 3 being true.

We illustrate this by the following example.

***Example 2.8:*** Let $V_s = [V_1 \mid V_2 \mid V_3]$ be a special super vector space over Q where

$$V_1 = Q \times Q,\ V_2 = Q \times Q \times Q \text{ and } V_3 = Q \times Q \times Q.$$

Let $W_s = [W_1 \mid W_2 \mid W_3] \subseteq V_s$ be a proper special super subspace of $V_s$ where

$$W_1 = Q \times \{0\};\ W_2 = Q \times Q \times \{0\} \text{ and } W_3 = Q \times \{0\} \times Q.$$

Now coset of $W_s$ related to $x_s = \{[7\ 3 \mid 1\ 2\ 3 \mid 5\ 7\ 1]\}$ in $V_s$ is given by $x_s + W_s = \{x_s + w_s \mid w_s \in W_s\}$. It is easily verified that $x_s + W_s$ with varying $x_s \in V_s$ partitions $V_s$ as cosets. This property will also be used in super special codes.

Since we are interested in only finite dimensional super special vector spaces that too defined over finite characteristic field we give some examples of them. As linear codes in most of the cases are binary codes we will be giving examples only using the field $Z_2 = \{0, 1\}$, the prime field of characteristic two.

***Example 2.9:*** Let $V_s = [V_1 \mid V_2 \mid V_3]$ be a super special vector space over $Z_2 = \{0, 1\}$, where $V_1 = Z_2 \times Z_2 \times Z_2$, $V_2 = Z_2 \times Z_2 \times Z_2$ and $V_3 = Z_2 \times Z_2 \times Z_2$ are vector spaces over $Z_2$. Each of the spaces are of dimension 3 over $Z_2$. Let $W_s = [W_1 \mid W_2 \mid W_3]$ be a super special subvector space of $V_s$ over $Z_2$ where $W_1 = Z_2 \times \{0\} \times Z_2$, $W_2 = \{0\} \times Z_2 \times Z_2$ and $W_3 = Z_2 \times Z_2 \times \{0\}$. Clearly

$$
\begin{aligned}
W_s = \{&[0\ 0\ 0 \mid 0\ 0\ 0 \mid 0\ 0\ 0], [1\ 0\ 0 \mid 0\ 0\ 0 \mid 0\ 0\ 0], \\
&[1\ 0\ 0 \mid 0\ 0\ 0 \mid 0\ 1\ 0], [1\ 0\ 0 \mid 0\ 1\ 1 \mid 0\ 0\ 0], \\
&[1\ 0\ 0 \mid 0\ 1\ 0 \mid 0\ 0\ 0], [1\ 0\ 0 \mid 0\ 0\ 1 \mid 0\ 0\ 0], \\
&[1\ 0\ 0 \mid 0\ 1\ 0 \mid 0\ 1\ 0], [1\ 0\ 0 \mid 0\ 1\ 1 \mid 1\ 0\ 0], \\
&[1\ 0\ 0 \mid 0\ 1\ 1 \mid 1\ 1\ 0], [1\ 0\ 0 \mid 0\ 1\ 1 \mid 0\ 1\ 0], \\
&[1\ 0\ 0 \mid 0\ 0\ 1 \mid 1\ 1\ 0], [1\ 0\ 0 \mid 0\ 0\ 1 \mid 1\ 0\ 0], \\
&[1\ 0\ 0 \mid 0\ 0\ 1 \mid 0\ 1\ 0] \ \ldots \text{ and so on}\}.
\end{aligned}
$$

We see the number of elements in $W_s$ is 64. Suppose

$$x_s = [1\ 1\ 1 \mid 1\ 1\ 1 \mid 1\ 1\ 1] \in V_s .$$

We see $x_s \notin W_s$. We can find $x_s + W_s = \{x_s + w_s \mid w_s \in W_s\}$. Clearly $V_s$ is partitioned into 8 disjoint sets relative to the super special subspace $W_s$ of $V_s$.

Let

$$
\begin{aligned}
B_s = \{&[1\ 0\ 0 \mid 1\ 0\ 0 \mid 1\ 0\ 0], [1\ 0\ 0 \mid 1\ 0\ 0 \mid 0\ 1\ 0], \\
&[1\ 0\ 0 \mid 1\ 0\ 0 \mid 0\ 0\ 1], [1\ 0\ 0 \mid 0\ 1\ 0 \mid 1\ 0\ 0], \\
&[1\ 0\ 0 \mid 0\ 1\ 0 \mid 0\ 1\ 0], [1\ 0\ 0 \mid 0\ 1\ 0 \mid 0\ 0\ 1], \\
&[1\ 0\ 0 \mid 0\ 0\ 1 \mid 1\ 0\ 0], [1\ 0\ 0 \mid 0\ 0\ 1 \mid 0\ 1\ 0], \\
&[1\ 0\ 0 \mid 0\ 0\ 1 \mid 0\ 0\ 1], [0\ 1\ 0 \mid 1\ 0\ 0 \mid 1\ 0\ 0], \\
&[0\ 1\ 0 \mid 1\ 0\ 0 \mid 0\ 1\ 0], [0\ 1\ 0 \mid 1\ 0\ 0 \mid 0\ 0\ 1], \\
&[0\ 1\ 0 \mid 0\ 1\ 0 \mid 1\ 0\ 0], [0\ 1\ 0 \mid 0\ 1\ 0 \mid 0\ 0\ 1],
\end{aligned}
$$

[0 1 0 | 0 1 0 | 0 1 0], [0 1 0 | 0 0 1 | 1 0 0],
[0 1 0 | 0 0 1 | 0 1 0], [0 1 0 | 0 0 1 | 0 0 1],
[0 0 1 | 1 0 0 | 0 0 1], [0 0 1| 1 0 0 | 0 1 0],
[0 0 1 | 1 0 0 | 1 0 0], [0 0 1 | 0 1 0 | 0 0 1],
[0 0 1| 0 1 0 | 0 1 0], [0 0 1 | 0 1 0 | 1 0 0],
[0 0 1 | 0 0 1 | 1 0 0], [0 0 1 | 0 0 1 | 0 1 0],
and [0 0 1 | 0 0 1 | 0 0 1]}.

$B_s$ is a super special basis of the super special vector space $V_s$. Let

$T_s$ = {[1 0 0 | 0 1 0 | 1 0 0], [1 0 0 | 0 1 0 | 0 1 0],
[1 0 0 | 0 0 1 | 1 0 0], [1 0 0 | 0 1 0 | 0 1 0],
[0 0 1 | 0 1 0 | 1 0 0], [0 0 1 | 0 1 0 | 0 1 0],
[0 0 1 | 0 0 1 | 1 0 0] and [0 0 1 | 0 0 1 | 0 1 0]}

$T_s$ is a super special basis of $W_s$ and the number of elements in $T_s$ is 8.

Let us take $R_s = [R_1 \mid R_2 \mid R_3]$ where $R_1 = \{0\} \times Z_2 \times Z_2$, $R_2 = \{0\} \times \{0\} \times Z_2$ and $R_3 = \{0\} \times Z_2 \times Z_2$. Clearly $R_s$ is a super special mixed dimension subspace of $V_s$. Let

$M_s$ = {[0 0 1 | 0 0 1 | 0 0 1], [0 0 1 | 0 0 1| 0 1 0],
[0 1 0 | 0 0 1 | 0 0 1] and [0 1 0 | 0 0 1 | 0 1 0],

$M_s$ is a super special basis of $R_s$. We see $R_s$ is a subspace of dimension 4. Further we see all the super special base elements are only super row vectors.

We give yet another example of a super special mixed dimensional vector space over $Z_2$.

***Example 2.10:*** Let $V_s = [V_1 \mid V_2 \mid V_3]$ be a super special vector space over the field $Z_2 = \{0, 1\}$
where

$$V_1 = Z_2 \times Z_2,$$
$$V_2 = Z_2 \times Z_2 \times Z_2$$

and

$$V_3 = Z_2 \times Z_2 \times Z_2 \times Z_2.$$

It is easily verified that

$$B_s = \{[1\ 0 \mid 1\ 0\ 0 \mid 1\ 0\ 0\ 0], [0\ 1 \mid 1\ 0\ 0 \mid 1\ 0\ 0\ 0],$$
$$[1\ 0 \mid 0\ 1\ 0 \mid 1\ 0\ 0\ 0], [0\ 1 \mid 0\ 1\ 0 \mid 1\ 0\ 0\ 0],$$
$$[1\ 0 \mid 0\ 0\ 1 \mid 1\ 0\ 0\ 0], [0\ 1 \mid 0\ 0\ 1 \mid 1\ 0\ 0\ 0],$$
$$[1\ 0 \mid 1\ 0\ 0 \mid 0\ 1\ 0\ 0], [1\ 0 \mid 0\ 0\ 1 \mid 0\ 1\ 0\ 0],$$
$$[0\ 1 \mid 1\ 0\ 0 \mid 0\ 1\ 0\ 0], [0\ 1 \mid 0\ 1\ 0 \mid 0\ 1\ 0\ 0],$$
$$[1\ 0 \mid 0\ 1\ 0 \mid 0\ 1\ 0\ 0], [0\ 1 \mid 0\ 0\ 1 \mid 0\ 1\ 0\ 0],$$
$$[1\ 0 \mid 1\ 0\ 0 \mid 0\ 0\ 1\ 0], [0\ 1 \mid 1\ 0\ 0 \mid 0\ 0\ 1\ 0],$$
$$[1\ 0 \mid 0\ 1\ 0 \mid 0\ 0\ 1\ 0], [0\ 1 \mid 0\ 1\ 0 \mid 0\ 0\ 1\ 0],$$
$$[1\ 0 \mid 0\ 0\ 1 \mid 0\ 0\ 1\ 0], [0\ 1 \mid 0\ 0\ 1 \mid 0\ 0\ 1\ 0],$$
$$[1\ 0 \mid 1\ 0\ 0 \mid 0\ 0\ 0\ 1], [0\ 1 \mid 1\ 0\ 0 \mid 0\ 0\ 0\ 1],$$
$$[1\ 0 \mid 0\ 1\ 0 \mid 0\ 0\ 0\ 1], [0\ 1 \mid 0\ 1\ 0 \mid 0\ 0\ 0\ 1],$$
$$[1\ 0 \mid 0\ 0\ 1 \mid 0\ 0\ 0\ 1] \text{ and } [0\ 1 \mid 0\ 0\ 1 \mid 0\ 0\ 0\ 1]\} \subseteq V_s;$$

is a super special basis of $V_s$ and the number of elements in $B_s$ is 24 = Base elements of $V_1$ × Base elements of $V_2$ × Base elements of $V_3$ = 2 × 3 × 4 = 24.

Let us consider a super special subspace of $V_s$ say

$$W_s = [W_1 \mid W_2 \mid W_3]$$

where

$$W_1 = Z_2 \times \{0\},$$
$$W_2 = Z_2 \times Z_2 \times \{0\}$$

and

$$W_3 = \{0\} \times \{0\} \times Z_2 \times Z_2$$

are vector subspace of $V_1$, $V_2$ and $V_3$ respectively.

Let

$$R_s = \{[1\ 0 \mid 1\ 0\ 0 \mid 0\ 0\ 0\ 1], [1\ 0 \mid 1\ 0\ 0 \mid 0\ 0\ 1\ 0],$$
$$[1\ 0 \mid 0\ 1\ 0 \mid 0\ 0\ 0\ 1], [1\ 0 \mid 0\ 1\ 0 \mid 0\ 0\ 1\ 0]$$
$$\subset W_s.$$

Clearly $R_s$ is a super special subvector space of $W_s$ and the super special dimension of $W_s$ is 4.

Now the following factors are easily verified to be true. If $V_s$ = [$V_1$ | $V_2$ | … | $V_n$] is a super special mixed dimension vector space over a field F and if $V_i$ is of dimension $n_i$ over F, i = 1, 2, …, n, then the dimension of $V_s$ = dimension of $V_1$ × dimension of $V_2$ × … × dimension of $V_n$ = $n_1 \times n_2 \times \ldots \times n_n$.

In the next chapter we define the notion of super special codes, which are built mainly using these super special vector spaces.

## Chapter Three

# SUPER SPECIAL CODES

In this chapter for the first time we define new classes of super special codes; using super matrices and enumerate some of their error correcting and error detecting techniques. However their uses and applications would be given only in chapter four. This chapter has three sections. Section one introduces super special row codes and super special column codes are introduced in section two. Section three defines the new notion of super special codes and discusses their properties.

## 3.1 Super Special Row Codes

In this section we define two new classes of super special row codes using super row matrix and super mixed row matrix.

We just say a super matrix M = $[V_1 \mid V_2 \mid \ldots \mid V_r]$ is known as the super row vector or matrix if each $V_i$ is a $n \times m$ matrix so that M can be visualized as a $n \times \underbrace{(m + \ldots + m)}_{r-\text{times}}$ matrix where partitions are done vertically between the m and $(m + 1)^{th}$ column, 2m and $(2m + 1)^{th}$ column and so on and lastly $(r – 1)m$ and $\{(r – 1)m +1\}^{th}$ column.

For example

$$\left[\begin{array}{cccc|cccc|cccc} 1 & 2 & 3 & 1 & 4 & 5 & 6 & 7 & 1 & 0 & 0 & 1 \\ 4 & 5 & 6 & 2 & 8 & 9 & 0 & 1 & 2 & 1 & 0 & 3 \\ 1 & 0 & 1 & 3 & 1 & 1 & 0 & 0 & 1 & 1 & 0 & 2 \end{array}\right]$$

is a super row vector or matrix, here r = 3, n = 3 and m = 4.

A super mixed row vector or matrix $V = [V_1 | V_2 | \ldots | V_s]$ is a super matrix such that each $V_i$ is a $n \times m_i$ matrix $m_i \neq m_j$ for at least one $i \neq j$, $1 \leq i, j \leq s$ i.e., V is a $n \times (m_1 + \ldots + m_s)$ matrix where vertical partitions are made between $m_1$ and $(m_1 + 1)^{th}$ column of V, $m_2$ and $(m_2 + 1)^{th}$ column and so on. Lastly partitioned between the $m_{s-1}$ and $(m_{s-1} + 1)^{th}$ column.

For example

$$V = \left[\begin{array}{cc|cccc|ccc|ccccc} 1 & 0 & 1 & 1 & 0 & 1 & 1 & 7 & 2 & 1 & 2 & 3 & 4 & 5 \\ 2 & 1 & 2 & 0 & 5 & 0 & 3 & 5 & 1 & 0 & 1 & 0 & 1 & 0 \\ 3 & 1 & 1 & 1 & 0 & 1 & 1 & 2 & 3 & 1 & 1 & 1 & 0 & 1 \\ 4 & 2 & 0 & 0 & 1 & 0 & 4 & 5 & 6 & 0 & 0 & 1 & 1 & 0 \end{array}\right]$$

where n = 4, $m_1 = 2$, $m_2 = 4$, $m_3 = 3$ and $m_4 = 5$ and s = 4. Clearly V is a super mixed row matrix or vector. For more refer chapter one of this book.

Now we proceed on to define the new class of super special row code.

**DEFINITION 3.1.1:** *Suppose we have to transform some n set of $k_1, \ldots, k_n$ message symbols $a_1^1 a_2^1 \ldots a_{k_1}^1$, $a_1^2 a_2^2 \ldots a_{k_2}^2$, …, $a_1^n a_2^n \ldots a_{k_n}^n$, $a_i^t \in F_q$; $1 \leq t \leq n$ and $1 \leq i \leq k_i$ (q a power of a prime) as a set of code words simultaneously into n-code words such that each code word is of length $n_i$, i = 1, 2, …, n and $n_1 – k_1 = n_2 – k_2 = \ldots = n_n – k_n = m$ say i.e., the number check symbols of every* *code word is the same i.e., the number of message symbols and the length of the code word may not be the same. That is the code word consisting of n code words can be represented as a super row vector;*

$$x_s = \left[ x_1^1 x_2^1 \ldots x_{n_1}^1 \;\middle|\; x_1^2 x_2^2 \ldots x_{n_2}^2 \;\middle|\; \ldots \;\middle|\; x_1^n x_2^n \ldots x_{n_n}^n \right]$$

*$n_i > k_i$, $1 \le i \le n$. In this super row vector $x_j^i = a_j^i$, $1 \le j \le k_i$; i = 1, 2, …, n and the remaining $n_i - k_i$ elements $x_{k_i+1}^i\ x_{k_i+2}^i \ldots x_{n_i}^i$ are check symbols or control symbols; i = 1, 2, …, n.*

These n code words denoted collectively by $x_s$ will be known as the super special row code word.

As in case of usual code, the check symbols can be obtained in such a way that the super special code words $x_s$ satisfy a super system of linear equations; $H_s x_s^T = (0)$ where $H_s$ is a super mixed row matrix given by $H_s = [H_1 \mid H_2 \mid \ldots \mid H_n]$ where each $H_i$ is a m × $n_i$ matrix with elements from $F_q$, i = 1, 2, …, n, i.e.,

$$\begin{aligned} H_s x_s^T &= [H_1 \mid H_2 \mid \ldots \mid H_n] \left[ x_s^1\ \ x_s^2\ \ldots\ x_s^n \right]^T \\ &= \left[ H_1 \left(x_s^1\right)^T \;\middle|\; H_2 \left(x_s^2\right)^T \;\middle|\; \ldots \;\middle|\; H_2 \left(x_s^n\right)^T \right] \\ &= [|(0) \mid (0) \mid \ldots \mid (0)] \end{aligned}$$

i.e., each $H_i$ is the partity check matrix of the code words $x_s^i$; i = 1, 2, …, n. $H_s = [H_1 \mid H_2 \mid \ldots \mid H_s]$ will be known as the super special parity check super special matrix of the super special row code $C_s$. $C_s$ will also be known as the linear [($n_1$ $n_2$ … $n_n$), ($k_1$ $k_2$ … $k_n$)] or [($n_1$, $k_1$), ($n_2$, $k_2$), …, ($n_n$, $k_n$)] super special row code.

If each of the parity check matrix $H_i$ is of the form $\left(A_i, I_{n_i - k_i}\right)$; i = 1, 2, …, n.

$$H_s = [H_1 \mid H_2 \mid \ldots \mid H_n]$$

$$= \left[ \left(A_1, I_{n_1-k_1}\right) \;\middle|\; \left(A_2, I_{n_2-k_2}\right) \;\middle|\; \ldots \;\middle|\; \left(A_n, I_{n_n-k_n}\right) \right] \quad ---- \quad \text{(I)}$$

$C_s$ is then also called a systematic linear $((n_1\ n_2 \ldots n_n), (k_1\ k_2 \ldots k_n))$ super special code.

If q = 2 then $C_s$ is a super special binary row code; $(k_1 + \ldots + k_n) \mid (n_1 + n_2 + \ldots + n_n)$ is called the super transmission(or super information) rate.

It is important and interesting to note the set $C_s$ of solutions $x_s$ of $H_s x_s^T = (0)$ i.e., known as the super solution space of the super system of equations. Clearly this will form the super special subspace of the super special vector space over $F_q$ of super special dimension $(k_1\ k_2 \ldots k_n)$.

$C_s$ being a super special subspace can be realized to be a group under addition known as the super special group code, where $H_s$ is represented in the form given in equation I will be known as the standard form.

Now we will illustrate this super special row codes by some examples.

***Example 3.1.1:*** Suppose we have a super special binary row code given by the super special parity check matrix $H_S = [H_1 \mid H_2 \mid H_3]$ where

$$H_1 = \begin{bmatrix} 0 & 1 & 1 & 1 & 0 & 0 \\ 1 & 0 & 1 & 0 & 1 & 0 \\ 1 & 1 & 0 & 0 & 0 & 1 \end{bmatrix},$$

$$H_2 = \begin{bmatrix} 0 & 0 & 0 & 1 & 1 & 0 & 0 \\ 0 & 1 & 1 & 0 & 0 & 1 & 0 \\ 1 & 1 & 0 & 1 & 0 & 0 & 1 \end{bmatrix}$$

and

$$H_3 = \begin{bmatrix} 1 & 1 & 0 & 0 & 0 & 1 & 0 & 0 \\ 0 & 0 & 1 & 1 & 0 & 0 & 1 & 0 \\ 1 & 0 & 1 & 0 & 1 & 0 & 0 & 1 \end{bmatrix}$$

i.e., the super row matrix associated with the super special code is given by

$$H_s = \left[\begin{array}{cccccc|ccccccc|cccccccc} 0 & 1 & 1 & 1 & 0 & 0 & 0 & 0 & 0 & 1 & 1 & 0 & 0 & 1 & 1 & 0 & 0 & 0 & 1 & 0 & 0 \\ 1 & 0 & 1 & 0 & 1 & 0 & 0 & 1 & 1 & 0 & 0 & 1 & 0 & 0 & 0 & 1 & 1 & 0 & 0 & 1 & 0 \\ 1 & 1 & 0 & 0 & 0 & 1 & 1 & 1 & 0 & 1 & 0 & 0 & 1 & 1 & 0 & 1 & 0 & 1 & 0 & 0 & 1 \end{array}\right]$$

$$= [(A_1, I_3) \mid (A_2, I_3) \mid (A_3, I_3)] ;$$

i.e., we have given the super special code which is a binary row code. The super special code words are given by

$$x_s = \left[a_1^1\ a_2^1\ a_3^1\ x_4^1\ x_5^1\ x_6^1 \mid a_1^2\ a_2^2\ a_3^2\ a_4^2\ x_5^2\ x_6^2\ x_7^2 \mid a_1^3\ a_2^3\ a_3^3\ a_4^3\ a_5^3\ x_6^3\ x_7^3\ x_8^3\right] = \left[x_s^1 \mid x_s^2 \mid x_s^3\right].$$

$H_s x_s^T$ = (0) gives 3 sets of super linear equations i.e., $H_s x_s^T$ = (0) is read as

$$[H_1 \mid H_2 \mid H_3] \left[x_s^1 \mid x_s^2 \mid x_s^3\right]^T = \left[H_1 \left(x_s^1\right)^T \mid H_2 \left(x_s^2\right)^T \mid H_3 \left(x_s^3\right)^T\right]$$

$$= [(0) \mid (0) \mid (0)];$$

i.e.,

$$H_1 \left(x_s^1\right)^T = (0)$$

gives

$$a_2^1 + a_3^1 + x_4^1 = 0$$

$$a_1^1 + a_3^1 + x_5^1 = 0$$

$$a_1^1 + a_2^1 + x_6^1 = 0 .$$

Therefore

$$x_4^1 = a_2^1 + a_3^1$$

$$x_5^1 = a_1^1 + a_3^1 \text{ and } x_6^1 = a_1^1 + a_2^1 ;$$

i.e., we have {0 0 0 0 0 0, 1 0 0 0 1 1, 1 1 1 0 0 0, 0 1 0 1 0 1, 1 0 1 1 0 1, 0 0 1 1 1 0, 1 1 0 1 1 0, 0 1 1 0 1 1} = $C_s^1$. We define $C_s^1$ to be the subcode of the super special rowcode $C_s$. $H_2 (x_s^2)^T$ = (0) gives

$$a_4^2 + x_5^2 = 0$$
$$a_2^2 + a_3^2 + x_6^2 = 0$$
$$a_1^2 + a_2^2 + a_4^2 + x_7^2 = 0\,.$$

The set of codewords of the super special row subcode is given by

{0 0 0 0 0 0 0, 1 0 0 0 0 0 1, 0 1 0 0 0 1 1, 0 0 1 0 0 1 0,
0 0 0 1 1 0 1, 1 1 0 0 0 1 0, 0 1 1 0 0 0 1, 0 0 1 1 1 1 1,
1 0 1 0 0 1 1, 1 0 0 1 1 0 0, 0 1 0 1 1 1 0, 1 1 1 0 0 0 0,
0 1 1 1 1 0 0, 1 1 0 1 1 1 1, 1 0 1 1 1 1 0, 1 1 1 1 1 0 1} = $C_s^2$.

Now the sublinear equation $H_3\left(x_s^3\right)^T$ = (0) gives

$$a_1^3 + a_2^3 + x_6^3 = 0$$
$$a_3^3 + a_4^3 + x_7^3 = 0$$
$$a_1^3 + a_3^3 + a_5^3 + x_8^3 = 0\,.$$

{0 0 0 0 0 0 0 0, 1 0 0 0 0 1 0 1, 0 1 0 0 0 1 0 1, 0 0 1 0 0 0 1 1,
0 0 0 1 0 01 0, 0 0 0 0 1 0 0 1, 1 1 0 0 0 0 0 1, 0 1 1 0 0 1 1 1,
0 0 1 1 0 0 0 1, 0 0 0 1 1 0 1 1, 1 0 1 0 0 1 1 0, 1 0 0 1 0 1 1 1,
1 0 0 0 1 1 0 1, 0 1 0 1 0 1 1 1, 0 1 0 0 1 1 0 1, 0 0 1 0 1 0 1 1,
1 1 1 0 0 0 1 0, 0 1 1 1 0 1 0 1, 0 0 1 1 1 0 0 0, 1 1 0 1 0 0 1 1,
1 1 0 0 1 0 0 0, 1 0 1 1 0 1 0 0, 1 0 0 1 1 1 1 0, 0 1 1 0 1 1 1 0,
0 1 0 1 1 1 1 1, 1 0 1 0 1 1 1 0, 1 1 1 1 0 0 0 0, 1 1 1 0 1 0 1 1,
0 1 1 1 1 1 0 0, 1 1 0 1 1 0 1 0, 1 0 1 1 1 1 0 1, 1 1 1 1 1 0 0 1} =
$C_s^3$.

$C_s^3$ is a subcode of the super special row code $C_s$. Any element $x_s = \left(x_s^1 \mid x_s^2 \mid x_s^3\right)$ is formed by taking one element from the subcode $C_s^1$, one element from the set of the subcode $C_s^2$ and one from the set $C_s^3$ i.e., an element $x_s$ = [0 0 0 0 1 | 1 0 0 0 0 0 1 | 1 1 1 1 1 0 0 1] which is clearly a super mixed row vector. Thus the number of super special row code words in this example of the super special code $C_s$ is $8 \times 16 \times 32$.

The super transmission rate is 12/21. Thus this code has several advantages which will be enumerated in the last chapter of this book. We give yet another example of super special code in which every super special code word is a super row vector and not a super mixed row vector.

***Example 3.1.2:*** Let $H_s = [H_1 \mid H_2 \mid H_3]$ be the super special parity check super matrix associated with the super special code $C_s$. Here

$$H_1 = \begin{bmatrix} 0 & 1 & 1 & 0 & 1 & 0 & 0 & 0 \\ 1 & 0 & 0 & 1 & 0 & 1 & 0 & 0 \\ 1 & 1 & 1 & 0 & 0 & 0 & 1 & 0 \\ 1 & 0 & 0 & 0 & 0 & 0 & 0 & 1 \end{bmatrix},$$

$$H_2 = \begin{bmatrix} 1 & 1 & 0 & 0 & 1 & 0 & 0 & 0 \\ 1 & 1 & 1 & 0 & 0 & 1 & 0 & 0 \\ 0 & 1 & 1 & 0 & 0 & 0 & 1 & 0 \\ 0 & 1 & 0 & 1 & 0 & 0 & 0 & 1 \end{bmatrix}$$

and

$$H_3 = \begin{bmatrix} 0 & 1 & 1 & 1 & 1 & 0 & 0 & 0 \\ 0 & 0 & 1 & 0 & 0 & 1 & 0 & 0 \\ 0 & 0 & 1 & 1 & 0 & 0 & 1 & 0 \\ 1 & 0 & 1 & 0 & 0 & 0 & 0 & 1 \end{bmatrix};$$

i.e.,

$$H_s = [H_1 \mid H_2 \mid H_3] =$$

$$\left[\begin{array}{cccccccc|cccccccc} 0 & 1 & 1 & 0 & 1 & 0 & 0 & 0 & 1 & 1 & 0 & 0 & 1 & 0 & 0 & 0 \\ 1 & 0 & 0 & 1 & 0 & 1 & 0 & 0 & 1 & 1 & 1 & 0 & 0 & 1 & 0 & 0 \\ 1 & 1 & 1 & 0 & 0 & 0 & 1 & 0 & 0 & 1 & 1 & 0 & 0 & 0 & 1 & 0 \\ 1 & 0 & 0 & 0 & 0 & 0 & 0 & 1 & 0 & 1 & 0 & 1 & 0 & 0 & 0 & 1 \end{array}\right.$$

$$\begin{bmatrix} 0 & 1 & 1 & 1 & 1 & 0 & 0 & 0 \\ 0 & 0 & 1 & 0 & 0 & 1 & 0 & 0 \\ 0 & 0 & 1 & 1 & 0 & 0 & 1 & 0 \\ 1 & 0 & 1 & 0 & 0 & 0 & 0 & 1 \end{bmatrix}$$

is the super special parity check matrix of the super special code $C_s$. Now the super special system of equations is given by

$$H_s x_s^T = (0) \text{ i.e.,}$$

$$\left[H_1 \,|H_2\, |H_3\right]\left[x_s^1 \,|x_s^2\, |x_s^3\right]^T = [(0) \mid (0) \mid (0)]$$

$$= \left[H_1\left(x_s^1\right)^T \middle| H_2\left(x_s^2\right)^T \middle| H_3\left(x_s^3\right)^T\right].$$

We call the linear equations given by $H_i\left(x_s^i\right)^T = (0)$ to be subequations of the super linear equations. Now the sublinear equations given by $H_1\left(x_s^1\right)^T = (0)$ is

$$a_2^1 + a_3^1 + x_1^1 = 0$$
$$a_1^1 + a_4^1 + x_2^1 = 0$$
$$a_1^1 + a_2^1 + a_3^1 + x_3^1 = 0$$
$$a_1^1 + x_4^1 = 0\,;$$

where $a_1^1\ a_2^1\ a_3^1\ a_4^1$ is the set of message symbols and $x_1^1\ x_2^1\ x_3^1\ x_4^1$ the check symbols and the check equations using $H_1$ is given above. The subcode associated with super code $C_s$ is

$C_s^1$ = {0 0 0 0 0 0 0 0, 1 0 0 0 0 1 1 1, 0 1 0 0 1 0 1 0,
0 0 1 0 1 0 1 0, 0 0 0 1 0 1 0 0, 1 1 0 0 1 1 0 1,
0 1 1 0 0 1 0 0, 0 0 1 1 1 1 0 0, 1 0 0 1 0 0 1 1,
1 0 1 0 1 1 0 1, 0 1 0 1 1 1 1 0, 1 1 1 0 0 1 1 1,
0 1 1 1 0 1 0 0, 1 1 0 1 1 0 1 1, 1 0 1 1 1 0 0 1,
1 1 1 1 0 0 1 1} .

Now we use the sublinear equation of the super linear equation got from $H_s x_s^T$ = (0), we get $H_2\left(x_s^2\right)^T$ = (0). This gives

$$a_1^2 + a_2^2 + x_1^2 = 0$$
$$a_1^2 + a_2^2 + a_3^2 + x_2^2 = 0$$
$$a_2^2 + a_3^2 + x_3^2 = 0$$
$$a_2^2 + a_4^2 + x_4^2 = 0 .$$

The super special subcode $C_s^2$ associated with the above set of equations is given by :

$C_s^2$ = {0 0 0 0 0 0 0 0, 1 0 0 0 1 1 0 0, 0 1 0 0 1 1 1 1,
0 0 1 0 0 1 1 0, 0 0 0 1 0 0 0 1, 1 1 0 0 0 0 1 1,
0 1 1 0 1 0 0 1, 0 0 1 1 0 1 1 1, 1 0 1 0 1 0 0 1,
1 0 0 1 1 1 0 1, 0 1 0 1 1 1 1 0, 1 1 1 0 0 1 0 1,
0 1 1 1 1 0 0 0, 1 1 0 1 0 0 1 1, 1 0 1 1 1 0 1 1,
1 1 1 1 0 1 0 0}.

Now using the sublinear equation $H_3\left(x_s^3\right)^T$ = (0) we get

$$a_2^3 + a_3^3 + a_4^3 + x_1^3 = 0$$
$$a_3^2 + x_2^3 = 0$$
$$a_3^3 + a_4^3 + x_3^3 = 0$$
$$a_1^3 + a_3^3 + x_4^3 = 0 .$$

$C_s^3$ = {0 0 0 0 0 0 0 0, 1 0 0 0 0 0 0 1, 0 1 0 0 1 0 0 0,
0 0 1 0 1 1 1 1, 0 0 0 1 1 0 1 0, 1 1 0 0 1 0 0 1,
0 1 1 0 0 1 1 1, 0 0 1 1 0 1 0 1, 1 0 0 1 1 0 1 1,
0 1 0 1 0 0 1 0, 1 0 1 0 1 1 1 0, 1 1 1 0 0 1 1 0,
0 1 1 1 1 1 0 1, 1 1 0 1 0 1 1 1, 1 0 1 1 0 1 0 0,
1 1 1 1 1 1 0 0}.

Thus the super special code word of the super special code $C_s$ will be $C_s = \left[C_s^1 \mid C_s^2 \mid C_s^3\right]$ i.e., it is formed by taking one code word from each one of the $C_s^i$; i = 1, 2, 3. Thus if $x_s \in C_s$ then any $x_s$ = [1 1 1 1 0 0 1 1 | 1 1 1 1 0 1 0 0| 1 1 1 1 1 1 0 0].

So we can realize the super special row code to be codes in which each and every subcode $C_s^i$ of $C_s$ have the same number of check symbols.

Now we proceed on to define the notion of super special row repetition code.

**DEFINITION 3.1.2:** *Let $C_s = \left[C_s^1 \mid C_s^2 \mid \ldots \mid C_s^n\right]$ be a super special row code in which each of the $C_s^i$ is a repetition code, i = 1, 2, …, n, then we define $C_s$ to be a super special repetition row code. Here if $H_s = [H_1|H_2| \ldots|H_n]$ is the super special parity check matrix of $C_s$, then each $H_i$ is a $t – 1 \times t$ matrix that is we have*

$$H_1 = H_2 = \ldots = H_n = \begin{bmatrix} 1 & 1 & 0 & \cdots & 0 \\ 1 & 0 & 1 & \cdots & 0 \\ \vdots & \vdots & \vdots & & \vdots \\ 1 & 0 & 0 & \ldots & 1 \end{bmatrix}_{t-1 \times t}$$

*is the parity check matrix. The super special code words associated with $C_s$ are just super row vectors only and not super mixed row vectors. The number of super special code words in $C_s$ is $2^n$.*

We illustrate a super special row repetition code by the following example.

***Example 3.1.3:*** Let $C_s = \left[C_s^1 \mid C_s^2 \mid C_s^3 \mid C_s^4\right]$ be a super row repetition code with associated super special row matrix $H_s$ = $[H_1 \mid H_2 \mid H_3 \mid H_4]$

$$= \left[\begin{array}{cccccc|cccccc|cccccc|cccccc}
1 & 1 & 0 & 0 & 0 & 0 & 1 & 1 & 0 & 0 & 0 & 0 & 1 & 1 & 0 & 0 & 0 & 0 & 1 & 1 & 0 & 0 & 0 & 0 \\
1 & 0 & 1 & 0 & 0 & 0 & 1 & 0 & 1 & 0 & 0 & 0 & 1 & 0 & 1 & 0 & 0 & 0 & 1 & 0 & 1 & 0 & 0 & 0 \\
1 & 0 & 0 & 1 & 0 & 0 & 1 & 0 & 0 & 1 & 0 & 0 & 1 & 0 & 0 & 1 & 0 & 0 & 1 & 0 & 0 & 1 & 0 & 0 \\
1 & 0 & 0 & 0 & 1 & 0 & 1 & 0 & 0 & 0 & 1 & 0 & 1 & 0 & 0 & 0 & 1 & 0 & 1 & 0 & 0 & 0 & 1 & 0 \\
1 & 0 & 0 & 0 & 0 & 1 & 1 & 0 & 0 & 0 & 0 & 1 & 1 & 0 & 0 & 0 & 0 & 1 & 1 & 0 & 0 & 0 & 0 & 1
\end{array}\right].$$

Thus

$C_s$ = {[0 0 0 0 0 0 | 0 0 0 0 0 0 | 0 0 0 0 0 0 | 0 0 0 0 0 0],
[1 1 1 1 1 1 | 1 1 1 1 1 1 | 1 1 1 1 1 1 | 1 1 1 1 1 1],
[0 0 0 0 0 0 | 1 1 1 1 1 1 | 0 0 0 0 0 0 | 1 1 1 1 1 1],
[0 0 0 0 0 0 | 1 1 1 1 1 1 | 0 0 0 0 0 0 | 0 0 0 0 0 0],
[0 0 0 0 0 0 | 0 0 0 0 0 0 | 0 0 0 0 0 0 | 1 1 1 1 1 1],
[1 1 1 1 1 1 | 0 0 0 0 0 0 | 0 0 0 0 0 0 | 0 0 0 0 0 0],
[0 0 0 0 0 0 | 0 0 0 0 0 0 | 1 1 1 1 1 1 | 0 0 0 0 0 0],
[1 1 1 1 1 1 | 1 1 1 1 1 1 | 0 0 0 0 0 0 | 0 0 0 0 0 0],
[1 1 1 1 1 1 | 0 0 0 0 0 0 | 1 1 1 1 1 1 | 0 0 0 0 0 0],
[1 1 1 1 1 1 | 0 0 0 0 0 0 | 0 0 0 0 0 0 | 1 1 1 1 1 1],
[0 0 0 0 0 0 | 1 1 1 1 1 1 | 1 1 1 1 1 1 | 0 0 0 0 0 0],
[0 0 0 0 0 0 | 0 0 0 0 0 0 | 1 1 1 1 1 1 | 1 1 1 1 1 1],
[1 1 1 1 1 1 | 0 0 0 0 0 0 | 1 1 1 1 1 1 | 0 0 0 0 0 0],
[1 1 1 1 1 1 | 1 1 1 1 1 1 | 0 0 0 0 0 0 | 1 1 1 1 1 1],
[1 1 1 1 1 1 | 1 1 1 1 1 1 | 1 1 1 1 1 1 | 0 0 0 0 0 0],
[0 0 0 0 0 0 | 1 1 1 1 1 1 | 1 1 1 1 1 1 | 1 1 1 1 1 1]}.

Clearly $|C_s| = 2^4 = 16$.

Now having seen an example of a super special repetition row code we proceed on to define the super special parity check row code. We have two types of super special row parity check codes.

**DEFINITION 3.1.3:** *Let $C_s$ be a super special parity check mixed row code i.e.,* $C_s = \left[C_s^1 \mid C_s^2 \mid \ldots \mid C_s^n\right]$ *where $C_s$ is obtained using the super special mixed row matrix $H_s = [H_1 \mid H_2 \mid \ldots \mid H_n]$ where each $H_i$ is a unit row vector having $t_i$ number of elements i.e.,*

$$H_s = \left[\underbrace{1\ \ 1\ \ \cdots\ \ 1}_{t_1 \text{ times}} \Bigg| \underbrace{1\ \ 1\ \ \cdots\ \ 1}_{t_2 \text{ times}} \Bigg| \cdots \Bigg| \underbrace{1\ \ 1\ \ \cdots\ \ 1}_{t_n \text{ times}}\right]$$

*where at least one $t_i \neq t_j$ for $i \neq j$. Any super special code word in $C_s$ would be of the form*

$$x_s = \left[x_1^1\ x_2^1 \ldots x_{t_1}^1 \;\middle|\; x_1^2\, x_2^2 \ldots x_{t_2}^2 \;\middle|\; \ldots \middle|\, x_1^n\, x_2^n \ldots x_{t_n}^n\right] = \left[x_s^1 \middle| x_s^2 \middle| \ldots \middle| x_s^n\right]$$

*with $H_s x_s^T$ = (0); i.e., each $x_s^i$ would contain only even number of ones and the rest are zeros.*

$C_s = [C_1 \mid C_2 \mid \ldots \mid C_n]$ is defined to be super special parity check row code. $C_s$ is obtained from the parity check row matrix / vector $H_s = [H_1 \mid H_2 \mid \ldots \mid H_n]$ where $H_1 = H_2 = \ldots = H_n$ $= \left[\underbrace{1\ \ 1\ \ \cdots\ \ 1}_{\text{m times}}\right]$. Here a super special codeword in $C_s$ would be a super row vector of the form $\left[x_s^1 \mid x_s^2 \mid \ldots \mid x_s^n\right]$ with each $x_s^i = \left[x_1^i\ x_2^i\ \ldots\ x_m^i\right]$ where only even number of $x_j^i$ are ones and the rest zero, $1 \leq j \leq m$ and $i = 1, 2, \ldots, n$.

Now we will illustrate the two types of super special parity check (mixed) row codes.

***Example 3.1.4:*** Let $C_s = \left[C_s^1 \middle| C_s^2 \middle| C_s^3\right]$ be a super special parity check code having $H_s = [H_1 \mid H_2 \mid H_3]$ = [1 1 1 | 1 1 1 | 1 1 1] to be the super special parity check matrix associated with it.

$C_s$ = {[0 0 0 | 0 0 0 | 0 0 0], [0 0 0 | 0 0 0 | 1 1 0],

[0 0 0 | 0 0 0 | 1 0 1], [0 0 0 | 0 0 0 | 0 1 1], [0 0 0 | 1 1 0 | 0 0 0],
[0 0 0 | 1 1 0 | 1 1 0], [0 0 0 | 1 1 0 | 1 0 1], [0 0 0 | 1 1 0 | 0 1 1],
[0 0 0 | 0 1 1 | 0 0 0], [0 0 0 | 0 1 1 | 1 0 1], [0 0 0 | 0 1 1 | 1 1 0],
[0 0 0 | 0 1 1 | 0 1 1], [0 0 0 | 1 0 1 | 1 1 0], [0 0 0 | 1 0 1 | 0 0 0],
[0 0 0 | 1 0 1 | 0 1 1], [0 0 0 | 1 0 1 | 1 0 1], [1 1 0 | 0 0 0 | 0 0 0],
[1 1 0 | 0 0 0 | 1 0 1], [1 1 0 | 0 0 0 | 1 1 0], [1 1 0 | 0 0 0 | 0 1 1],
[1 1 0 | 1 1 0 | 0 0 0], [1 1 0 | 1 1 0 | 0 1 1], [1 1 0 | 1 1 0 | 1 0 1],
[1 1 0 | 1 1 0 | 1 1 0], [1 1 0 | 1 0 1 | 0 0 0], [1 1 0 | 1 0 1 | 1 0 1],
[1 1 0 | 1 0 1 | 0 1 1], [1 1 0 | 1 0 1 | 1 1 0], [1 1 0 | 0 1 1 | 0 0 0],
[1 1 0 | 0 1 1 | 1 1 0], [1 1 0 | 0 1 1 | 1 0 1], [1 1 0 | 0 1 1 | 0 1 1],
[0 1 1 | 0 0 0 | 0 0 0], [0 1 1 | 0 0 0 | 0 1 1], [0 1 1 | 0 0 0 | 1 1 0],
[0 1 1 | 0 0 0 | 1 0 1], [0 1 1 | 1 1 0 | 0 0 0], [0 1 1 | 1 1 0 | 1 1 0],
[0 1 1 | 1 1 0 | 0 1 1], [0 1 1 | 1 1 0 | 1 0 1], [0 1 1 | 0 1 1 | 0 0 0],
[0 1 1 | 0 1 1 | 1 0 1], [0 1 1 | 0 1 1 | 1 1 0], [0 1 1 | 0 1 1 | 0 1 1],
[0 1 1 | 1 0 1 | 0 0 0], [0 1 1 | 1 0 1 | 0 1 1], [0 1 1 | 1 0 1 | 1 0 1],
[0 1 1 | 1 0 1 | 1 1 0], [1 0 1 | 0 0 0 | 0 0 0], [1 0 1 | 0 0 0 | 0 1 1],
[1 0 1 | 0 0 0 | 1 1 0], [1 0 1 | 0 0 0 | 1 0 1], [1 0 1 | 0 1 1 | 0 0 0],
[1 0 1 | 0 1 1 | 1 1 0], [1 0 1 | 0 1 1 | 1 0 1], [1 0 1 | 0 1 1 | 0 1 1],
[1 0 1 | 1 0 1 | 0 0 0], [1 0 1 | 1 0 1 | 1 0 1], [1 0 1 | 1 0 1 | 0 1 1],
[1 0 1 | 1 0 1 | 1 1 0], [1 0 1 | 1 1 0 | 0 0 0], [1 0 1 | 1 1 0 | 1 1 0],
[1 0 1 | 1 1 0 | 0 1 1],[1 0 1 | 1 1 0 | 1 0 1]}.

Thus $|C_s| = 4^3 = 64$. We see in every super row vector the number of non zero ones is even.

Next we give an example of a super special parity check row code $C_s$.

***Example 3.1.5:*** Let $C_s = \left[C_s^1 \mid C_s^2\right]$ be a super special parity check mixed code with the associated super special parity check mixed row vector $H_s = [H_1 \mid H_2] = [1\ 1\ 1\ 1 \mid 1\ 1\ 1]$. The super special codewords given by $H_s$ is

$C_s$ = {[0 0 0 0 | 0 0 0], [0 0 0 0 | 1 1 0], [0 0 0 0 | 1 0 1],
[0 0 0 0 | 0 1 1], [1 0 1 0 | 0 0 0], [1 0 1 0 | 0 1 1],
[1 0 1 0 | 1 0 1], [1 0 1 0 | 1 1 0], [1 0 0 1 | 0 0 0],
[1 0 0 1 | 0 1 1], [1 0 0 1 | 1 1 0], [1 0 0 1 | 1 0 1],
[0 1 0 1 | 0 1 1], [0 1 0 1 | 0 0 0], [0 1 0 1 | 1 0 1],
[0 1 0 1 | 1 1 0], [1 1 0 0 | 0 0 0], [1 1 0 0 | 0 1 1],
[1 1 0 0 | 1 0 1], [1 1 0 0 | 1 1 0], [0 1 1 0 | 0 0 0],

[0 1 1 0 | 1 1 0], [0 1 1 0 | 1 0 1], [0 1 1 0 | 0 1 1],
[0 0 1 1 | 0 0 0], [0 0 1 1 | 1 0 1], [0 0 1 1 | 1 1 0],
[0 0 1 1 | 0 1 1], [1 1 1 1 | 0 0 0], [1 1 1 1 | 1 1 0],
[1 1 1 1 | 0 1 1], [1 1 1 1 | 1 0 1]}.

Clearly $|C_s| = 8 \times 4 = 32$.

Now having seen examples of the two types of super special parity check codes, we now proceed on to define super special Hamming distance, super special Hamming weight and super special errors in super special codes $C_s$.

**DEFINITION 3.1.4:** *Let* $C_s = \left[C_s^1 \mid C_s^2 \mid \ldots \mid C_s^n\right]$ *be a super special row code. Suppose* $x_s = \left[x_s^1 \mid x_s^2 \mid \ldots \mid x_s^n\right]$ *is a transmitted super code word and* $y_s = \left[y_s^1 \mid y_s^2 \mid \ldots \mid y_s^n\right]$ *is the received supercode word then* $e_s = y_s - x_s = \left[y_s^1 - x_s^1 \mid y_s^2 - x_s^2 \mid \ldots \mid y_s^n - x_s^n\right] = \left[e_s^1 \mid e_s^2 \mid \ldots \mid e_s^n\right]$ *is called the super error word or the super error vector.*

We first illustrate how the super error is determined.

***Example 3.1.6:*** Let $C_s = \left[C_s^1 \mid C_s^2 \mid C_s^3 \mid C_s^4\right]$ be a super special code with associated super parity check row matrix $H_s = [H_1 \mid H_2 \mid H_3 \mid H_4]$

$$= \left[\begin{array}{ccccccc|cccccccc} 1&0&0&1&0&0&0&1&0&1&0&1&0&0&0 \\ 1&1&0&0&1&0&0&0&1&0&1&0&1&0&0 \\ 0&1&0&0&0&1&0&1&1&0&0&0&0&1&0 \\ 0&0&1&0&0&0&1&0&1&1&0&0&0&0&1 \end{array}\right|$$

$$\left|\begin{array}{cccccc|ccccccccc} 1&0&1&0&0&0&1&0&1&0&0&1&0&0&0 \\ 1&1&0&1&0&0&1&1&0&1&0&0&1&0&0 \\ 0&1&0&0&1&0&0&1&1&0&1&0&0&1&0 \\ 1&1&0&0&0&1&0&0&1&1&0&0&0&0&1 \end{array}\right].$$

Let $x_s = \left[x_s^1 \mid x_s^2 \mid x_s^3 \mid x_s^4\right]$ = [1 1 0 1 0 1 0 | 1 1 1 1 0 0 0 0 | 1 1 1 0 1 0 | 1 1 1 1 0 0 1 0 0] $\in C_s$ be the sent super special code. Suppose $y_s = \left[y_s^1 \mid y_s^2 \mid y_s^3 \mid y_s^4\right]$ = [1 0 1 0 1 1 0 | 1 1 1 1 1 1 0 0| 1 1 0 1 0 1| 1 1 1 0 1 0 1 0 0] be the received super code word. The super error vector is given by

$$\begin{aligned} y_s - x_s &= \left[y_s^1 \mid y_s^2 \mid y_s^3 \mid y_s^4\right] - \left[x_s^1 \mid x_s^2 \mid x_s^3 \mid x_s^4\right] \\ &= \left[y_s^1 - x_s^1 \middle| y_s^2 - x_s^2 \middle| y_s^3 - x_s^3 \middle| y_s^4 - x_s^4\right] \\ &= [0\ 1\ 1\ 1\ 1\ 0\ 0 \mid 0\ 0\ 0\ 0\ 1\ 1\ 0\ 0 \mid 0\ 0\ 1\ 1\ 1\ 1 \mid 0\ 0\ 0\ 1\ 1\ 0\ 0\ 0\ 0] \\ &= \left[e_s^1 \mid e_s^2 \mid e_s^3 \mid e_s^4\right] \\ &= e_s. \end{aligned}$$

Clearly $y_s + e_s = x_s$.

**DEFINITION 3.1.5:** *The super Hamming distance $d_s(x_s, y_s)$ between two super row vectors of the super special vector space $V_s$, where $x_s = \left[x_s^1 \mid x_s^2 \mid \ldots \mid x_s^n\right]$ and $y_s = \left[y_s^1 \mid y_s^2 \mid \ldots \mid y_s^n\right]$ is the number of coordinates in which $x_s^i$ and $y_s^i$ differ for i = 1, 2, …, n. The super Hamming weight $w_s(x_s)$ of the super vector $x_s = \left[x_s^1 \mid x_s^2 \mid \ldots \mid x_s^n\right]$ in $V_s$ is the number of non zero coordinates in each $x_s^i$; i = 1, 2, …, n. In short $w_s(x_s) = d(x_s, (0))$.*

As in case of usual linear codes we define super minimum distance $d_{min}^s$ of a super special linear row code $C_s$ as

$$d_{min}^s = \min_{\substack{u_s, \nu_s \in C_s \\ u_s \neq \nu_s}} d_s(u_s, \nu_s),$$

$$d_s(u_s, \nu_s) = d_s(u_s - \nu_s, (0)) = w_s(u_s - \nu_s).$$

Thus the super minimum distance of $C_s$ is equal to the least super weights of all non zero super special code words.

Now the value of

$$d^s_{min} = \min_{\substack{u_s, \nu_s \in C_s \\ u_s \neq \nu_s}} d_s\,(u_s, \nu_s)$$

$$= \min_{\substack{u_s, \nu_s \in C_s \\ u_s \neq \nu_s}} d_s\left(\left[u_s^1 \middle| u_s^2 \middle| \ldots \middle| u_s^n\right], \left[v_s^1 \middle| v_s^2 \middle| \ldots \middle| v_s^n\right]\right)$$

$$= \min\ \left[d\left(u_s^1,\ \nu_s^1\right) + d\left(u_s^2,\ \nu_s^2\right) + \ldots + d\left(u_s^n,\ \nu_s^n\right)\right].$$

Now $d^s_{min}$ of the super special row code given in example 3.1.2 is 7 verified using the fact in $C_s = \left[C_s^1 \mid C_s^2 \mid C_s^3\right]$; $d^s_{min} C_s^1 = 3$, $d^s_{min} C_s^2 = 2$ and $d^s_{min} C_s^3 = 2$. Hence $d^s_{min}\left(C_s\right) = 3 + 2 + 2 = 7$. So we will denote $d^s_{min} = \min\ d_s\left(u_s, \nu_s\right)$ by $d^s_{min}\left(C_s\right)$, $u_s, \nu_s \in C_s$, $u_s \neq \nu_s$.

$$d^S_{min}\left[C_s^1 \mid C_s^2 \mid \ldots \mid C_s^n\right] = d^S_{min}\left[\left(C_s^1\right) + \left(C_s^2\right) + \ldots + \left(C_s^n\right)\right]$$

$$= \min_{\substack{x_s^1, y_s^1 \in C_s^1 \\ x_s^1 \neq y_s^1}} d\left(x_s^1,\ y_s^1\right) + \min_{\substack{x_s^2, y_s^2 \in C_s^2 \\ x_s^2 \neq y_s^2}} d\left(x_s^2,\ y_s^2\right) + \ldots + \min_{\substack{x_s^n, y_s^n \in C_s^n \\ x_s^n \neq y_s^n}} d\left(x_s^n,\ y_s^n\right).$$

Now we proceed on to define the dual of a super special row code.

**DEFINITION 3.1.6:** *Let* $C_s = \left[C_s^1 \mid C_s^2 \mid \ldots \mid C_s^n\right]$ *be a super special row [($n_1$, …, $n_n$), ($k_1$, …, $k_n$)] binary code. The super special dual row code of $C_s$ denoted by*

$$C_s^{\perp} = \left[\left(C_s^1\right)^{\perp} \,\middle|\, \left(C_s^2\right)^{\perp} \,\middle|\, \ldots \,\middle|\, \left(C_s^n\right)^{\perp}\right]$$

*where* $\left(C_s^i\right)^{\perp} = \{\, u_s^i \mid u_s^i \cdot v_s^i = 0$ *for all* $v_s^i \in C_s^i \,\}$, *i = 1, 2, …, n.*

*Since in $C_s$ we have $n_1 - k_1 = n_2 - k_2 = \ldots = n_n - k_n$ i.e., the number of check symbols of each and every code in $C_s^i$ is the same for i = 1, 2, …, n. Thus we see $n = 2k_i$ alone can give us a dual, in all other cases we will have problem with the compatibility for the simple reason the dual code of $C_s^i$ being the orthogonal complement will have $n_i - k_i$ to be the dimension, where as $C_s^i$ will be of dimension $k_i$, i = 1, 2, …, n. Hence we* *can say the super special dual code would be defined if and only if $n_i = 2k_i$ and such that $n_1 = n_2 = \ldots = n_n$.*

We can define the new notion of super special syndrome to super code words of a super special row code which is analogous to syndrome of the usual codes.

**DEFINITION 3.1.7:** *Let $C_s$ be a super special row code. Let $H_s$ be the associated super special parity check matrix of $C_s$ the super special syndrome of any element $y_s \in V_s$ where $C_s$ is a super special subspace of the super special vector space $V_s$ is given by $S(y_s) = H_s y_s^T$. $S(y_s) = (0)$ if and only if $y_s \in C_s$.*

*Thus this gives us a condition to find out whether the received super code word is a right message or not. Suppose $y_s$ is the received super special code word, we find $S(y_s) = H_s y_s^T$; if $S(y_s) = (0)$ then we accept $y_s$ as the correct message if $S(y_s) = H_s y_s^T \neq (0)$ then we can declare the received word has error.*

We can find the correct word by the following method. Before we give this method we illustrate how the super special syndrome is calculated.

***Example 3.1.7:*** Let $C_s = \left[C_s^1 \mid C_s^2 \mid C_s^3 \mid C_s^4\right]$ be a super special row code. Let $H_s = [H_1 \mid H_2 \mid H_3 \mid H_4]$ be the super special parity check matrix of $C_s$.

Let

$$H_s = \left[\begin{array}{ccccccc|cccccc} 1 & 0 & 0 & 0 & 1 & 0 & 0 & 0 & 1 & 0 & 1 & 0 & 0 \\ 1 & 0 & 0 & 1 & 0 & 1 & 0 & 1 & 0 & 1 & 0 & 1 & 0 \\ 0 & 1 & 1 & 0 & 0 & 0 & 1 & 0 & 1 & 1 & 0 & 0 & 1 \end{array}\right|$$

$$\left.\begin{array}{ccccc|cccccccc} 0 & 1 & 1 & 0 & 0 & 0 & 0 & 1 & 1 & 0 & 1 & 0 & 0 \\ 1 & 0 & 0 & 1 & 0 & 1 & 0 & 1 & 0 & 1 & 0 & 1 & 0 \\ 1 & 1 & 0 & 0 & 1 & 0 & 1 & 0 & 0 & 0 & 0 & 0 & 1 \end{array}\right].$$

Suppose

$$x_s = \left[x_s^1 \mid x_s^2 \mid x_s^3 \mid x_s^4\right]$$
$$= [1\ 0\ 0\ 0\ 1\ 1\ 0 \mid 0\ 1\ 0\ 1\ 0\ 1 \mid 1\ 1\ 1\ 1\ 0 \mid 0\ 0\ 0\ 0\ 1\ 0\ 1\ 0] \in C_s.$$

Now the super special syndrome of $x_s$ is given by

$$\begin{aligned} S(x_s) &= H_s x_s^T \\ &= [H_1 \mid H_2 \mid H_3 \mid H_4] \left[x_s^1 \mid x_s^2 \mid x_s^3 \mid x_s^4\right]^T \\ &= \left[H_1 \left(x_s^1\right)^T \mid H_2 \left(x_s^2\right)^T \mid H_3 \left(x_s^3\right)^T \mid H_4 \left(x_s^4\right)^T\right] \\ &= [0\ 0\ 0 \mid 0\ 0\ 0 \mid 0\ 0\ 0 \mid 0\ 0\ 0]. \end{aligned}$$

Let $y_s = [1\ 1\ 1\ 0\ 0\ 1\ 1 \mid 0\ 1\ 1\ 1\ 1\ 0 \mid 1\ 0\ 1\ 0\ 1 \mid 1\ 1\ 1\ 0\ 0\ 1\ 1\ 1] \in V_s$ where $C_s$ is a proper super special subspace of the super special vector space $V_s$.

Now

$$\begin{aligned} S(y_s) &= H_s y_s^T \\ &= [H_1 \mid H_2 \mid H_3 \mid H_4] \left[y_s^1 \mid y_s^2 \mid y_s^3 \mid y_s^4\right]^T \\ &= \left[H_1 \left(y_s^1\right)^T \mid H_2 \left(y_s^2\right)^T \mid H_3 \left(y_s^3\right)^T \mid H_4 \left(y_s^4\right)^T\right] \\ &= [1\ 0\ 1 \mid 0\ 0\ 0 \mid 1\ 1\ 0 \mid 0\ 1\ 0] \\ &\neq [0\ 0\ 0 \mid 0\ 0\ 0 \mid 0\ 0\ 0 \mid 0\ 0\ 0]. \end{aligned}$$

Thus $y_s \notin C_s$.

Now we have to shown how to find whether the received super special code word is correct or otherwise. It is important to note that what ever be the super special code word $x_s \in C_s$ (i.e., it may be a super special mixed row vector or not) but the syndrome $S(x_s) = H_s x_s^T$ is always a super special row vector which is not mixed and each row vector is of length equal to the number of rows of $H_s$.

Now we know that every super special row code $C_s$ is a subgroup of the super special vector space $V_s$ over $Z_2 = \{0, 1\}$. Now we can for any $x_s \in V_s$ define super special cosets as

$$x_s + C_s = \{x_s + c_s \mid c_s \in C_s\}.$$

Thus

$$\begin{aligned} V_s &= \{Z_2 \times Z_2 \times \ldots \times Z_2 \mid Z_2 \times \ldots \times Z_2 \mid \ldots \mid Z_2 \times Z_2 \times \ldots \times Z_2\} \\ &= C_s \cup [x_s^1 + C_s] \cup \ldots \cup [x_s^t + C_s] \end{aligned}$$

where $C_s = \left[C_s^1 \mid C_s^2 \mid \ldots \mid C_s^n\right]$

and

$$x_s = \left[x_1^1 \ldots x_{n_1}^1 \middle| x_1^2 \ldots x_{n_2}^2 \middle| \ldots \middle| x_1^n \ldots x_{n_n}^n\right]$$

$$= \left[x_s^1 \mid x_s^2 \mid \ldots \mid x_s^n\right]$$

and

$$x_s + C_s = \left[x_s^1 + C_s^1 \mid x_s^2 + C_s^2 \mid \ldots \mid x_s^n + C_s^n\right].$$

Now we can find the coset leader of every $x_s^i + C_s^i$ as in case of usual codes described in chapter one of this book. Now if

$$y_s = \left[y_s^1 \mid y_s^2 \mid \ldots \mid y_s^n\right]$$

is the received message and $\left[e_s^1 + (0) \middle| e_s^2 + (0) \middle| \ldots \middle| e_s^n + (0)\right]$ is a special super coset leaders then using the relation $y_s - e_s$ we get $y_s - e_s$ to be super special corrected code word. It is interesting to note that each $e_s^i$ + (0) has a stipulated number of coset leaders depending on $n_i$, i = 1, 2, …, n.

We will illustrate this by the following example.

***Example 3.1.8:*** Let $C_s = \left[C_s^1 \middle| C_s^2\right]$ be a super special row code. Suppose $H_s$ = [$H_1$ | $H_2$] be the super special row matrix associated with $C_s$. Let

$$H_s = [H_1 \mid H_2] = \left[\begin{array}{cccc|ccccc} 1 & 0 & 1 & 0 & 1 & 0 & 1 & 1 & 0 \\ 1 & 1 & 0 & 1 & 0 & 1 & 1 & 0 & 1 \end{array}\right].$$

Now $C_s = \left[C_s^1 \middle| C_s^2\right]$ with $C_s^1$ = {(0 0 0 0), (1 0 1 1), (0 1 0 1), (1 1 1 0)}and $C_s^2$ = {( 0 0 0 0 0), (1 0 0 1 0), (0 1 0 0 1), (0 0 1 1 1), (1 1 0 1 1), (0 1 1 1 0), (1 0 1 0 1), (1 1 1 0 0)}.

$C_s$ = {[0 0 0 0 | 0 0 0 0 0], [1 0 1 1 | 0 0 0 0 0],
[0 1 0 1 | 0 0 0 0 0], [1 1 1 0 | 0 0 0 0 0], [0 0 0 0 | 1 0 0 1 0],
[0 0 0 0 | 0 1 0 0 1], [1 1 1 0 | 0 1 0 0 1], [0 1 0 1 | 0 1 0 0 1],
[1 0 1 1 | 1 0 0 0 1], [0 0 0 0 | 0 0 1 1 1], [1 1 1 0 | 0 0 1 1 1],
[0 1 0 1 | 0 0 1 1 1], [1 0 1 1 | 0 0 1 1 1], [0 0 0 0 | 1 1 0 1 1],

[1 1 1 0 | 1 1 0 1 1], [1 0 1 1 | 1 1 0 1 1], [0 1 0 1 | 1 1 0 1 1],
[0 0 0 0 | 0 1 1 1 0], [1 0 1 1 | 0 1 1 1 0], [0 1 0 1 | 0 1 1 1 0],
[0 0 0 0 | 1 0 1 0 1], [1 0 1 1 | 1 0 1 0 1], [0 1 0 1 | 1 0 1 0 1],
[0 0 0 0 | 1 1 1 0 0], [1 1 1 0 | 1 1 1 0 0], [1 0 1 1 | 1 1 1 0 0],
[0 1 0 1 | 1 1 1 0 0] and so on}.

Clearly $|C_s| = 32$ . Now the coset table of $C_s^1$ is given by

| Message | code | | words |
|---|---|---|---|
| 0 0 | 1 0 | 0 1 | 1 1 |
| 0 0 0 0 | 1 0 1 1 | 0 1 0 1 | 1 1 1 0 |
| Other cosets | | | |
| 1 0 0 0 | 0 0 1 1 | 1 1 0 1 | 0 1 1 0 |
| 0 1 0 0 | 1 1 1 1 | 0 0 0 1 | 1 0 1 0 . |
| 0 0 1 0 | 1 0 0 1 | 0 1 1 1 | 1 1 0 0 |
| coset leaders | | | |

Now the coset table of $C_s^2$ is given by

| | | |
|---|---|---|
| message | 0 0 0 | 1 0 0 |
| codewords | 0 0 0 0 0 | 1 0 0 1 0 |
| other cosets | 1 0 0 0 0 | 0 0 0 1 0 |
| | 0 1 0 0 0 | 1 1 0 1 0 |
| | 0 0 1 0 0 | 1 0 1 1 0 |
| message | 0 1 0 | 0 0 1 |
| codewords | 0 1 0 0 1 | 0 0 1 1 1 |
| other cosets | 1 1 0 0 1 | 1 0 1 1 1 |
| | 1 0 0 0 1 | 0 1 1 1 1 |
| | 0 1 1 0 1 | 0 0 0 1 1 |
| message | 1 1 0 | 0 1 1 |
| codewords | 1 1 0 1 1 | 0 1 1 1 0 |
| other cosets | 0 1 0 1 1 | 1 1 1 1 0 |
| | 1 0 0 1 1 | 0 0 1 1 0 |
| | 1 1 1 1 1 | 0 1 0 1 0 |
| message | 1 0 1 | 1 1 1 |
| codewords | 1 0 1 0 1 | 1 1 1 0 0 |
| other cosets | 0 0 1 0 1 | 0 1 1 0 0 |
| | 1 1 1 0 1 | 1 0 1 0 0 |
| | 1 0 0 0 1 | 1 1 0 0 0 . |

Suppose $y_s$ = [1 1 1 1 | 1 1 1 1 1] is the received word then $S(y_s)$ = $H_s y_s^T \neq [(0) \mid (0)]$. $e_s$ = [0 1 0 0 | 0 0 1 0 0] is the super set coset leader. Thus $x_s = y_s + e_s$ = [1 0 1 1 | 1 1 0 1 1] $\in C_s$.

With the advent of computers calculating the super special coset leaders is not a very tedious job. Appropriate programs will yield the result in no time.

Now we proceed on to describe/define the super special row cyclic code.

**DEFINITION 3.1.8:** *Let $C_s = [C_1 \mid C_2 \mid \ldots \mid C_n]$ be a super special row code. If every $C_s^i$ is a cyclic code in $C_s$ we call $C_s$ to be a super special cyclic row code. $H_s = [H_1 \mid H_2 \mid \ldots \mid H_n]$ denotes the super special parity check row matrix of the super special cyclic code.*

We illustrate this by the following example.

***Example 3.1.9:*** Let $C_s = \left[C_s^1 \mid C_s^2 \mid C_s^3\right]$ be a super special cyclic code with an associated super special parity check matrix

$$H_s = [H_1 \mid H_2 \mid H_3]$$

$$= \left[\begin{array}{ccccccc} 1 & 1 & 1 & 0 & 1 & 0 & 0 \\ 0 & 1 & 1 & 1 & 0 & 1 & 0 \\ 0 & 0 & 1 & 1 & 1 & 0 & 1 \end{array}\right.$$

$$\left.\begin{array}{ccccccc|cccccc} 0 & 0 & 1 & 0 & 1 & 1 & 1 & 0 & 0 & 1 & 0 & 0 & 1 \\ 0 & 1 & 0 & 1 & 1 & 1 & 0 & 0 & 1 & 0 & 0 & 1 & 0 \\ 1 & 0 & 1 & 1 & 1 & 0 & 0 & 1 & 0 & 0 & 1 & 0 & 0 \end{array}\right].$$

We see each of $C_s^1$, $C_s^2$ and $C_s^3$ are cyclic codes.

Now we see in general for any super special mixed row code with an associated super special parity check matrix $H_s$ which happens to be a super mixed row matrix we cannot define the super special generator row matrix $G_s$. The simple reason being if the code words in each of the $C_s^i$ in $C_s$ where $C_s = \left[C_s^1 \mid C_s^2 \mid \ldots \mid C_s^n\right]$, i = 1, 2, …, n happens to be of different length then it would be impossible to define a super generator matrix. So we shall first define the notion of super special generator row matrix of a super special row code(mixed row code).

**DEFINITION 3.1.9:** *Let $C_s = \left[C_s^1 \mid C_s^2 \mid \ldots \mid C_s^n\right]$ be a super special row code. A super special row matrix which generates $C_s$ exists if and only if in each $C_s^i$ the codes in $C_s$ have the same number of message symbols, that is if $C_s$ has a super special parity check row matrix $H_s = [H_1 \mid H_2 \mid \ldots \mid H_n]$ then we demanded each $C_s^i$ must have the same number of check symbols. Likewise for the super special generator row matrix to exist we must have $G_s = [G_1 \mid G_2 \mid \ldots \mid G_n]$ where $C_s^i$ have the same number of message symbols which forms the number of rows of the super row generator matrix $G_s$.*

We shall first illustrate this by an example.

***Example 3.1.10:*** $C_s = \left[C_s^1 \mid C_s^2\right]$ be a super special row code.

Let $G_s = \left[ G_s^1 \mid G_s^2 \right]$

$$= \left[\begin{array}{ccccccc|ccccccc} 1 & 1 & 0 & 1 & 0 & 0 & 0 & 1 & 0 & 0 & 0 & 1 & 0 & 1 \\ 0 & 1 & 1 & 0 & 1 & 0 & 0 & 0 & 1 & 0 & 0 & 1 & 1 & 1 \\ 0 & 0 & 1 & 1 & 0 & 1 & 0 & 0 & 0 & 1 & 0 & 1 & 1 & 0 \\ 0 & 0 & 0 & 1 & 1 & 0 & 1 & 0 & 0 & 0 & 1 & 1 & 0 & 1 \end{array}\right]$$

be the super row generator matrix which generates $C_s$.

The code words of $C_s^1$ generated by $G_1$ is given by

$C_s^1$ = {(0 0 0 0 0 0 0), (1 1 0 1 0 0 0), (0 1 1 0 1 0 0),
(0 0 1 1 0 1 0), (0 0 0 1 1 0 1), (1 0 1 1 1 0 0),
(0 1 0 1 1 1 0), (0 0 1 0 1 1 1), (1 1 1 0 0 1 0),
(1 1 0 0 1 0 1), (0 1 1 1 0 0 1), (1 0 0 0 1 1 0),
(0 1 0 0 0 1 1), (1 0 1 0 0 0 1), (1 1 1 1 1 0 1), (1 0 0 1 0 1 1)}.

$C_s^2$ = {(0 0 0 0 0 0 0), (1 0 0 0 1 0 1), (0 1 0 0 1 1 1),
(0 0 1 0 1 1 0), (0 0 0 1 0 1 1), (1 1 0 0 0 1 0),
(0 1 1 0 0 0 1), (0 0 1 1 1 0 1), (1 0 0 1 1 1 0),
(1 0 1 0 0 1 1), (0 1 0 1 1 0 0), (1 1 1 0 1 0 0),
(0 1 1 1 0 1 0), (1 1 0 1 0 0 1), (1 0 1 1 0 0 0), (1 1 1 1 1 0 1)}.

If $x_s = \left[ x_s^1 \mid x_s^2 \right]$ by taking $x_s^1 \in C_s^1$ and $x_s^2 \in C_s^2$ we get $C_s$. Clearly elements in $C_s$ are super row vectors.

Now we proceed to define super special mixed row code of $C_s$ and its super special generator mixed row matrix.

**DEFINITION 3.1.10:** *Let $C_s = \left[ C_s^1 \mid C_s^2 \mid \ldots \mid C_s^n \right]$ be a super special mixed row code. If each of the codes $C_s^i$ have the same number of message symbols then we have the super special generator mixed row matrix $G_s = \left[ G_s^1 \mid G_s^2 \mid \ldots \mid G_s^n \right]$ associated with $C_s$. Number of message symbols in each of the $G_s^i$ are equal and is the super special mixed row matrix $G_s$; $1 \leq i \leq n$.*

We illustrate this by the following example.

***Example 3.1.11:*** Let $C_s = \left[C_s^1 \mid C_s^2 \mid C_s^3 \mid C_s^4\right]$ be a super special mixed row code. Let $G_s = \left[G_s^1 \mid G_s^2 \mid G_s^3 \mid G_s^4\right]$ be the associated super special mixed row generator matrix given by

$$G_s = \left[G_s^1 \mid G_s^2 \mid G_s^3 \mid G_s^4\right]$$

$$= \left[\begin{array}{ccccc|cccccc} 1 & 0 & 0 & 0 & 1 & 1 & 0 & 0 & 1 & 1 & 0 \\ 0 & 1 & 0 & 1 & 0 & 0 & 1 & 0 & 0 & 1 & 1 \\ 0 & 0 & 1 & 1 & 1 & 0 & 0 & 1 & 1 & 0 & 1 \end{array}\right.$$

$$\left.\begin{array}{|ccccccc|ccccccc} 1 & 0 & 0 & 0 & 0 & 1 & 0 & 1 & 0 & 0 & 1 & 0 & 0 & 1 \\ 0 & 1 & 0 & 1 & 0 & 0 & 0 & 0 & 1 & 0 & 0 & 1 & 1 & 0 \\ 0 & 0 & 1 & 0 & 1 & 0 & 1 & 0 & 0 & 1 & 1 & 0 & 1 & 0 \end{array}\right].$$

Clearly $G_s$ is a super mixed row matrix. All the codes generated by $G_1$, $G_2$, $G_3$ and $G_4$ have the same number of message symbols. The code $C_s^1$ = {(0 0 0 0 0), (1 0 0 0 1), (0 1 0 1 0), (0 0 1 1 1), (1 1 0 1 1), (0 1 1 0 1), (1 0 1 1 0), (1 1 1 0 0)}. The codewords given by $C_s^2$ = {(0 0 0 0 0 0), (1 0 0 1 1 0), (0 1 0 0 1 1), (0 0 1 1 0 1), (1 1 0 1 0 1), (0 1 1 1 1 0), (1 0 1 0 1 1), (1 1 1 0 0 0)}. The codes associated with $C_s^3$ = {1 0 0 0 0 1 0), (0 0 0 0 0 0 0), (0 1 0 1 0 0 0), (0 0 1 0 1 0 1), (1 1 0 1 0 1 0), (1 0 1 0 1 1 1), (0 1 1 1 1 0 1), (1 1 1 1 1 1 1)} and $C_s^4$ = {(0 0 0 0 0 0 0), (1 0 0 1 0 0 1), (0 1 0 0 1 1 0), (0 0 1 1 0 1 0), (1 1 0 1 1 1 1), (0 1 1 1 1 0 0), (1 0 1 0 0 1 1), (1 1 1 0 1 0 1)}.

We see the number of code words in each and every code $C_s^i$, i = 1, 2, 3, 4 is the same; equal to 8. Further the number of super code words in $C_s$ is 8 × 8 × 8 × 8 i.e., $|C_s| = 8^4$. We give a necessary and sufficient condition for $H_s$ to have an associated generator matrix $G_s$ for a super special row code $C_s$.

**THEOREM 3.1.1:** *Let* $C_s = \left[C_s^1 \mid C_s^2 \mid \ldots \mid C_s^n\right]$ *be a super special row code with* $H_s = [H_1 \mid H_2 \mid \ldots \mid H_n]$ *, the super special parity check matrix. If each* $H_i = (A_i, I_{n-k})$, *i = 1, 2, …, n then* $G_s = [G_1 \mid G_2 \mid \ldots \mid G_n]$ *with* $G_i = (I_k, -A^T)$; $1 \leq i \leq n$ *if and only if the length of each code word in* $C_s^i$ *is the same for i = 1, 2, …, n.*

*Proof:* Suppose we are given $H_s = [H_1 \mid H_2 \mid \ldots \mid H_n]$ to be the super special parity check matrix of the super special row code $C_s = [C_1 \mid C_2 \mid \ldots \mid C_n]$. We know every subcode $C_i$ of $C_s$ have the same number of check symbols. Suppose we have for this super special row code $C_s$ the super special generator row matrix $G_s$ with $G_s H_s^T = [(0) \mid (0) \mid \ldots \mid (0)]$ . Then we have $n_1 - k_1 = n_2 - k_2 = \ldots = n_n - k_n$ and $k_1 = k_2 = k_3 = \ldots = k_n$ this is possible if and only if $n_1 = n_2 = \ldots = n_n$. Hence the result.

Thus we see in this situation we have both the super special generator row matrix of the code $C_s$ as well as the super special parity check matrix of the code $C_s$ are not super special mixed row matrices; we see both of them have the same length n for each row.

We illustrate this situation by an example before we proceed on to define more concepts.

***Example 3.1.12:*** Let $C_s = \left[C_s^1 \mid C_s^2 \mid C_s^3 \mid C_s^4\right]$ be a super special row code. $H_s = [H_1 \mid H_2 \mid H_3 \mid H_4]$

$$= \left[\begin{array}{ccccccc|ccccccc} 0 & 1 & 1 & 1 & 0 & 0 & 0 & 1 & 0 & 0 & 1 & 0 & 0 & 0 \\ 1 & 0 & 1 & 0 & 1 & 0 & 0 & 0 & 1 & 1 & 0 & 1 & 0 & 0 \\ 1 & 1 & 1 & 0 & 0 & 1 & 0 & 1 & 1 & 1 & 0 & 0 & 1 & 0 \\ 1 & 1 & 0 & 0 & 0 & 0 & 1 & 1 & 0 & 1 & 0 & 0 & 0 & 1 \end{array}\right.$$

$$\left.\begin{array}{|ccccccc|ccccccc} 1 & 1 & 1 & 1 & 0 & 0 & 0 & 1 & 0 & 1 & 1 & 0 & 0 & 0 \\ 0 & 0 & 1 & 0 & 1 & 0 & 0 & 0 & 1 & 0 & 0 & 1 & 0 & 0 \\ 1 & 0 & 0 & 0 & 0 & 1 & 0 & 1 & 1 & 1 & 0 & 0 & 1 & 0 \\ 0 & 1 & 0 & 0 & 0 & 0 & 1 & 1 & 1 & 0 & 0 & 0 & 0 & 1 \end{array}\right]$$

the super special row parity check matrix; then the related super special row generator matrix

$$G_s = [G_1 \mid G_2 \mid G_3 \mid G_4]$$

$$= \left[\begin{array}{ccccccc|ccccccc} 1 & 0 & 0 & 0 & 1 & 1 & 1 & 1 & 0 & 0 & 1 & 0 & 1 & 1 \\ 0 & 1 & 0 & 1 & 0 & 1 & 1 & 0 & 1 & 0 & 0 & 1 & 1 & 0 \\ 0 & 0 & 1 & 1 & 1 & 1 & 0 & 0 & 0 & 1 & 0 & 1 & 1 & 1 \end{array}\right.$$

$$\left|\begin{array}{ccccccc|ccccccc} 1 & 0 & 0 & 1 & 0 & 1 & 0 & 1 & 0 & 0 & 1 & 0 & 1 & 1 \\ 0 & 1 & 0 & 1 & 0 & 0 & 1 & 0 & 1 & 0 & 0 & 1 & 1 & 1 \\ 0 & 0 & 1 & 1 & 1 & 0 & 0 & 0 & 0 & 1 & 1 & 0 & 1 & 0 \end{array}\right].$$

Now we find

$$G_s H_s^T = [G_1 \mid G_2 \mid G_3 \mid G_4] \times [H_1 \mid H_2 \mid H_3 \mid H_4]^T$$

$$= [G_1 \mid G_2 \mid G_3 \mid G_4] \times \left[ H_1^T \middle| H_2^T \middle| H_3^T \middle| H_4^T \right]$$

$$= \left[ G_1 H_1^T \middle| G_2 H_2^T \middle| G_3 H_3^T \middle| G_4 H_4^T \right]$$

$$= \left[ \begin{bmatrix} 1 & 0 & 0 & 0 & 1 & 1 & 1 \\ 0 & 1 & 0 & 1 & 0 & 1 & 1 \\ 0 & 0 & 1 & 1 & 1 & 1 & 0 \end{bmatrix} \begin{bmatrix} 0 & 1 & 1 & 1 \\ 1 & 0 & 1 & 1 \\ 1 & 1 & 1 & 0 \\ 1 & 0 & 0 & 0 \\ 0 & 1 & 0 & 0 \\ 0 & 0 & 1 & 0 \\ 0 & 0 & 0 & 1 \end{bmatrix} \right|$$

$$\begin{bmatrix} 1 & 0 & 0 & 1 & 0 & 1 & 1 \\ 0 & 1 & 0 & 0 & 1 & 1 & 0 \\ 0 & 0 & 1 & 0 & 1 & 1 & 1 \end{bmatrix}$$

$$
= \left.\begin{bmatrix} 1 & 0 & 1 & 1 \\ 0 & 1 & 1 & 0 \\ 0 & 1 & 1 & 1 \\ 1 & 0 & 0 & 0 \\ 0 & 1 & 0 & 0 \\ 0 & 0 & 1 & 0 \\ 0 & 0 & 0 & 1 \end{bmatrix}\right| \begin{bmatrix} 1 & 0 & 0 & 1 & 0 & 1 & 0 \\ 0 & 1 & 0 & 1 & 0 & 0 & 1 \\ 0 & 0 & 1 & 1 & 1 & 0 & 0 \end{bmatrix} \left| \begin{bmatrix} 1 & 0 & 1 & 0 \\ 1 & 0 & 0 & 1 \\ 1 & 1 & 0 & 0 \\ 1 & 0 & 0 & 0 \\ 0 & 1 & 0 & 0 \\ 0 & 0 & 1 & 0 \\ 0 & 0 & 0 & 1 \end{bmatrix}\right.
$$

$$
\left| \begin{bmatrix} 1 & 0 & 0 & 1 & 0 & 1 & 1 \\ 0 & 1 & 0 & 0 & 1 & 1 & 1 \\ 0 & 0 & 1 & 1 & 0 & 1 & 0 \end{bmatrix} \right| \left. \begin{bmatrix} 1 & 0 & 1 & 1 \\ 0 & 1 & 1 & 1 \\ 1 & 0 & 1 & 0 \\ 1 & 0 & 0 & 0 \\ 0 & 1 & 0 & 0 \\ 0 & 0 & 1 & 0 \\ 0 & 0 & 0 & 1 \end{bmatrix} \right]
$$

$$
= \left[\begin{array}{cccc|cccc|cccc|cccc} 0 & 0 & 0 & 0 & 0 & 0 & 0 & 0 & 0 & 0 & 0 & 0 & 0 & 0 & 0 & 0 \\ 0 & 0 & 0 & 0 & 0 & 0 & 0 & 0 & 0 & 0 & 0 & 0 & 0 & 0 & 0 & 0 \\ 0 & 0 & 0 & 0 & 0 & 0 & 0 & 0 & 0 & 0 & 0 & 0 & 0 & 0 & 0 & 0 \end{array}\right]
$$

which is a super special zero row vector.

Now we shall illustrate by an example, in which $G_s H_s^T \neq$ [(0) | (0) | … | (0)].

***Example 3.1.13:*** Let $C_s = [C_1 \mid C_2 \mid C_3]$ be a super special row code where $H_s = [H_1 \mid H_2 \mid H_3]$, is the associated super special parity check row matrix.

$$H_s = \left[\begin{array}{ccccccc|cccccccc}
1 & 0 & 1 & 1 & 0 & 0 & 0 & 1 & 1 & 0 & 0 & 1 & 0 & 0 & 0 \\
1 & 1 & 1 & 0 & 1 & 0 & 0 & 0 & 0 & 1 & 1 & 0 & 1 & 0 & 0 \\
0 & 1 & 0 & 0 & 0 & 1 & 0 & 1 & 0 & 0 & 1 & 0 & 0 & 1 & 0 \\
1 & 0 & 0 & 0 & 0 & 0 & 1 & 1 & 1 & 1 & 1 & 0 & 0 & 0 & 1
\end{array}\right.$$

$$\left|\begin{array}{cccccc}
1 & 0 & 1 & 0 & 0 & 0 \\
0 & 1 & 0 & 1 & 0 & 0 \\
1 & 1 & 0 & 0 & 1 & 0 \\
1 & 0 & 0 & 0 & 0 & 1
\end{array}\right].$$

Now

$$G_1 = \begin{bmatrix} 1 & 0 & 0 & 1 & 1 & 0 & 1 \\ 0 & 1 & 0 & 0 & 1 & 1 & 0 \\ 0 & 0 & 1 & 1 & 1 & 0 & 0 \end{bmatrix}$$

a $3 \times 7$ matrix, with $G_1H_1^T$ = (0).

$$G_2 = \begin{bmatrix} 1 & 0 & 0 & 0 & 1 & 0 & 1 & 1 \\ 0 & 1 & 0 & 0 & 1 & 0 & 0 & 1 \\ 0 & 0 & 1 & 0 & 0 & 1 & 0 & 1 \\ 0 & 0 & 0 & 1 & 0 & 1 & 1 & 1 \end{bmatrix}$$

a $4 \times 8$ generator matrix got from $H_2$ with $G_2H_2^T$ = (0) and

$$G_3 = \begin{bmatrix} 1 & 0 & 1 & 0 & 1 & 1 \\ 0 & 1 & 0 & 1 & 1 & 0 \end{bmatrix}$$

is a $2 \times 6$ matrix with $G_3H_3^T$ = (0).

We see $G_1$, $G_2$ and $G_3$ cannot be formed into a super mixed row matrix. So this example clearly shows to us that even if $C_s$ = [$C_1$ | $C_2$ | $C_3$] is a super special mixed row code with $H_s$ a super

special parity check mixed row matrix of $C_s$ yet $G_s$ is undefined for this $C_s$. This is in keeping with the theorem.

Likewise if we have a super special row code $C_s$ we may have the super special row matrix which generates $C_s$ yet $H_s$ may not exist.

This is the marked difference between the super special row codes and usual linear codes.

***Example 3.1.14:*** Let $C_s = [C_1 \mid C_2 \mid C_3]$ be a super special code with super special row generator matrix

$$G_s = \left[\begin{array}{cccccc|ccccccc|ccccc} 1 & 0 & 0 & 0 & 1 & 1 & 1 & 0 & 0 & 1 & 1 & 1 & 1 & 1 & 0 & 0 & 1 & 1 \\ 0 & 1 & 0 & 1 & 0 & 1 & 0 & 1 & 0 & 0 & 0 & 1 & 1 & 0 & 1 & 0 & 0 & 1 \\ 0 & 0 & 1 & 1 & 1 & 0 & 0 & 0 & 1 & 1 & 0 & 0 & 1 & 0 & 0 & 1 & 1 & 0 \end{array}\right]$$

$$= [G_1 \mid G_2 \mid G_3].$$

Clearly $G_s$ is a super special row mixed matrix. Now

$$H_1 = \begin{bmatrix} 0 & 1 & 1 & 1 & 0 & 0 \\ 1 & 0 & 1 & 0 & 1 & 0 \\ 1 & 1 & 0 & 0 & 0 & 1 \end{bmatrix}$$

got from $G_1$ and we have $G_1H_1^T = (0)$. The parity check matrix got from $G_2$ is

$$H_2 = \begin{bmatrix} 1 & 0 & 1 & 1 & 0 & 0 & 0 \\ 1 & 0 & 0 & 0 & 1 & 0 & 0 \\ 1 & 1 & 0 & 0 & 0 & 1 & 0 \\ 1 & 1 & 1 & 0 & 0 & 0 & 1 \end{bmatrix}$$

where $G_2H_2^T = (0)$. Also

$$H_3 = \begin{bmatrix} 1 & 0 & 1 & 1 & 0 \\ 1 & 1 & 0 & 0 & 1 \end{bmatrix}$$

and is such that $G_3H_3^T = (0)$. We see $H_1$, $H_2$ and $H_3$ cannot be made into a super special row matrix. Hence the claim.

Now having defined the new class of super special row (mixed row) codes we will now define new classes of mixed super classes of mixed super special row codes $C_s$ i.e., we may have the super special row code to contain classical subcodes as Hamming code or cyclic code or code and its orthogonal complement and so on.

***DEFINITION 3.1.11:*** *Let $C_s = [C_1 \mid C_2 \mid \ldots \mid C_n]$ be a super special row code. If some of the $C_i$'s are Hamming codes, some $C_j$'s are cyclic codes $i \neq j$, some $C_k$'s are repetition codes and some $C_t$'s are codes and $C_p$'s are dual codes of $C_t$'s ; $1 \leq j, k, t, i, p < n$ then we call $C_s$ to be a mixed super special row code.*

It is important to mention here that even if two types of classical codes are present still we call $C_s$ as a mixed super special row code.

We will illustrate them by the following examples.

***Example 3.1.15:*** Let $C_s = [C_1 \mid C_2 \mid C_3 \mid C_4]$ be a mixed super special row code. Here $C_1$ is a Hamming code, $C_2$ the repetition code, $C_3$ a code of no specific type and $C_4$ a cyclic code.

Let the mixed super special parity check matrix $H_s$ associated with $C_s$ be given by

$$H_s = [H_1 \mid H_2 \mid H_3 \mid H_4]$$

$$= \left[\begin{array}{ccccccc|cccc|} 0 & 0 & 0 & 1 & 1 & 1 & 1 & 1 & 1 & 0 & 0 \\ 0 & 1 & 1 & 0 & 0 & 1 & 1 & 1 & 0 & 1 & 0 \\ 1 & 0 & 1 & 0 & 1 & 0 & 1 & 1 & 0 & 0 & 1 \end{array}\right.$$

$$\left.\begin{array}{cccccc|ccccccc} 0 & 1 & 1 & 1 & 0 & 0 & 1 & 1 & 1 & 0 & 1 & 0 & 0 \\ 1 & 0 & 1 & 0 & 1 & 0 & 0 & 1 & 1 & 1 & 0 & 1 & 0 \\ 1 & 1 & 0 & 0 & 0 & 1 & 0 & 0 & 1 & 1 & 0 & 0 & 1 \end{array}\right].$$

Clearly $C_s$ is a mixed super special mixed row code. Any super code word $x_s$ of $C_s$ will be of a form $x_s$ = [1 0 0 0 1 0 0 | 1 1 1 1 | 0 1 1 0 1 1 | 1 1 1 1 1 1 0] which is clearly a super mixed row vector.

***Example 3.1.16:*** Let $C_s = [C_1 \mid C_2 \mid C_3]$ be a mixed super special row code. Let $H_s = [H_1 \mid H_2 \mid H_3]$ be the associated super special parity check mixed row matrix. $C_1$ is the Hamming code, $C_2$ any code and $C_3$ a repetition code.

$$H_s = \left[\begin{array}{cccccccc|ccccccc} 1 & 1 & 1 & 1 & 1 & 1 & 1 & 1 & 1 & 0 & 1 & 1 & 0 & 0 & 0 \\ 0 & 0 & 0 & 1 & 1 & 1 & 1 & 0 & 1 & 1 & 1 & 0 & 1 & 0 & 0 \\ 0 & 1 & 1 & 0 & 0 & 1 & 1 & 0 & 1 & 0 & 0 & 0 & 0 & 1 & 0 \\ 1 & 0 & 1 & 0 & 1 & 0 & 1 & 0 & 0 & 1 & 1 & 0 & 0 & 0 & 1 \end{array}\right.$$

$$\left|\begin{array}{ccccc} 1 & 1 & 0 & 0 & 0 \\ 1 & 0 & 1 & 0 & 0 \\ 1 & 0 & 0 & 1 & 0 \\ 1 & 0 & 0 & 0 & 1 \end{array}\right]$$

is the mixed super special parity check mixed row matrix for which $G_s$ does not exist.

We define the new notion of super special Hamming row code.

**DEFINITION 3.1.12:** *Let* $C_S = \left[C_S^1 \mid C_S^2 \mid \ldots \mid C_S^n\right]$ *where each* $C_s^i$ *is a* $(2^m - 1, 2^m - 1 - m)$ *Hamming code for i = 1, 2, …, n. Then we call* $C_s$ *to be a super special Hamming row code. If* $H_s = [H_1 \mid H_2 \mid \ldots \mid H_n]$ *be the super special parity check matrix associated with* $C_s$ *we see* $H_s$ *is a super special row matrix having m rows and each parity check matrix* $H_i$ *has m rows and* $2^m - 1$ *columns, i =1, 2, …, n.*

Further the transmission rate can never be equal to ½. If m > 2 then will the transmission rate be always greater than ½ ?

We will just illustrate a super special Hamming row code by the following example.

***Example 3.1.17:*** Let $C_s = \left[C_s^1 \mid C_s^2 \mid C_s^3\right]$ be a super special Hamming row code where

$$H_s = \left[\begin{array}{ccccccc|ccccccc|ccccccc} 0 & 0 & 0 & 1 & 1 & 1 & 1 & 1 & 1 & 1 & 0 & 1 & 0 & 0 & 1 & 0 & 0 & 1 & 1 & 0 & 1 \\ 0 & 1 & 1 & 0 & 0 & 1 & 1 & 0 & 1 & 1 & 1 & 0 & 1 & 0 & 0 & 1 & 0 & 1 & 0 & 1 & 1 \\ 1 & 0 & 1 & 0 & 1 & 0 & 1 & 0 & 0 & 1 & 1 & 1 & 0 & 1 & 0 & 0 & 1 & 0 & 1 & 1 & 1 \end{array}\right]$$

$$= [H_1 \mid H_2 \mid H_3]$$

is the super special matrix associated with $C_s$. We see $C_s^i$ is a (7, 4) Hamming code; i = 1, 2, 3. Here m = 3 and $n = 2^3 - 1 = 7$.

The code words associated with $C_s^1$ is

{(0 0 0 0 0 0 0), (1 0 0 0 0 0 1 1), (0 1 0 0 1 0 1), (0 0 1 0 1 1 0),
(0 0 0 1 1 1 1), (1 1 0 0 1 1 0), (0 1 1 0 0 1 1), (0 0 1 1 0 0 1),
(1 0 1 0 1 0 1), (0 1 0 1 1 1 1), (1 0 0 1 1 0 0), (1 1 1 0 0 0 0),
(0 1 1 1 1 0 0), (1 1 0 1 0 0 1), (1 0 1 1 1 0 1), (1 1 1 1 1 1 1)}.

$C_s^2$ = {(0 0 0 0 0 0 0), (10 0 0 1 0 1), (0 1 0 0 1 1 1),
(0 0 1 0 1 1 1), (0 0 0 1 0 1 1), (1 1 0 0 0 1 0),
(0 1 1 0 0 0 1), (0 0 1 1 1 0 1), (1 0 0 1 1 1 1),
(1 0 1 0 0 1 1), (0 1 0 1 1 0 1), (1 1 1 0 1 0 1),
(1 1 0 1 0 0 1), (0 1 1 1 0 1 1), (1 0 1 1 0 0 0), (1 1 1 1 1 1 1)}.

$C_s^3$ = {(0 0 0 0 0 0 0), (1 0 0 0 0 1 1), (0 1 0 0 1 0 1),

(0 0 1 0 1 1 1), (0 0 0 1 1 1 0), (1 1 0 0 1 1 0),
(0 1 1 0 0 1 0), (0 0 1 1 0 0 1), (1 0 0 1 1 1 0),
(1 0 1 0 1 0 0), (0 1 0 1 0 1 1), (1 1 1 0 0 0 1),
(0 1 1 1 1 0 0), (1 1 0 1 0 1 1), (1 0 1 1 0 1 0), (1 1 1 1 1 1 1)}.

By taking one code word from $C_s^1$, one code word from $C_s^2$ and one from $C_s^3$ we form the super special Hamming row code. Any $x_s = (x_s^1 \mid x_s^2 \mid x_s^3)$ where $x_s^1 \in C_s^1$, $x_s^2 \in C_s^2$ and $x_s^3 = C_s^3$.

We see

$$H_s\ x_s^T = [H_1 \mid H_2 \mid H_3]\ (x_s^1 \mid x_s^2 \mid x_s^3)$$
$$= [H_1 (x_s^1)^T \mid H_2 (x_s^2)^T \mid H_3 (x_s^3)^T]$$
$$= [(0) \mid (0) \mid (0)]$$

for every $x_s \in C_s$.

## 3.2 New Classes of Super Special Column Codes

Suppose we are interested in finding super special codes in which the number of check symbols will be different for the subcodes. In such a situation we see certainly we cannot work with the super special row codes $C_s$, for in this case we demand always the number of check symbols to be the same for every subcode in $C_s$, so we are forced to define this special or new classes of super special codes.

**DEFINITION 3.2.1:** *Suppose we have to describe n codes each of same length say m but with varying sets of check symbols by a single matrix. Then we define it using super column matrix as the super code parity check matrix. Let*

$$C_s = \begin{bmatrix} C_1 \\ \hline C_2 \\ \hline \vdots \\ \hline C_m \end{bmatrix}$$

*be a set of m codes, $C_1$, $C_2$, ..., $C_m$ where all of them have the same length n but have $n - k_1$, $n - k_2$, ..., $n - k_m$ to be the number of check symbols and $k_1$, $k_2$, ..., $k_m$ are the number of message symbols associated with each of the codes $C_1$, $C_2$, ..., $C_m$ respectively.*

Let us consider

$$H^s = \begin{bmatrix} H^1 \\ \hline H^2 \\ \hline \vdots \\ \hline H^m \end{bmatrix}$$

where each $H^i$ is the $n - k_i \times n$ parity check matrix of the code $C_i$; i = 1, 2, ..., m. We call $H^s$ to be the super special parity check mixed column matrix of $C_s$ and $C_s$ is defined as the super special mixed column code.

The main difference between the super special row code and the super special column code is that in super special row codes always the number of check symbols in every code in $C_s$ is the same as the number of message symbols in $C_i$ and the length of the code $C_i$ can vary where as in the super special column code, we will always have the same length for every code $C_i$ in $C_s$ but the number of message symbols and the number check symbols for each and every code $C_i$ in $C_s$ need not be the same. In case if the number of check symbols in each and every code $C_i$ is the same. Then we call $C_s$ to be a super special column code.

In case when we have varying number of check symbols then we call the code $C_s$ to be a super special mixed column code.

In the case of super special column code $C_s = [C_1 \mid C_2 \mid \ldots \mid C_m]^t$ we see every code $C_i$ in $C_s$ have the same number of message symbols. Thus every code is a (n, k) code. It may so happen that some of the $C_i$ and $C_j$ are identical codes. Now we proceed on to give examples of the two types of super special column codes.

***Example 3.2.1:*** Let

$$C_s = \left[\begin{array}{c} C_1 \\ \hline C_2 \\ \hline C_3 \\ \hline C_4 \end{array}\right]$$

be a super special column code.
Let

$$H_s = \left[\begin{array}{c} H_1 \\ \hline H_2 \\ \hline H_3 \\ \hline H_4 \end{array}\right]$$

where each $H_i$ is a 3 × 7 parity check matrix.

$$H_s = \left[\begin{array}{ccccccc} 1 & 0 & 1 & 0 & 1 & 0 & 0 \\ 0 & 1 & 0 & 1 & 0 & 1 & 0 \\ 1 & 1 & 1 & 0 & 0 & 0 & 1 \\ \hline 0 & 1 & 1 & 0 & 1 & 0 & 0 \\ 1 & 0 & 0 & 1 & 0 & 1 & 0 \\ 1 & 0 & 1 & 0 & 0 & 0 & 1 \\ \hline 1 & 0 & 0 & 1 & 1 & 0 & 0 \\ 1 & 1 & 1 & 0 & 0 & 1 & 0 \\ 1 & 0 & 0 & 1 & 0 & 0 & 1 \\ \hline 0 & 1 & 1 & 0 & 1 & 0 & 0 \\ 1 & 1 & 1 & 1 & 0 & 1 & 0 \\ 1 & 0 & 1 & 1 & 0 & 0 & 1 \end{array}\right]$$

is a super column matrix. We see if

$$H_s\,[x_1 \mid x_2 \mid x_3 \mid x_4]^T = (0) = \left[\begin{array}{c} H_1 \\ \hline H_2 \\ \hline H_3 \\ \hline H_4 \end{array}\right] \left[x_1^T \mid x_2^T \mid x_3^T \mid x_4^T\right]$$

$$= \begin{bmatrix} H_1 x_1^T \\ \hline H_2 x_2^T \\ \hline H_3 x_3^T \\ \hline H_4 x_4^T \end{bmatrix} = \begin{bmatrix} 0 \\ 0 \\ 0 \\ \hline 0 \\ 0 \\ 0 \\ \hline 0 \\ 0 \\ 0 \\ \hline 0 \\ 0 \\ 0 \end{bmatrix}.$$

As in case of super special row codes we get any super special column code, the main criteria being every code $C_i$ in $C_s = [C_1 \mid C_2 \mid \ldots \mid C_m]$ would have the same length. It can have any number of message symbols and any arbitrary number of check symbols.

Just we saw in example 3.2.1 a super special column code. One of the main reason for us to have this new class of codes is that when we have super special row codes we see if the length of the code words are the same then automatically it is such that the number of message symbols become equal to the number of check symbols Thus the transmission rate becomes fixed equal to 1/2. But if we wish to have lesser transmission rate we cannot get it from these codes, the super special column codes comes handy.

We now describe yet another super special column code $C_s$.

***Example 3.2.2:*** Let $C_s = [C_1 \mid C_2 \mid \ldots \mid C_n]^t$ be a super special column code. Here n = 3 i.e., $C_s = [C_1 \mid C_2 \mid C_3]^t$. Let

$$H_s = \left[\begin{array}{ccccccc} 0 & 0 & 1 & 1 & 0 & 0 & 0 \\ 0 & 1 & 0 & 0 & 1 & 0 & 0 \\ 1 & 1 & 1 & 0 & 0 & 1 & 0 \\ 1 & 0 & 0 & 0 & 0 & 0 & 1 \\ \hline 1 & 1 & 0 & 0 & 0 & 1 & 0 \\ 1 & 1 & 0 & 1 & 0 & 0 & 1 \\ \hline 1 & 0 & 1 & 0 & 1 & 0 & 0 \\ 0 & 1 & 1 & 0 & 0 & 1 & 0 \\ 1 & 1 & 1 & 1 & 0 & 0 & 1 \end{array}\right].$$

We see every code word in $C_i$ are of length 7. But the number of message symbols in $C_1$ is 3 and the number of check symbols is 4 i.e., $C_1$ is a (7, 3) code where as $C_2$ is a (7, 5) code and $C_3$ is a (7, 4) code. The number of super code words in the code $C_s$ is $2^3 \times 2^5 \times 2^4 = 2^{12}$.

The code words in

$C_1$ = {(0 0 0 0 0 0 0), (1 0 0 0 0 1 1), (0 1 0 0 1 1 0), (0 0 1 1 0 1 0), (1 1 0 0 1 0 1), (0 1 1 1 1 0 0), (1 0 1 1 0 0 1), (1 1 1 1 1 1 1)}.

The code

$C_2$ = {(0 0 0 0 0 0 0), (1 0 0 0 0 1 1), (0 1 0 0 0 1 1), (0 0 1 0 0 0 0), (0 0 0 1 0 0 1), (0 0 0 0 1 0 0), (1 1 0 0 0 0 0), (0 1 1 0 0 1 1), (0 0 1 1 0 0 1), (0 0 0 1 1 0 1), (1 0 1 0 0 1 1), (1 0 0 1 0 1 0), (1 0 0 0 1 1 1), (0 1 0 1 0 0 1), (0 0 1 0 1 0 0), (0 1 0 0 1 1 1), (1 1 1 0 0 0 0), (0 1 1 1 0 1 1), (0 0 1 1 1 0 1), (1 1 0 1 0 0 1), (1 1 0 0 1 0 0), (1 0 1 1 0 1 0), (1 0 0 1 1 1 0), (0 1 1 0 1 1 0), (1 0 1 0 1 1 1), (0 1 0 1 1 1 0), (1 1 1 1 0 0 1), (0 1 1 1 1 1 0), (1 1 1 0 1 0 0), (1 1 0 1 1 0 1), (1 0 1 1 1 1 0), (1 1 1 1 1 0 1)}

and

$C_3$ = {(0 0 0 0 0 0 0), (1 0 0 0 1 0 1), (0 1 0 0 0 1 1), (0 0 1 0 1 1 1), (0 0 0 1 0 0 1), (1 1 0 0 1 1 0), (0 1 1 0 1 0 0),

(0 0 1 1 1 1 0), (1 0 1 0 0 1 0), (0 1 0 1 0 1 0), (1 0 0 1 1 0 0),
(1 1 1 0 0 0 1), (0 1 1 1 1 0 1), (1 1 0 1 1 1 1),
(1 0 1 1 0 1 1), (1 1 1 1 0 0 0)}.

Thus by taking one code word from each of the $C_i$'s we get the super special code word which is super row vector of length 7. Thus $C_s$ has $8 \times 32 \times 16 = 4096$ super code words.

Unlike in the case of super special row codes given a super special column code $C_s$ with associated super special parity check column matrix $H_s$ we will always be in a position to a get the super special generator matrix $G_s$. We see in the example 3.2.2 just given we do not have for that $H_s$ an associated generator matrix $G_s$ though clearly each $H_i$ in $H_s$ is only in the standard form.

However for the example 3.2.1 we have an associated generator matrix $G_s$ for the

$$H_s = \left[\begin{array}{c} H_1 \\ \hline H_2 \\ \hline H_3 \\ \hline H_4 \end{array}\right]$$

$$= \left[\begin{array}{ccccccc} 1 & 0 & 1 & 0 & 1 & 0 & 0 \\ 0 & 1 & 0 & 1 & 0 & 1 & 0 \\ 1 & 1 & 1 & 0 & 0 & 0 & 1 \\ \hline 0 & 1 & 1 & 0 & 1 & 0 & 0 \\ 1 & 0 & 0 & 1 & 0 & 1 & 0 \\ 1 & 0 & 1 & 0 & 0 & 0 & 1 \\ \hline 1 & 0 & 0 & 1 & 1 & 0 & 0 \\ 1 & 1 & 1 & 0 & 0 & 1 & 0 \\ 1 & 0 & 0 & 1 & 0 & 0 & 1 \\ \hline 0 & 1 & 1 & 0 & 1 & 0 & 0 \\ 1 & 1 & 1 & 1 & 0 & 1 & 0 \\ 1 & 0 & 1 & 1 & 0 & 0 & 1 \end{array}\right].$$

$$G_5 = \left[\begin{array}{ccccccc}
1 & 0 & 0 & 0 & 1 & 0 & 1 \\
0 & 1 & 0 & 0 & 0 & 1 & 1 \\
0 & 0 & 1 & 0 & 1 & 0 & 1 \\
0 & 0 & 0 & 1 & 0 & 1 & 0 \\
\hline
1 & 0 & 0 & 0 & 0 & 1 & 1 \\
0 & 1 & 0 & 0 & 1 & 0 & 0 \\
0 & 0 & 1 & 0 & 1 & 0 & 1 \\
0 & 0 & 0 & 1 & 0 & 1 & 0 \\
\hline
1 & 0 & 0 & 0 & 1 & 1 & 1 \\
0 & 1 & 0 & 0 & 0 & 1 & 0 \\
0 & 0 & 1 & 0 & 0 & 1 & 0 \\
0 & 0 & 0 & 1 & 1 & 0 & 1 \\
\hline
1 & 0 & 0 & 0 & 0 & 1 & 1 \\
0 & 1 & 0 & 0 & 1 & 1 & 0 \\
0 & 0 & 1 & 0 & 1 & 1 & 1 \\
0 & 0 & 0 & 1 & 0 & 1 & 1
\end{array}\right] = \left[\begin{array}{c} G_1 \\ \hline G_2 \\ \hline G_3 \\ \hline G_4 \end{array}\right].$$

Clearly

$$G_s H_s^T = \left[\begin{array}{c} G_1 \\ \hline G_2 \\ \hline G_3 \\ \hline G_4 \end{array}\right] \left[ H_1^T \mid H_2^T \mid H_3^T \mid H_4^T \right].$$

$$= \left[\begin{array}{c} G_1 H_1^T \\ \hline G_2 H_2^T \\ \hline G_3 H_3^T \\ \hline G_4 H_4^T \end{array}\right]$$

$$
= \left[\begin{array}{ccccccc}
1 & 0 & 0 & 0 & 1 & 0 & 1 \\
0 & 1 & 0 & 0 & 0 & 1 & 1 \\
0 & 0 & 1 & 0 & 1 & 0 & 1 \\
0 & 0 & 0 & 1 & 0 & 1 & 0 \\
\hline
1 & 0 & 0 & 0 & 0 & 1 & 1 \\
0 & 1 & 0 & 0 & 1 & 0 & 0 \\
0 & 0 & 1 & 0 & 1 & 0 & 1 \\
0 & 0 & 0 & 1 & 0 & 1 & 0 \\
\hline
1 & 0 & 0 & 0 & 1 & 1 & 1 \\
0 & 1 & 0 & 0 & 0 & 1 & 0 \\
0 & 0 & 1 & 0 & 0 & 1 & 0 \\
0 & 0 & 0 & 1 & 1 & 0 & 1 \\
\hline
1 & 0 & 0 & 0 & 0 & 1 & 1 \\
0 & 1 & 0 & 0 & 1 & 1 & 0 \\
0 & 0 & 1 & 0 & 1 & 1 & 1 \\
0 & 0 & 0 & 1 & 0 & 1 & 1
\end{array}\right]
$$

$$
\times \left[\begin{array}{ccc|ccc|ccc|ccc}
1 & 0 & 1 & 0 & 1 & 1 & 1 & 1 & 1 & 0 & 1 & 1 \\
0 & 1 & 1 & 1 & 0 & 0 & 0 & 1 & 0 & 1 & 1 & 0 \\
1 & 0 & 1 & 1 & 0 & 1 & 0 & 1 & 0 & 1 & 1 & 1 \\
0 & 1 & 0 & 0 & 1 & 0 & 1 & 0 & 1 & 0 & 1 & 1 \\
1 & 0 & 0 & 1 & 0 & 0 & 1 & 0 & 0 & 1 & 0 & 0 \\
0 & 1 & 0 & 0 & 1 & 0 & 0 & 1 & 0 & 0 & 1 & 0 \\
0 & 0 & 1 & 0 & 0 & 1 & 0 & 0 & 1 & 0 & 0 & 1
\end{array}\right]
$$

$$= \left[\begin{array}{ccc} 0 & 0 & 0 \\ 0 & 0 & 0 \\ 0 & 0 & 0 \\ 0 & 0 & 0 \\ \hline 0 & 0 & 0 \\ 0 & 0 & 0 \\ 0 & 0 & 0 \\ 0 & 0 & 0 \\ \hline 0 & 0 & 0 \\ 0 & 0 & 0 \\ 0 & 0 & 0 \\ 0 & 0 & 0 \\ \hline 0 & 0 & 0 \\ 0 & 0 & 0 \\ 0 & 0 & 0 \\ 0 & 0 & 0 \end{array}\right].$$

We find conditions under which we have a super special column code $C_s$ with an associated super special parity check column matrix $H_s$ in the standard form to have a super special generator column matrix $G_s$ with $G_s H_s^T = \left[\begin{array}{c} 0 \\ \hline \vdots \\ \hline 0 \end{array}\right]$.

Before we prove a result of this nature we define the two types of super special generator column matrix for a super special column code $C_s$.

**DEFINITION 3.2.2**: *Let $C_s = [C_1 | C_2 | \ldots | C_n]^t$ where $C_i$'s are codes of same length m. Suppose each $C_i$ is generated by a matrix $G_i$, i = 1, 2, …, n, then*

$$G_s = \left[\begin{array}{c} G_1 \\ \hline \vdots \\ \hline G_n \end{array}\right]$$

*generates the super special column code $C_s$. We call $G_s$ the super special generator column matrix which generates $C_s$. If in each of the codes $C_i$ in $C_s$, we have same number of message symbols then we call $G_s$ to be a super special generator column matrix; i = 1, 2, …, n. If each of the codes $C_i$'s in $C_s$ have different number of message symbols then we call $G_s$ to be a super special generator mixed column matrix.*

We say $G_s$ is in the standard form only if each $G_i$ is in the standard form. Further only when $G_s$ is in the standard form and $G_s$ is a super special column matrix which is not a mixed matrix we have $H_s$ the super special parity check column matrix of the same $C_s$ with

$$G_s H_s^T = \begin{bmatrix} 0 \\ \hline 0 \\ \hline \vdots \\ \hline 0 \end{bmatrix}.$$

Now we illustrate both the situations by examples.

***Example 3.2.3:*** Let $C_s = [C_1 \mid C_2 \mid C_3 \mid C_4]^t$ be a super special column code generated by

$$G_s = \begin{bmatrix} 1 & 0 & 0 & 1 & 1 & 0 & 1 \\ 0 & 1 & 0 & 0 & 1 & 1 & 0 \\ 0 & 0 & 1 & 1 & 0 & 0 & 1 \\ \hline 1 & 0 & 0 & 0 & 1 & 0 & 1 \\ 0 & 1 & 0 & 1 & 0 & 0 & 0 \\ 0 & 0 & 1 & 0 & 0 & 1 & 1 \\ \hline 1 & 0 & 0 & 1 & 1 & 1 & 1 \\ 0 & 1 & 0 & 0 & 1 & 0 & 1 \\ 0 & 0 & 1 & 1 & 0 & 1 & 0 \\ \hline 1 & 0 & 0 & 1 & 0 & 1 & 1 \\ 0 & 1 & 0 & 0 & 1 & 1 & 0 \\ 0 & 0 & 1 & 0 & 0 & 1 & 1 \end{bmatrix} = \begin{bmatrix} G_1 \\ \hline G_2 \\ \hline G_3 \\ \hline G_4 \end{bmatrix}$$

the super special column matrix. Now the related $H_s$ of $G_s$ exists and is given by

$$H = \left[\begin{array}{c} H_1 \\ \hline H_2 \\ \hline H_3 \\ \hline H_4 \end{array}\right] = \left[\begin{array}{ccccccc} 1 & 0 & 1 & 1 & 0 & 0 & 0 \\ 1 & 1 & 0 & 0 & 1 & 0 & 0 \\ 0 & 1 & 0 & 0 & 0 & 1 & 0 \\ 1 & 0 & 1 & 0 & 0 & 0 & 1 \\ \hline 0 & 1 & 0 & 1 & 0 & 0 & 0 \\ 1 & 0 & 0 & 0 & 1 & 0 & 0 \\ 0 & 0 & 1 & 0 & 0 & 1 & 0 \\ 1 & 0 & 1 & 0 & 0 & 0 & 1 \\ \hline 1 & 0 & 1 & 1 & 0 & 0 & 0 \\ 1 & 1 & 0 & 0 & 1 & 0 & 0 \\ 1 & 0 & 1 & 0 & 0 & 1 & 0 \\ 1 & 1 & 0 & 0 & 0 & 0 & 1 \\ \hline 1 & 0 & 0 & 1 & 0 & 0 & 0 \\ 0 & 1 & 0 & 0 & 1 & 0 & 0 \\ 1 & 1 & 1 & 0 & 0 & 1 & 0 \\ 1 & 0 & 1 & 0 & 0 & 0 & 1 \end{array}\right].$$

Now

$$G_s\ H_s^T = \left[\begin{array}{c} G_1 \\ \hline G_2 \\ \hline G_3 \\ \hline G_4 \end{array}\right] \left[\begin{array}{c} H_1 \\ \hline H_2 \\ \hline H_3 \\ \hline H_4 \end{array}\right]^T = \left[\begin{array}{c} G_1H_1^T \\ \hline G_2H_2^T \\ \hline G_3H_3^T \\ \hline G_4H_4^T \end{array}\right] = \left[\begin{array}{cccc} 0 & 0 & 0 & 0 \\ 0 & 0 & 0 & 0 \\ 0 & 0 & 0 & 0 \\ \hline 0 & 0 & 0 & 0 \\ 0 & 0 & 0 & 0 \\ 0 & 0 & 0 & 0 \\ \hline 0 & 0 & 0 & 0 \\ 0 & 0 & 0 & 0 \\ 0 & 0 & 0 & 0 \\ \hline 0 & 0 & 0 & 0 \\ 0 & 0 & 0 & 0 \\ 0 & 0 & 0 & 0 \end{array}\right].$$

***Example 3.2.4:*** Let $C_s = [C_1 \mid C_2 \mid C_3]^t$ be a super special column code. Suppose

$$G_s = \left[\begin{array}{ccccccc} 1 & 0 & 0 & 1 & 0 & 0 & 0 \\ 0 & 1 & 0 & 1 & 1 & 0 & 0 \\ 0 & 0 & 1 & 0 & 1 & 1 & 0 \\ \hline 1 & 0 & 1 & 1 & 0 & 0 & 1 \\ 0 & 1 & 0 & 0 & 1 & 1 & 0 \\ \hline 1 & 0 & 0 & 0 & 0 & 1 & 0 \\ 0 & 1 & 0 & 0 & 1 & 1 & 0 \\ 0 & 0 & 1 & 0 & 1 & 1 & 1 \\ 0 & 0 & 0 & 1 & 0 & 1 & 1 \end{array}\right] = \left[\begin{array}{c} G_1 \\ \hline G_2 \\ \hline G_3 \end{array}\right]$$

be the super special column matrix which generates $C_s$.

$$H_1 = \begin{bmatrix} 1 & 1 & 0 & 1 & 0 & 0 & 0 \\ 0 & 1 & 1 & 0 & 1 & 0 & 0 \\ 0 & 0 & 1 & 0 & 0 & 1 & 0 \\ 0 & 0 & 0 & 0 & 0 & 0 & 1 \end{bmatrix},$$

$$H_2 = \begin{bmatrix} 1 & 0 & 1 & 0 & 0 & 0 & 0 \\ 1 & 0 & 0 & 1 & 0 & 0 & 0 \\ 0 & 1 & 0 & 0 & 1 & 0 & 0 \\ 0 & 1 & 0 & 0 & 0 & 1 & 0 \\ 1 & 0 & 0 & 0 & 0 & 0 & 1 \end{bmatrix}$$

and

$$H_3 = \begin{bmatrix} 0 & 1 & 1 & 0 & 1 & 0 & 0 \\ 1 & 1 & 1 & 1 & 0 & 1 & 0 \\ 0 & 0 & 1 & 1 & 0 & 0 & 1 \end{bmatrix}.$$

Thus

$$H_s = \left[\begin{array}{ccccccc} 1 & 1 & 0 & 1 & 0 & 0 & 0 \\ 0 & 1 & 1 & 0 & 1 & 0 & 0 \\ 0 & 0 & 1 & 0 & 0 & 1 & 0 \\ 0 & 0 & 0 & 0 & 0 & 0 & 1 \\ \hline 1 & 0 & 1 & 0 & 0 & 0 & 0 \\ 1 & 0 & 0 & 1 & 0 & 0 & 0 \\ 0 & 1 & 0 & 0 & 1 & 0 & 0 \\ 0 & 1 & 0 & 0 & 0 & 1 & 0 \\ 1 & 0 & 0 & 0 & 0 & 0 & 1 \\ \hline 0 & 1 & 1 & 0 & 1 & 0 & 0 \\ 1 & 1 & 1 & 1 & 0 & 1 & 0 \\ 0 & 0 & 1 & 1 & 0 & 0 & 1 \end{array}\right]$$

is the super special column matrix which is the super parity check column matrix of $C_s$. Clearly $G_s H_s^T = (0)$, super column zero matrix. As in case of usual codes given the parity check matrix H or the generator matrix G in the standard form we can always get G from H or H from G and we have $GH^T = (0)$. Like wise in the case of super special (mixed column) column code $C_s$ if $G_s$ is the super special column generator matrix in the standard form we can always get the super special column parity check matrix $H_s$ from $H_s$ and we have $G_s H_s^T$ is always a super special zero column matrix.

As in case of super special row code we can define classical super special column codes.

**DEFINITION 3.2.3:** *Let $C_s = [C_1 \mid C_2 \mid \ldots \mid C_n]^t$ be a super special column code if each of the code $C_i$ is a repetition code of length n then $C_s$ is a super special repetition column code with $C_1 = C_2 = \ldots = C_n$. The super special column parity check matrix*

$$H_s = \left[\begin{array}{c} H \\ \hline \vdots \\ \hline H \end{array}\right]$$

$$= \left[\begin{array}{cccccc} 1 & 1 & 0 & 0 & \ldots & 0 \\ 1 & 0 & 1 & 0 & \ldots & 0 \\ \vdots & \vdots & \vdots & \vdots & & \vdots \\ 1 & 0 & 0 & 0 & \ldots & 1 \\ \hline & & & \vdots & & \\ \hline 1 & 1 & 0 & 0 & \ldots & 0 \\ 1 & 0 & 1 & 0 & \ldots & 0 \\ \vdots & \vdots & \vdots & \vdots & & \vdots \\ 1 & 0 & 0 & 0 & \ldots & 1 \end{array}\right]$$

where

$$H = \left[\begin{array}{cccccc} 1 & 1 & 0 & 0 & \ldots & 0 \\ 1 & 0 & 1 & 0 & \ldots & 0 \\ \vdots & \vdots & \vdots & \vdots & & \vdots \\ 1 & 0 & 0 & 0 & \ldots & 1 \end{array}\right].$$

It is important to note that unlike the super special repetition row code which can have different lengths the super special repetition column code can have only a fixed length and any super special code word

$$x_s = \left[ x_s^1 \mid x_s^2 \mid \cdots \mid x_s^n \right]$$

where $x_s^j$ = (1 1 … 1), n-times or (0 0 … 0) n-times only; $1 \leq j \leq n$.

***Example 3.2.5:*** Let

$$C_s = \left[\begin{array}{c} C_1 \\ \hline C_2 \\ \hline C_3 \end{array}\right]$$

be a super special repetition column code where the related super special column parity check matrix

$$H_s = \begin{bmatrix} 1 & 1 & 0 & 0 & 0 & 0 \\ 1 & 0 & 1 & 0 & 0 & 0 \\ 1 & 0 & 0 & 1 & 0 & 0 \\ 1 & 0 & 0 & 0 & 1 & 0 \\ 1 & 0 & 0 & 0 & 0 & 1 \\ \hline 1 & 1 & 0 & 0 & 0 & 0 \\ 1 & 0 & 1 & 0 & 0 & 0 \\ 1 & 0 & 0 & 1 & 0 & 0 \\ 1 & 0 & 0 & 0 & 1 & 0 \\ 1 & 0 & 0 & 0 & 0 & 1 \\ \hline 1 & 1 & 0 & 0 & 0 & 0 \\ 1 & 0 & 1 & 0 & 0 & 0 \\ 1 & 0 & 0 & 1 & 0 & 0 \\ 1 & 0 & 0 & 0 & 1 & 0 \\ 1 & 0 & 0 & 0 & 0 & 1 \end{bmatrix}.$$

Thus

$$\begin{aligned} C_s = \{&[0\ 0\ 0\ 0\ 0\ 0 \mid 0\ 0\ 0\ 0\ 0\ 0 \mid 0\ 0\ 0\ 0\ 0\ 0], \\ &[1\ 1\ 1\ 1\ 1\ 1 \mid 1\ 1\ 1\ 1\ 1\ 1 \mid 1\ 1\ 1\ 1\ 1\ 1], \\ &[0\ 0\ 0\ 0\ 0\ 0 \mid 0\ 0\ 0\ 0\ 0\ 0 \mid 1\ 1\ 1\ 1\ 1\ 1], \\ &[0\ 0\ 0\ 0\ 0\ 0 \mid 1\ 1\ 1\ 1\ 1\ 1 \mid 1\ 1\ 1\ 1\ 1\ 1], \\ &[1\ 1\ 1\ 1\ 1\ 1 \mid 0\ 0\ 0\ 0\ 0\ 0 \mid 0\ 0\ 0\ 0\ 0\ 0], \\ &[1\ 1\ 1\ 1\ 1\ 1 \mid 1\ 1\ 1\ 1\ 1\ 1 \mid 0\ 0\ 0\ 0\ 0\ 0], \\ &[1\ 1\ 1\ 1\ 1\ 1 \mid 0\ 0\ 0\ 0\ 0\ 0 \mid 1\ 1\ 1\ 1\ 1\ 1], \\ &[0\ 0\ 0\ 0\ 0\ 0 \mid 1\ 1\ 1\ 1\ 1\ 1 \mid 0\ 0\ 0\ 0\ 0\ 0]\}. \end{aligned}$$

Thus if $C_s = [C_1| \, C_2 \mid \ldots \mid C_n]$ then $| C_s | = 2^n$ where every super code word in $C_s$ takes entries from the two code words {(0 0 0 … 0), (1 1 1 1 … 1)}.

Thus we see in case of super special column code it is impossible to get repetition codes of different lengths. However this is possible as seen earlier using super special row codes.

Similarly we can get only parity check super special codes to have same length.

**DEFINITION 3.2.4:** *Let $C_s = [C_1 \mid C_2 \mid \ldots \mid C_n]^t$ be a super special parity check column code . Let the super special parity check column matrix associated with $C_s$ be given by*

$$H_s = \left[\begin{array}{c} H_1 \\ \hline H_2 \\ \hline \vdots \\ \hline H_n \end{array}\right]$$

*where*

$$H_1 = H_2 = \ldots = H_n = \underbrace{\left(1\ 1\ \ldots\ 1\right)}_{m-times}.$$

Thus we see we cannot get different lengths of parity check codes using the super special column code. However using super special row code we can get super special parity check codes of different lengths.

We illustrate super special parity check column codes.

***Example 3.2.6:*** Let $C_s = [C_1 | C_2 \mid C_3 \mid C_4]^t$ be a super special parity check column code with super special parity check matrix

$$H = \left[\begin{array}{c} H_1 \\ \hline H_2 \\ \hline H_3 \\ \hline H_4 \end{array}\right] = \left[\begin{array}{ccccc} 1 & 1 & 1 & 1 & 1 \\ \hline 1 & 1 & 1 & 1 & 1 \\ \hline 1 & 1 & 1 & 1 & 1 \\ \hline 1 & 1 & 1 & 1 & 1 \end{array}\right].$$

The codes associated with

$H_1$ is $C_1$ = {(0 0 0 0 0), (1 1 0 0 0), (0 1 1 0 0), (0 0 1 1 0), (0 0 0 1 1), (1 0 1 0 0), (1 0 0 1 0), (1 0 0 0 1), (0 1 0 1 0), (0 1 0 0 1), (0 0 1 0 1), (1 1 1 1 0), (1 1 1 0 1), (1 1 0 1 1), (1 0 1 1 1), (0 1 1 1 1)} = $C_2 = C_3 = C_4$.

Thus $C_s$ has $16 \times 16 \times 16 \times 16$ super special code words in it.

***Example 3.2.7:*** Let $C_s = [C_1 \mid C_2 \mid C_3]^t$ be a super special column parity check code. The super special column parity check matrix

$$H_s = \begin{bmatrix} H_1 \\ \hline H_2 \\ \hline H_3 \end{bmatrix} = \begin{bmatrix} 1 & 1 & 1 & 1 \\ \hline 1 & 1 & 1 & 1 \\ \hline 1 & 1 & 1 & 1 \end{bmatrix}.$$

The codes related with $H_1$ is $C_1$ = {(0 0 0 0), (1 1 0 0), (0 1 1 0), (0 0 1 1), (1 0 0 1), (1 0 1 0), (0 1 0 1), (1 1 1 1)} = $C_2$ = $C_3$. Thus $C_s$ contains 8 × 8 × 8 = 512 number of super special code words.

Thus only when the user wants to send messages generated by the same parity check matrix he can use it. However the main advantage of this special parity check column code $C_s = [C_1 \mid C_2 \mid \ldots \mid C_n]$; has each code word in $C_i$ which is of length m then $C_s$ has $n2^{m-1}$ number of code words hence one can use it in channels were one is not very much concerned with the transmission rate; for the transmission rate increases with increase in the length of the code words in $C_s$.

Now having seen this new class of codes using the parity check column code. We proceed on to built another new class of column codes using Hamming codes.

**DEFINITION 3.2.5:** *Let $C_s = [C_1 \mid C_2 \mid \ldots \mid C_n]^t$ be a super special column code if each of the codes $C_i$ in $C_s$ is a $(2^m - 1, 2^m - 1 - m)$ Hamming code for i = 1, 2, … , n then we call $C_s$ to be a super special column Hamming code . It is pertinent to mention that each code $C_i$ is a Hamming code of same length; i = 1, 2, … , n.*

Now we shall illustrate this by the following example.

***Example 3.2.8:*** Let $C_s = [C_1 \mid C_2 \mid C_3]^t$ be a super special column Hamming code.

Suppose

$$H_s = \begin{bmatrix} H_1 \\ \hline H_2 \\ \hline H_3 \end{bmatrix}$$

be the associated super special parity check column matrix of $C_s$, where

$$H_1 = H_2 = H_3 = \begin{bmatrix} 0 & 0 & 0 & 1 & 1 & 1 & 1 \\ 0 & 1 & 1 & 0 & 0 & 1 & 1 \\ 1 & 0 & 1 & 0 & 1 & 0 & 1 \end{bmatrix}.$$

The code

$C_i$ = {(0 0 0 0 0 0 0), (1 0 0 0 0 1 1), (0 1 0 0 1 0 1),
(0 0 1 0 1 1 0), (0 0 0 1 1 0 1), (1 1 0 0 1 0 1), (0 1 1 0 0 1 1),
(0 0 1 1 1 0 0), (1 0 1 0 1 0 1), (1 0 0 1 1 0 0), (0 1 0 1 0 1 0),
(1 1 1 0 0 0 0), (0 1 1 1 1 0 0), (1 1 0 1 0 1 1),
(1 0 1 1 1 1 1), (1 1 1 1 1 1 1)};

i = 1, 2, 3. Thus $C_s$ has $16^3$ super special code words in it.

We can yet have different codes using different parity check matrices.

We just illustrate this by the following example.

***Example 3.2.9:*** Let $C_s = [C_1 \mid C_2]^t$ be a super special Hamming column code of length seven.

Suppose

$$H_s = \begin{bmatrix} H_1 \\ \hline H_2 \end{bmatrix} = \begin{bmatrix} 0 & 0 & 0 & 1 & 1 & 1 & 1 \\ 0 & 1 & 1 & 0 & 0 & 1 & 1 \\ 1 & 0 & 1 & 0 & 1 & 0 & 1 \\ \hline 1 & 1 & 1 & 0 & 1 & 0 & 0 \\ 0 & 1 & 1 & 1 & 0 & 1 & 0 \\ 0 & 0 & 1 & 1 & 1 & 0 & 1 \end{bmatrix}$$

be the associated super special parity check column matrix of $C_s$. Then

$C_1$ = {(0 0 0 0 0 0 0), (1 0 0 0 0 1 1), (0 1 0 0 1 0 1),
(0 0 1 0 1 1 0), (0 0 0 1 1 0 1), (1 1 0 0 1 0 1), (0 1 1 0 0 1 1),
(0 0 1 1 1 0 0), (1 0 1 0 1 0 1), (1 0 0 1 1 0 0), (0 1 0 1 0 1 0),
(1 1 1 0 0 0 0), (0 1 1 1 1 0 0), (1 1 0 1 0 1 1),
(1 0 1 1 1 1 1), (1 1 1 1 1 1 1)}

and

$C_2$ = {(0 0 0 0 0 0 0), (1 0 0 0 1 0 1), (0 1 0 0 1 1 1),
(0 0 1 0 1 1 0), (0 0 0 1 0 1 1), (1 1 0 0 0 1 0), (0 1 1 0 0 0 1),
(0 0 1 1 1 0 1), (1 0 1 0 0 1 1), (1 0 0 1 1 1 1), (0 1 0 1 1 0 0),
(1 1 1 0 1 0 1), (0 1 1 1 0 1 0), (1 1 0 1 0 0 1),
(1 0 1 1 0 0 0), (1 1 1 1 1 1 1)}.

Thus we see $C_1$ and $C_2$ are different Hamming codes but $C_1 \cap C_2 \neq \phi$ as well as $C_1 \neq C_2$.

So we can get different sets of codes and the number of elements in $C_s$ is 256.

Now as in case of super special row codes we can in case of super special column codes have mixed super special column codes. The only criteria being is that each code $C_i$ in $C_s$ will be of same length.

Now we proceed on to define them.

**DEFINITION 3.2.6:** *$C_s = [C_1 \mid C_2 \mid \ldots \mid C_n]^t$ is a mixed super special column code if some $C_i$'s are repetition codes of length n some $C_j$'s are Hamming codes of length n, some $C_k$'s parity check codes of length n and others are arbitrary codes, $1 \leq i, j, k \leq n$.*

We illustrate this by the following example.

***Example 3.2.10:*** Let $C_s = [C_1 \mid C_2 \mid C_3 \mid C_4]^t$ be a mixed super special column code. Let

$$H_s = \left[\begin{array}{c} H_1 \\ \hline H_2 \\ \hline H_3 \\ \hline H_4 \end{array}\right]$$

be the super special parity check column matrix associated with $C_s$. Here

$$H_s = \begin{bmatrix} H_1 \\ \hline H_2 \\ \hline H_3 \\ \hline H_4 \end{bmatrix} = \begin{bmatrix} 1 & 1 & 1 & 1 & 1 & 1 & 1 \\ \hline 1 & 1 & 0 & 0 & 0 & 0 & 0 \\ 1 & 0 & 1 & 0 & 0 & 0 & 0 \\ 1 & 0 & 0 & 1 & 0 & 0 & 0 \\ 1 & 0 & 0 & 0 & 1 & 0 & 0 \\ 1 & 0 & 0 & 0 & 0 & 1 & 0 \\ 1 & 0 & 0 & 0 & 0 & 0 & 1 \\ \hline 0 & 0 & 0 & 1 & 1 & 1 & 1 \\ 0 & 1 & 1 & 0 & 0 & 1 & 1 \\ 1 & 0 & 1 & 0 & 1 & 0 & 1 \\ \hline 1 & 0 & 0 & 1 & 0 & 0 & 0 \\ 0 & 1 & 1 & 0 & 1 & 0 & 0 \\ 1 & 0 & 1 & 0 & 0 & 1 & 0 \\ 1 & 1 & 1 & 0 & 0 & 0 & 1 \end{bmatrix};$$

here $H_1$ is a parity check code of length 7, $H_2$ is a repetition code with 6 check symbols, $H_3$ is the Hamming code of length 7 and $H_4$ a code with 3 message symbols. The transmission rate of

$$C_s = \frac{6+1+4+3}{7+7+7+7} = \frac{14}{28} = \frac{1}{2}.$$

Now we proceed on to define the super special column cyclic code.

**DEFINITION 3.2.7:** *Let $C_s = [C_1 \mid C_2 \mid \ldots \mid C_m]^t$ be a super special column code if each of the codes $C_i$ is a cyclic code then we call $C_s$ to be a super special cyclic column code. However the length of each code $C_i$; i = 1, 2, …, n will only be a cyclic code of length n, but the number of message symbols and check symbols can be anything.*

Now we illustrate this by an example.

***Example 3.2.11:*** Let $C_s = [C_1 \mid C_2 \mid C_3]^t$ be a super special cyclic column code where each $C_i$ is of length six.

Now the associated super special parity check column matrix

$$H_s = \begin{bmatrix} H_1 \\ \hline H_2 \\ \hline H_3 \end{bmatrix} = \left[\begin{array}{cccccc} 0 & 0 & 1 & 0 & 0 & 1 \\ 0 & 1 & 0 & 0 & 1 & 0 \\ 1 & 0 & 0 & 1 & 0 & 0 \\ \hline 1 & 1 & 0 & 0 & 0 & 0 \\ 1 & 0 & 1 & 0 & 0 & 0 \\ 1 & 0 & 0 & 1 & 0 & 0 \\ 1 & 0 & 0 & 0 & 1 & 0 \\ 1 & 0 & 0 & 0 & 0 & 1 \\ \hline 0 & 0 & 0 & 0 & 1 & 1 \\ 0 & 0 & 0 & 1 & 1 & 0 \\ 0 & 0 & 1 & 1 & 0 & 0 \\ 0 & 1 & 1 & 0 & 0 & 0 \\ 1 & 1 & 0 & 0 & 0 & 0 \end{array}\right].$$

The cyclic codes given by the parity check matrix $H_1$ is

$$C_1 = \{(0\ 0\ 0\ 0\ 0\ 0), (1\ 0\ 0\ 1\ 0\ 0), (0\ 1\ 0\ 0\ 1\ 0), (0\ 0\ 1\ 0\ 0\ 1), (1\ 1\ 0\ 1\ 1\ 0), (0\ 1\ 1\ 0\ 1\ 1), (1\ 0\ 1\ 1\ 0\ 1), (1\ 1\ 1\ 1\ 1\ 1)\},$$

$$C_2 = \{(1\ 1\ 1\ 1\ 1\ 1), (0\ 0\ 0\ 0\ 0\ 0)\}$$

and

$$C_3 = \{(0\ 0\ 0\ 0\ 0\ 0), (1\ 1\ 1\ 1\ 1\ 1)\}.$$

Thus $C_s$ contains $8 \times 2 \times 2 = 32$ elements.

Now in the mixed super special column code $C_s = [C_1 \mid C_2 \mid \ldots \mid C_n]^t$ code we can have some of the $C_i$'s to be cyclic codes also. Now in case of super special column code for any given $C_s$ if we have a super special parity check matrix $H_s$ to be in the standard form we can always get $G_s$ and we have $G_s H_s^T$ = a zero super special column matrix.

## 3.3 New Classes of Super Special Codes

We have given the basic definition and properties of super matrices in chapter one. In this section we proceed on to define new classes of supper special codes and discuss a few properties about them.

**DEFINITION 3.3.1:** *Let*

$$C(S) = \left[\begin{array}{c|c|c|c} C_1^1 & C_2^1 & \cdots & C_n^1 \\ \hline C_1^2 & C_2^2 & \cdots & C_n^2 \\ \hline \vdots & \vdots & \cdots & \vdots \\ \hline C_1^m & C_2^m & \cdots & C_n^m \end{array}\right]$$

*where $C_j^i$ are codes $1 \leq i \leq m$ and $1 \leq j \leq n$. Further all codes $C_1^1, C_1^2, \ldots, C_1^m$ are of same length $C_2^1, C_2^2, \ldots, C_2^m$ are of same length and $C_n^1, C_n^2, \ldots, C_2^m$ are of same length. $C_1^1, C_2^1, \ldots, C_n^1$ have same number of check symbols, $C_1^2, C_2^2, \ldots, C_n^2$ have same number of check symbols and $C_1^m, C_2^m, \ldots, C_n^m$ have same number of check symbols.*

*We call C(S) to be a super special code. We can have the super parity check matrix*

$$H(S) = \left[\begin{array}{c|c|c|c} H_1^1 & H_2^1 & \cdots & H_n^1 \\ \hline H_1^2 & H_2^2 & \cdots & H_n^2 \\ \hline \vdots & \vdots & \cdots & \vdots \\ \hline H_1^m & H_2^m & \cdots & H_n^m \end{array}\right]$$

*where $H_j^i$'s are parity check matrices $1 \leq i \leq m$ and $1 \leq j \leq n$. Further $H_1^1, H_2^1, \ldots, H_n^1$ have the same number of rows, $H_1^2, H_2^2, \ldots, H_n^2$ have same number of rows and so on. $H_1^m, H_2^m, \ldots, H_n^m$ have the same number of rows. Likewise $H_1^1, H_1^2, \ldots, H_1^m$ have the same number of columns, $H_2^1, H_2^2, \ldots, H_2^m$ have the same number of columns and so on. $H_n^1, H_n^2, \ldots, H_n^m$ have same number of columns.*

Now this super special code has two types of messages i.e.,

1. Array of super row vector messages i.e.,

$$\begin{array}{c}\left[C_1^1 \middle| C_2^1 \middle| \ldots \middle| C_n^1\right] \\ \left[C_1^2 \middle| C_2^2 \middle| \ldots \middle| C_n^2\right] \\ \vdots \quad \vdots \\ \left[C_1^m \middle| C_2^m \middle| \ldots \middle| C_n^m\right];\end{array}$$

we have m rows of super row vectors and one by one rows are sent at the receiving end; one after decoding puts it only in the form of the same array of rows.

2. Array of super column vector messages

$$\begin{bmatrix} C_1^1 \\ C_1^2 \\ \vdots \\ C_1^m \end{bmatrix} \begin{bmatrix} C_2^1 \\ C_2^2 \\ \vdots \\ C_2^m \end{bmatrix} \ldots \begin{bmatrix} C_n^1 \\ C_n^2 \\ \vdots \\ C_n^m \end{bmatrix}.$$

Now we have n columns of super column vectors and these after taking transposes of each column the messages are sent one by one as

$\begin{bmatrix} C_1^1 \\ C_1^2 \\ \vdots \\ C_1^m \end{bmatrix}^t$ which is a super row vector, $\begin{bmatrix} C_1^2 \\ C_2^2 \\ \vdots \\ C_2^m \end{bmatrix}^t$ a super row vector and so on. At the receiving end once again they are arranged back as the array of the column vectors which will be termed as the received message. Thus we can have transmission of two types at the source and the received message will accordingly be of two types viz. array of super row vectors or array of super column vectors.

We will illustrate this by the following example.

***Example 3.3.1:*** Let

$$C(S) = \left[\begin{array}{c|c} C_1^1 & C_2^1 \\ \hline C_1^2 & C_2^2 \\ \hline C_1^3 & C_2^3 \end{array}\right]$$

be a super special code.
Suppose

$$H(S) = \left[\begin{array}{c|c} H_1^1 & H_2^1 \\ \hline H_1^2 & H_2^2 \\ \hline H_1^3 & H_2^3 \end{array}\right]$$

be the super special parity check matrix associated with C(S).

$$H(S) = \left[\begin{array}{cccccc|ccccccc} 0 & 0 & 1 & 1 & 0 & 0 & 1 & 0 & 0 & 1 & 1 & 0 & 0 \\ 0 & 1 & 1 & 0 & 1 & 0 & 0 & 1 & 0 & 1 & 0 & 1 & 0 \\ 1 & 1 & 1 & 0 & 0 & 1 & 1 & 1 & 1 & 0 & 0 & 0 & 1 \\ \hline 1 & 0 & 1 & 0 & 0 & 0 & 1 & 1 & 1 & 1 & 0 & 0 & 0 \\ 1 & 1 & 0 & 1 & 0 & 0 & 0 & 1 & 1 & 0 & 1 & 0 & 0 \\ 0 & 1 & 0 & 0 & 1 & 0 & 1 & 0 & 1 & 0 & 0 & 1 & 0 \\ 1 & 0 & 0 & 0 & 0 & 1 & 1 & 1 & 0 & 0 & 0 & 0 & 1 \\ \hline 1 & 0 & 0 & 1 & 0 & 0 & 1 & 1 & 0 & 1 & 1 & 0 & 0 \\ 1 & 1 & 0 & 0 & 1 & 0 & 0 & 1 & 1 & 0 & 0 & 1 & 0 \\ 1 & 0 & 1 & 0 & 0 & 1 & 1 & 1 & 1 & 1 & 0 & 0 & 1 \end{array}\right].$$

Now using the parity check matrices one can get the codes associated with each $H_t^i$; $1 \leq i \leq 3$ and $1 \leq t \leq 2$.

A typical super code word is a block

$$\left[\begin{array}{cccccc|ccccccc} 1 & 0 & 0 & 0 & 0 & 1 & 0 & 1 & 0 & 0 & 1 & 0 & 1 \\ \hline 1 & 1 & 1 & 0 & 1 & 1 & 1 & 0 & 1 & 0 & 1 & 0 & 1 \\ \hline 1 & 1 & 1 & 1 & 0 & 0 & 1 & 1 & 1 & 1 & 1 & 0 & 0 \end{array}\right]$$

if we take it as array of super row vectors we have

$$[1 \quad 0 \quad 0 \quad 0 \quad 0 \quad 1 \,|\, 0 \quad 1 \quad 0 \quad 0 \quad 1 \quad 0 \quad 1],$$
$$[1 \quad 1 \quad 1 \quad 0 \quad 1 \quad 1 \,|\, 1 \quad 0 \quad 1 \quad 0 \quad 1 \quad 0 \quad 1],$$
$$[1 \quad 1 \quad 1 \quad 1 \quad 0 \quad 0 \,|\, 1 \quad 1 \quad 1 \quad 1 \quad 1 \quad 0 \quad 0].$$

As an array of the super column vector we have

$$\left[\begin{array}{c} 1\\0\\0\\0\\0\\1\\ \hline 1\\1\\1\\0\\1\\1\\ \hline 1\\1\\1\\1\\0 \end{array}\right], \quad \left[\begin{array}{c} 0\\1\\0\\0\\1\\0\\1\\ \hline 1\\0\\1\\0\\1\\0\\1\\ \hline 1\\1\\1\\1\\1\\0\\0 \end{array}\right].$$

When we send it as message the array of super row vectors it would be sent as [1 0 0 0 0 1 | 0 1 0 0 1 0 1], [1 1 1 1 0 0 | 1 1 1

1 1 0 0] and [1 1 1 0 0 1 | 1 1 1 1 1 0 0] after receiving the message; the received message would be given the array super row representation. While sending the array of super columns we send the message as

$$\left[\begin{array}{c}1\\0\\0\\0\\0\\1\\\hline 1\\1\\1\\0\\1\\1\\\hline 1\\1\\1\\1\\0\end{array}\right]^t, \quad \left[\begin{array}{c}0\\1\\0\\0\\1\\0\\1\\\hline 1\\0\\1\\0\\1\\0\\1\\\hline 1\\1\\1\\1\\0\\0\end{array}\right]^t$$

that is as [1 0 0 0 0 1 | 1 1 1 0 1 1 | 1 1 1 1 0], [0 1 0 0 1 0 1 | 1 0 1 0 1 0 1| 1 1 1 1 0 0] after receiving the message it is given the array super column representation.

Now we would show how the received message is verified or found to be correct one or not by the following method. Though one can have other methods also to check whether the received message is correct or not. Now we first describe the method using the general case before we proceed to explain with specific examples. Suppose

$$C(S) = \left[\begin{array}{c|c|c|c} C_1^1 & C_2^1 & \cdots & C_n^1 \\ \hline C_1^2 & C_2^2 & \cdots & C_n^2 \\ \hline \vdots & \vdots & \cdots & \vdots \\ \hline C_1^m & C_2^m & \cdots & C_n^m \end{array}\right]$$

be the super special super code. Any super special super code word of C(S) say x(S) ∈ C(S) would be of the form i.e.,

$$x(S) = \left[\begin{array}{c|c|c|c} x_1^1 & x_2^1 & \cdots & x_n^1 \\ \hline x_1^2 & x_2^2 & \cdots & x_n^2 \\ \hline \vdots & \vdots & \cdots & \vdots \\ \hline x_1^m & x_2^m & \cdots & x_n^m \end{array}\right]$$

where each $x_i^j$ is a row vector; $1 \le j \le m$ and $1 \le i \le n$.
Suppose

$$H(S) = \left[\begin{array}{c|c|c|c} H_1^1 & H_2^1 & \cdots & H_n^1 \\ \hline H_1^2 & H_2^2 & \cdots & H_n^2 \\ \hline \vdots & \vdots & \cdots & \vdots \\ \hline H_1^m & H_2^m & \cdots & H_n^m \end{array}\right]$$

be the super special parity check matrix associated with C(S). Now R(S) be the received super special super code word given by

$$R(S) = \left[\begin{array}{c|c|c|c} y_1^1 & y_2^1 & \cdots & y_n^1 \\ \hline y_1^2 & y_2^2 & \cdots & y_n^2 \\ \hline \vdots & \vdots & \cdots & \vdots \\ \hline y_1^m & y_2^m & \cdots & y_n^m \end{array}\right].$$

To check whether the received super code word R(S) is correct or not, we make use of the super special syndrome technique.

We define super special syndrome of C(S) as

$$S\left(x\left(S\right)\right) = H\left(S\right)\left[x\left(S\right)\right]^t =$$

$$\left[\begin{array}{c|c|c|c} H_1^1\left(x_1^1\right)^t & H_2^1\left(x_2^1\right)^t & \cdots & H_n^1\left(x_n^1\right)^t \\ \hline H_1^2\left(x_1^2\right)^t & H_2^2\left(x_2^2\right)^t & \cdots & H_1^2\left(x_1^2\right)^t \\ \hline \vdots & \vdots & \cdots & \vdots \\ \hline H_1^m\left(x_1^m\right)^t & H_2^m\left(x_2^m\right)^t & \cdots & H_n^m\left(x_n^m\right)^t \end{array}\right] = \left[\begin{array}{c|c|c|c} (0) & (0) & \cdots & (0) \\ \hline (0) & (0) & \cdots & (0) \\ \hline \vdots & \vdots & \cdots & \vdots \\ \hline (0) & (0) & \cdots & (0) \end{array}\right]$$

then we take x (S) $\in$ C(S); if

$$H(S)\ [x(S)]^t \neq \left[\begin{array}{c|c|c|c} (0) & (0) & \cdots & (0) \\ \hline (0) & (0) & \cdots & (0) \\ \hline \vdots & \vdots & \cdots & \vdots \\ \hline (0) & (0) & \cdots & (0) \end{array}\right]$$

then we declare x(S) $\notin$ C(S) so if y(S) is the received message one calculates

$$H(S)\ [y(S)]^t = \left[\begin{array}{c|c|c|c} (0) & (0) & \cdots & (0) \\ \hline (0) & (0) & \cdots & (0) \\ \hline \vdots & \vdots & \cdots & \vdots \\ \hline (0) & (0) & \cdots & (0) \end{array}\right],$$

then y(S) is a code word of C(S); other wise we can declare the received message has an error in it. We can find $S\left(y_j^i\right) = H_j^i\left(y_j^i\right)^t$; $1 \leq i \leq m$; $1 \leq j \leq n$. If $S\left(y_j^i\right) \neq 0$ then we use the technique of coset leader to find the error. We can also use the method of best approximations and find the approximate sent message.

**Chapter Four**

# APPLICATIONS OF THESE NEW CLASSES OF SUPER SPECIAL CODES

We enumerate in this chapter a few applications and advantages of using the super special codes. Now the following will show why this super special super code is better than other usual code. These codes can also be realised as a special type of concatenated codes. We enumerate them in the following:

1. Instead of using ARQ protocols we can use the same code C, in the super special super code

$$\left[\begin{array}{c|c|c|c} C_1^1 & C_2^1 & \cdots & C_n^1 \\ \hline C_1^2 & C_2^2 & \cdots & C_n^2 \\ \hline \vdots & \vdots & \cdots & \vdots \\ \hline C_1^m & C_2^m & \cdots & C_n^m \end{array}\right];$$

replace every $C_i^j$ by C, m can be equal to n or greater than n or less than n. So that if the same message is sent from

$$C(S) = \left[\begin{array}{c|c|c|c} C & C & \cdots & C \\ \hline C & C & \cdots & C \\ \hline \vdots & \vdots & \cdots & \vdots \\ \hline C & C & \cdots & C \end{array}\right].$$

We can take the correct message. This saves both money and time. We can also send it is as array of super row vectors i.e., if $x$ is the message to be sent then,

$$\begin{array}{c} [x \mid x \mid \ldots \mid x] \\ [x \mid x \mid \ldots \mid x] \\ \vdots \\ [x \mid x \mid \ldots \mid x]; \end{array}$$

As an array of row super vector or the array of the column super vectors as

$$\left[\begin{array}{c} x \\ \hline x \\ \hline \vdots \\ \hline x \end{array}\right] \left[\begin{array}{c} x \\ \hline x \\ \hline \vdots \\ \hline x \end{array}\right] \cdots \left[\begin{array}{c} x \\ \hline x \\ \hline \vdots \\ \hline x \end{array}\right]$$

where

$$\left[\begin{array}{c} x \\ \hline x \\ \hline \vdots \\ \hline x \end{array}\right]^t$$

is sent. In case the receiver wants to get the very correct message he can transmit both as an array of super row codes as well as, an array of super column codes the same code $x$ and get the correct sent message.

2. This super special super code C(S) has another advantage for if

$$C(S) = \left[\begin{array}{c|c|c|c} C & C & \cdots & C \\ \hline C & C & \cdots & C \\ \hline \vdots & \vdots & \cdots & \vdots \\ \hline C & C & \cdots & C \end{array}\right]$$

then if one wishes to study the changes in terms of time the same message can be sent in all cells and the gradual stage by stage transformation can be seen (observed) and the resultant can be got.

Here it is not sending a message and receiving a message but is a study of transformation from time to time. It can be from satellite or pictures of heavenly bodies. Even in medical field this will find an immense use. Also these types of super special super codes can be used in scientific experiments so that changes can be recorded very minutely or with high sensitivity. What one needs is a proper calibration linking these codes with those experiments were one is interested in observing the changes from time to time were the graphical representation is impossible due to more number of variables.

3. Another striking advantage of the super special super code C(S) is that if one has to be doubly certain about the accuracy of the received message or one cannot request for second time transmission in those cases the sender can send super array of row codes and send the same codes as the super array of column codes.

After receiving both the messages, if both the received codes are identical and without any error one can accept it; otherwise find the error in each cell and accept the message which has less number of errors. If both the received messages have the same number of errors then if the machine which sends the code has a provision for a single or two or any desired number of cells alone can be non empty and other cells have empty code and send the message; in that case the super special code is taken as

$$C(S) \cup \left[\begin{array}{c|c|c|c} \phi & \phi & \cdots & \phi \\ \hline \phi & \phi & \cdots & \phi \\ \hline \vdots & \vdots & \cdots & \vdots \\ \hline \phi & \phi & \cdots & \phi \end{array}\right]$$

or

$$C(S) = \left[C_1^i \mid C_2^i \mid \ldots \mid C_n^i\right]$$

itself has the empty code word $\phi$ in each of the $C_j^i$'s in C(S) for i = 1, 2, …, n and $1 \leq j \leq m$.

4. At times the super special code C(S) may be of a special type i.e., if

$$C(S) = \left[\begin{array}{c|c|c|c} C_1^1 & C_2^1 & \cdots & C_n^1 \\ \hline C_1^2 & C_2^2 & \cdots & C_n^2 \\ \hline \vdots & \vdots & \cdots & \vdots \\ \hline C_1^m & C_2^m & \cdots & C_n^m \end{array}\right]$$

is such that

$$C_1^1 = C_1^2 = \ldots C_1^m, C_2^1 = C_2^2 = \ldots C_2^m, \ldots, C_n^1 = C_n^2 = \ldots C_n^m$$

i.e., each row in C(S) i.e.,

$$= \left[C_1^1 \mid C_2^1 \mid \ldots \mid C_n^1\right] =$$

$$\left[C_1^2 \mid C_2^2 \mid \ldots \mid C_n^2\right] = \ldots =$$

$$\left[C_1^m \mid C_2^m \mid \ldots \mid C_n^m\right];$$

i.e.,

$$C(S) = \left[\begin{array}{c|c|c|c} C_1^1 & C_2^1 & \cdots & C_n^1 \\ \hline C_1^1 & C_2^1 & \cdots & C_n^1 \\ \hline \vdots & \vdots & \cdots & \vdots \\ \hline C_1^1 & C_2^1 & \cdots & C_n^1 \end{array}\right]$$

so that any one message is sent one can test if not all the super row vectors received as a message is not identical one can guess the error has occurred during transmission. The maximum number of super row vectors which happen to repeat would be accepted as the approximately correct message. If there is not even a single pair of coinciding super row vectors then we choose a super row which has a minimum number of errors. i.e., if

$$H(S)=\left[\begin{array}{c|c|c|c} H_1^1 & H_2^1 & \cdots & H_n^1 \\ \hline H_1^1 & H_2^1 & \cdots & H_n^1 \\ \hline \vdots & \vdots & \cdots & \vdots \\ \hline H_1^1 & H_2^1 & \cdots & H_n^1 \end{array}\right]$$

then we say $i^{th}$ super row has least number of errors if

$$x^i(s) = \left[x_1^1 \mid x_2^1 \mid \ldots \mid x_n^1\right]$$

we find out

$$H_1^1\left(x_1^1\right)^t, H_2^1\left(x_2^1\right)^t, \ldots, H_n^1\left(x_n^1\right)^t;\ 1 \le i \le m;$$

we choose the super row which has the maximum number of zeros i.e., that super row consequently has minimum number of error. We find the correct code word from those cells in the $i^{th}$ row and accept it as the approximately correct received row. If already we have a super row in which $H_1^1\left(x_1^1\right)^t = (0)$ for $1 \le i \le n$ then we accept that as the correct message. This form of transformation helps the receiver to study the real error pattern and the cells in C(S) which misbehave or that which always has an error message. Thus one can know not only more about the sent message but also know more about the problems (in the machine) while the message is transmitted; consequently corrections can be made so as to guarantee an approximately correct message is received.

5. Suppose we have a super special super code C(S) where

$$C(S)=\left[\begin{array}{c|c|c|c} C_1^1 & C_1^1 & \cdots & C_1^1 \\ \hline C_2^1 & C_2^1 & \cdots & C_2^1 \\ \hline \vdots & \vdots & \cdots & \vdots \\ \hline C_m^1 & C_m^1 & \cdots & C_m^1 \end{array}\right]$$

with n columns where $C_j^1$ are linear codes $1 \leq j \leq m$; i.e., we have only m distinct codes filled in the n column in a way shown above. The transmission can also take place in two ways simultaneously or by either array of super row transmission alone or array of super column transmission; only when we say the simultaneous transmission, we will first send the message x(S) in array of super row vector

$$= \begin{matrix} \left[x_1^1 \mid x_1^1 \mid \ldots \mid x_1^1\right] \\ \left[x_2^1 \mid x_2^1 \mid \ldots \mid x_2^1\right] \\ \vdots \quad \vdots \quad \quad \vdots \\ \left[x_m^1 \mid x_m^1 \mid \ldots \mid x_n^1\right] \end{matrix}$$

and then send the same x(S) in the array of super column vector as

$$\begin{bmatrix} x_1^1 \\ \hline x_2^1 \\ \hline \vdots \\ \hline x_m^1 \end{bmatrix} \begin{bmatrix} x_1^1 \\ \hline x_2^1 \\ \hline \vdots \\ \hline x_m^1 \end{bmatrix} \cdots \begin{bmatrix} x_1^1 \\ \hline x_2^1 \\ \hline \vdots \\ \hline x_m^1 \end{bmatrix}$$

and the transmission is made as

$$\begin{bmatrix} x_1^1 \\ \hline x_2^1 \\ \hline \vdots \\ \hline x_m^1 \end{bmatrix}^t \begin{bmatrix} x_1^1 \\ \hline x_2^1 \\ \hline \vdots \\ \hline x_m^1 \end{bmatrix}^t \cdots \begin{bmatrix} x_1^1 \\ \hline x_2^1 \\ \hline \vdots \\ \hline x_m^1 \end{bmatrix}^t$$

i.e., when we send both super column vector as well as the super row vector we call it as the simultaneous transmission here

$$x(S) = \left[\begin{array}{c|c|c|c} x_1^1 & x_1^1 & \cdots & x_1^1 \\ \hline x_2^1 & x_2^1 & \cdots & x_2^1 \\ \hline \vdots & \vdots & \cdots & \vdots \\ \hline x_m^1 & x_m^1 & \cdots & x_m^1 \end{array}\right].$$

When we send only array super row transmission, every received code word in the i$^{th}$ row (say) $\left[y_i^1 \mid y_i^1 \mid \ldots \mid y_i^1\right]$ must be the same. This is true for i = 1, 2, …, m. If they are different for any row just by observation we can conclude the received message has an error and choose the row which has least number of differences.

Now when array of column transmission takes place we see if y(S) is the received message then

$$y(S) = \begin{matrix} \left[y_1^1 \mid y_2^1 \mid \ldots \mid y_m^1\right]^t \\ \left[y_1^1 \mid y_2^1 \mid \ldots \mid y_m^1\right]^t \\ \vdots \\ \left[y_1^1 \mid y_2^1 \mid \ldots \mid y_m^1\right]^t \end{matrix} .$$

We take that column which has least number of errors or which has no error as the received message. The main advantage of simultaneous transmission is we can compare each cell of the received array of super row vectors and array of super column vectors.

6. This type of code can be used when the ARQ process is impossible so that the same message can be filled in each cell of the super special code C(S) i.e.,

$$x(S) = \left[\begin{array}{c|c|c|c} x & x & \cdots & x \\ \hline x & x & \cdots & x \\ \hline \vdots & \vdots & \cdots & \vdots \\ \hline x & x & \cdots & x \end{array}\right]$$

where x is the code word in C, which is in every cell of C(S) so that the approximate correct message is always retrieved without any difficulty. That is if we have say for example a code C of length 8 with 16 code words only then we can choose C(S) to be a super special code with 16 rows and 17 columns.

$$16 \text{ rows}\left\{\underbrace{\left[\begin{array}{c|c|c|c} C & C & \cdots & C \\ \hline C & C & \cdots & C \\ \hline \vdots & \vdots & \cdots & \vdots \\ \hline C & C & \cdots & C \end{array}\right]}_{17 \text{ columns}}\right.,$$

thus if x is sent message then

$$x(S) = \left[\begin{array}{c|c|c|c} x & x & \cdots & x \\ \hline x & x & \cdots & x \\ \hline \vdots & \vdots & \cdots & \vdots \\ \hline x & x & \cdots & x \end{array}\right]$$

so when

$$y(S) = \left[\begin{array}{c|c|c|c} y_1 & y_2 & \cdots & y_{17} \\ \hline y_1 & y_2 & \cdots & y_{17} \\ \hline \vdots & \vdots & \cdots & \vdots \\ \hline y_{16} & y_{16} & \cdots & y_{16} \end{array}\right]$$

is the received super special code then we see if every cell element in y(S) is the same or that element in y(S) which repeats itself is taken as the approximate correct message and ARQ protocols can be avoided thus saving both time and money. Also when the number of cells in C(S) is increased and is greater than that of the number of elements in the space $V^n$ where C is the code then we can easily be guaranteed that the same message is sent mn times (where m > n) we are sure to retrieve the correct message. When these super special codes C(S) are used the user is certain to get marvellous results be it in any discipline.

7. The class of super special super codes C(S) will be very beneficial in the cryptography for two main reasons.

a. It is easy to be operated or transmitted and

b. The intruder can be very easily misled and his guess can never give him/her the true transmitted message. We just show

how the super special code C(S) functions so that the intruder can never guess the same or break the message.
Suppose

$$C(S) = \left[\begin{array}{c|c|c|c|c} C_1^1 & C_2^1 & C_3^1 & \cdots & C_n^1 \\ \hline C_1^2 & C_2^2 & C_3^2 & \cdots & C_n^2 \\ \hline \vdots & \vdots & \vdots & \cdots & \vdots \\ \hline C_1^1 & C_2^1 & C_3^1 & \cdots & C_n^1 \\ \hline C_1^{m+1} & C_2^{m+1} & C_3^{m+1} & \cdots & C_n^{m+1} \\ \hline \vdots & \vdots & \vdots & \cdots & \vdots \\ \hline C_1^1 & C_2^1 & C_3^1 & \cdots & C_n^1 \end{array}\right].$$

The real message carrying code is say in 1st row, (m + 1)th row and so on and the last row. Very different codes which do not carry the message will also be repeated so that even by the frequency of the repetition intruder even cannot guess. Only the concerned who are the reliable part and parcel of the communication work knows the exact super rows which carry the messages so they would only look into that super row to guarantee the error freeness during transmission. Several super rows carry the same message. In our opinion this method with larger number of rows and columns in C(S) would make it impossible for the intruder to break it.

We give yet another super special code C(S), which is impossible for any intruder to break and which guarantees the maximum security.

8. Let

$$C(S) = \left[\begin{array}{c|c|c|c|c|c|c|c} C_1^1 & C_1^1 & C_1^1 & \cdots & C_2^1 & C_2^1 & \cdots & C_m^1 \\ \hline C_2^2 & C_2^2 & C_2^2 & \cdots & C_2^2 & C_2^2 & \cdots & C_m^2 \\ \hline \vdots & \vdots & \vdots & \cdots & \vdots & \vdots & \cdots & \vdots \\ \hline C_t^t & C_t^t & C_t^t & \cdots & C_{t+1}^t & C_{t+1}^t & \cdots & C_m^t \\ \hline \vdots & \vdots & \vdots & \cdots & \vdots & \vdots & \cdots & \vdots \\ \hline C_r^n & C_r^n & C_r^n & \cdots & C_{r+1}^n & C_{r+1}^n & \cdots & C_m^n \end{array}\right]$$

that C(S) has some number of codes repeated a few times (say $p_1$) then another code repeating (say some $p_2$) times and some other code repeating $p_3$ times and so on finally yet another code repeating some $p_r$ times. This is done for each and every row. Only the concerned persons who are part and parcel of the group know at what stages and which codes carry the message and such keys are present only with the group who exchange the information so it is impossible for any intruder to break it and high percentage of confidentiality and security is maintained.

We give yet another type of super special code C(S).

9. Let C(S) be a super special super code. Suppose there are n codes arranged in the column and m codes along the rows of the super matrix of the code C(S). Now the m × n codes are arranged in special super blocks where by a block we mean a p × q array of same code i.e., say if C is the code then the special super block has

$$\text{p-rows}\left\{\underbrace{\left[\begin{array}{c|c|c|c} C & C & \cdots & C \\ \hline C & C & \cdots & C \\ \hline \vdots & \vdots & \cdots & \vdots \\ \hline C & C & \cdots & C \end{array}\right]}_{\text{q columns}}\right. .$$

Thus the code C(S) has

$$\begin{array}{l} r_1^1 \times s_1^1, r_1^2 \times s_1^1, \ldots, r_1^{q_1} \times s_1^{q_1}, \\ r_2^1 \times s_2^1, r_2^2 \times s_2^2, \ldots, r_2^{q_2} \times s_2^{q_2}, \\ \ldots, \\ r_{p_1}^1 \times s_{p_1}^1, r_{p_2}^2 \times r_{p_2}^2, \ldots, r_{p_s}^{q_s} \times s_{p_s}^{q_s} \end{array}$$

block such that $r_1^1 + r_2^1 + \ldots + r_{p_1}^1 = m$, $s_1^1 + s_1^2 + \ldots + s_1^{q_1} = n$ and so on $r_{p_1}^1 + r_{p_2}^2 + \ldots + r_{p_s}^{q_s} = m$ and $s_{p_1}^1 + s_{p_2}^2 + \ldots + s_{p_s}^{q_s} = n$.

Now each of these blocks contain the same code i.e., codes can vary only with the varying blocks.

The cryptographer can choose some blocks in C(S) to carry the messages and rest of the blocks may be used to mislead the intruder. When this type of super special codes are used it is impossible for any one to break and get into the structure.

Now even in this block the cryptographer can use only certain rows and columns to carry true message and the rest only to mislead the intruder.

We will illustrate these by some simple examples.

***Example 4.1:*** Let C(S) be a super special code given by

$$C(S) = \left[\begin{array}{c|c|c|c|c|c} C & C & C & C & C & C \\ \hline C_1 & C_1 & C_1 & C_1 & C_1 & C_1 \\ \hline C_1 & C_1 & C_1 & C_1 & C_1 & C_1 \\ \hline C & C & C & C & C & C \\ \hline C & C & C & C & C & C \\ \hline C_1 & C_1 & C_1 & C_1 & C_1 & C_1 \end{array}\right].$$

Here C is say a (n, k) code. $C_1$ is also only a (n, k) code but $C \neq C_1$. The codes C is assumed to carry the true messages. $C_1$ also carries messages but only to mislead the intruder so $C_1$'s can be called as misleading codes.

Now this key will be known to every one in the group who is sending or receiving the messages. Thus any one in the group only will be interested in the rows 1, 4 and 5 and ignore the messages in the rows 2, 3 and 6.

Thus when the number of rows and columns used are arbitrarily large it will be impossible for the intruder to guess at the codes and their by brake the key.

Similarly one can use a

$$C(S) = \left[\begin{array}{c|c|c|c|c|c} C_1 & C & C_1 & C_1 & C & C \\ \hline C_1 & C & C_1 & C_1 & C & C \\ \hline C_1 & C & C_1 & C_1 & C & C \\ \hline C_1 & C & C_1 & C_1 & C & C \\ \hline C_1 & C & C & C_1 & C & C \\ \hline C_1 & C & C_1 & C_1 & C & C \\ \hline C_1 & C & C_1 & C_1 & C & C \end{array}\right]$$

where C(S) is a 7 × 6 super special super code. Now both C and $C_1$ are only (n, $k_1$) codes or (n, $k_2$) codes $k_1 \neq k_2$. The group which uses this super special codes can agree upon to use the code C to carry the messages and $C_1$ are misleading codes. So anyone in this group will analyse only the codes in columns 2, 5 and 6 ignore columns 1, 3 and 4.

We give now the example of a block and misleading block code for the cryptographist.

***Example 4.2:*** Let the super special super code

$$C(S) = \left[\begin{array}{c|c|c|c|c|c|c|c|c} C_1 & C_1 & C_1 & C_2 & C_2 & C_3 & C_3 & C_3 & C_3 \\ \hline C_1 & C_1 & C_1 & C_2 & C_2 & C_3 & C_3 & C_3 & C_3 \\ \hline C_1 & C_1 & C_1 & C_2 & C_2 & C_8 & C_8 & C_8 & C_9 \\ \hline C_4 & C_4 & C_5 & C_2 & C_2 & C_8 & C_8 & C_8 & C_9 \\ \hline C_4 & C_4 & C_5 & C_{10} & C_{10} & C_8 & C_8 & C_8 & C_9 \\ \hline C_4 & C_4 & C_5 & C_{10} & C_{10} & C_8 & C_8 & C_8 & C_9 \\ \hline C_4 & C_4 & C_5 & C_{11} & C_{11} & C_{11} & C_{11} & C_{12} & C_{12} \\ \hline C_4 & C_4 & C_5 & C_{11} & C_{11} & C_{11} & C_{11} & C_{12} & C_{12} \\ \hline C_6 & C_6 & C_5 & C_{11} & C_{11} & C_{11} & C_{11} & C_{12} & C_{12} \\ \hline C_6 & C_6 & C_7 & C_7 & C_7 & C_7 & C_7 & C_{12} & C_{12} \end{array}\right].$$

The related super matrix which is the parity check matrix H(S) of C(S) would be of the form

$$
H(S) = \left[\begin{array}{c|c|c|c|c|c|c|c|c}
H_1 & H_1 & H_1 & H_2 & H_2 & H_3 & H_3 & H_3 & H_3 \\ \hline
H_1 & H_1 & H_1 & H_2 & H_2 & H_3 & H_3 & H_3 & H_3 \\ \hline
H_1 & H_1 & H_1 & H_2 & H_2 & H_8 & H_8 & H_8 & H_9 \\ \hline
H_4 & H_4 & H_5 & H_2 & H_2 & H_8 & H_8 & H_8 & H_9 \\ \hline
H_4 & H_4 & H_5 & H_{10} & H_{10} & H_8 & H_8 & H_8 & H_9 \\ \hline
H_4 & H_4 & H_5 & H_{10} & H_{10} & H_8 & H_8 & H_8 & H_9 \\ \hline
H_4 & H_4 & H_5 & H_{11} & H_{11} & H_{11} & H_{11} & H_{12} & H_{12} \\ \hline
H_4 & H_4 & H_5 & H_{11} & H_{11} & H_{11} & H_{11} & H_{12} & H_{12} \\ \hline
H_6 & H_6 & H_5 & H_{11} & H_{11} & H_{11} & H_{11} & H_{12} & H_{12} \\ \hline
H_6 & H_6 & H_7 & H_7 & H_7 & H_7 & H_7 & H_{12} & H_{12}
\end{array}\right]
$$

where each $H_i$ is the parity check matrix of a code $C_i$ , i = 1, 2, …, 12. Further all the code $C_1$, $C_4$, $C_5$, $C_6$, $C_7$, $C_{11}$, $C_8$, $C_{10}$, $C_2$, $C_3$, $C_{12}$, $C_9$ have the same length.

Further all codes given by $H_1$, $H_2$, $H_8$, $H_3$, $H_9$, $H_4$, $H_5$, $H_{11}$, $H_6$, $H_{12}$, $H_7$ and $H_{10}$ have the same number of check symbols. Now any super code word x(S) in C(S) would be of the form

$$
x(S) = \left[\begin{array}{c|c|c|c|c|c|c|c|c}
x_1^1 & x_2^1 & x_3^1 & x_1^2 & x_2^2 & x_1^3 & x_2^3 & x_3^3 & x_4^3 \\ \hline
x_4^1 & x_5^1 & x_6^1 & x_3^2 & x_4^2 & x_5^3 & x_6^3 & x_7^3 & x_8^3 \\ \hline
x_7^1 & x_8^1 & x_9^1 & x_5^2 & x_6^2 & x_1^8 & x_2^8 & x_3^8 & x_1^9 \\ \hline
x_1^4 & x_2^4 & x_1^5 & x_7^2 & x_8^2 & x_4^8 & x_5^8 & x_6^8 & x_2^9 \\ \hline
x_3^4 & x_4^4 & x_2^5 & x_1^{10} & x_2^{10} & x_7^8 & x_8^8 & x_9^8 & x_3^9 \\ \hline
x_5^4 & x_6^4 & x_3^5 & x_3^{10} & x_4^{10} & x_{10}^8 & x_{11}^8 & x_{12}^8 & x_4^9 \\ \hline
x_7^4 & x_8^4 & x_4^5 & x_1^{11} & x_2^{11} & x_3^{11} & x_4^{11} & x_5^{12} & x_1^{12} \\ \hline
x_9^4 & x_{10}^4 & x_5^5 & x_5^{11} & x_6^{11} & x_7^{11} & x_8^{11} & x_6^{12} & x_2^{12} \\ \hline
x_1^6 & x_2^6 & x_6^5 & x_9^{11} & x_{10}^{11} & x_{11}^{11} & x_{12}^{11} & x_7^{12} & x_3^{12} \\ \hline
x_3^6 & x_4^6 & x_1^7 & x_2^7 & x_3^7 & x_4^7 & x_5^7 & x_8^{12} & x_4^{12}
\end{array}\right]
$$

where $x_j^i \in C_i$; i = 1, 2, …, 12 and j varies according to the number of code words used. If y(S) is the received message

$H_i\left(y_p^i\right)^t$ = (0) will make the receiver accept it otherwise correct it or find the error using some techniques discussed earlier.

Now only few blocks are real blocks for which we have to work and other blocks are misleading blocks. Since all the codes have the same number of check symbols and the length of all the 12 codes are the same the intruder will not be in a position to make any form of guess and break the message. In fact he will not even be in a position to find out which of the blocks are misleading block of C(S) and which of them really carries the message.

Thus this provides a very high percentage of confidentiality and it is very difficult to know or break the message. This code when properly used will be a boon to the cryptography.

We give yet another example of a super special code C(S) which would be of use to the cryptographist.

***Example 4.3:*** Let C(S) be a super special super code;

$$C(S) = \left[\begin{array}{c|c|c|c|c|c|c|c|c|c|c}
C_1 & C_1 & C_1 & C_1 & C_1 & C_2 & C_2 & C_2 & C_2 & C_3 & C_3 \\ \hline
C & C & C & C & C & C_2 & C_2 & C_2 & C_2 & C & C_3 \\ \hline
C_1 & C_1 & C_1 & C_1 & C_1 & C & C_2 & C & C_2 & C_3 & C \\ \hline
C_1 & C_1 & C_1 & C_1 & C_1 & C_2 & C & C_2 & C_3 & C_3 & C_3 \\ \hline
C_5 & C_5 & C_5 & C_6 & C_6 & C_6 & C_6 & C & C & C & C \\ \hline
C & C & C & C_6 & C_6 & C_6 & C_4 & C & C & C_4 & C_4 \\ \hline
C_5 & C & C_5 & C & C & C & C_4 & C_4 & C & C_4 & C_4 \\ \hline
C & C_5 & C & C_6 & C_6 & C_6 & C & C & C_4 & C & C
\end{array}\right].$$

.

Now we have 7 sets of codes given by C, $C_1$, $C_2$, $C_3$, $C_4$, $C_5$ and $C_6$. All the seven codes are of same length and same number of check symbols and message symbols. Here only C is true, all the other 6 codes $C_1, C_2, \ldots, C_6$ are only misleading codes. The super special parity check matrix H(S) is given by

$$H(S) = \begin{bmatrix} H_1 & H_1 & H_1 & H_1 & H_1 & H_2 & H_2 & H_2 & H_2 & H_3 & H_3 \\ H & H & H & H & H & H_2 & H_2 & H_2 & H_2 & H & H_3 \\ H_1 & H_1 & H_1 & H_1 & H_1 & H & H_2 & H & H_2 & H_3 & H \\ H_1 & H_1 & H_1 & H_1 & H_1 & H_2 & H & H_2 & H_3 & H_3 & H_3 \\ H_5 & H_5 & H_5 & H_6 & H_6 & H_6 & H & H & H & H & H \\ H & H & H & H_6 & H_6 & H_6 & H_4 & H & H & H_4 & H_4 \\ H_5 & H & H_5 & H & H & H & H_4 & H_4 & H & H_4 & H_4 \\ H & H_5 & H & H_6 & H_6 & H_6 & H & H & H_4 & H & H \end{bmatrix}.$$

Here only the parity check matrix H gives the needed message, all other parity check matrices $H_1$, $H_2$, $H_3$, $H_4$, $H_5$ and $H_6$ need not be even known to the owner of this system of cryptography. Now everyone in group will be given a true/false chart as follows:

| | | | | | | | | | | |
|---|---|---|---|---|---|---|---|---|---|---|
| F | F | F | F | F | F | F | F | F | F | F |
| T | T | T | T | T | F | F | F | F | T | F |
| F | F | F | F | F | T | F | T | F | F | T |
| F | F | F | F | F | F | T | F | F | F | F |
| F | F | F | F | F | F | T | T | T | T | T |
| T | T | T | F | F | F | F | T | T | F | F |
| F | T | F | T | T | T | F | F | T | F | F |
| T | F | T | F | F | F | T | T | F | T | T |

or equivalently they can be supplied with a true chart or key.

| | | | | | | | | | | |
|---|---|---|---|---|---|---|---|---|---|---|
| | | | | | | | | | | |
| * | * | * | * | * | | | | | * | |
| | | | | | * | | * | | | * |
| | | | | | | * | | | | |
| | | | | | | * | * | * | * | * |
| * | * | * | | | | | * | * | | |
| | * | | * | * | * | | | * | | |
| * | | * | | | | * | * | | * | * |

.

The receiver would only decode the *'s of the table and will ignore the blanks. This will be known as the true chart of C(S).

These types of super special super codes C(S) can be thought of as special steganography. Unlike in a stegnanography where a secret message would be hidden with other messages here only the secret message is the essential message and all other messages are sent as misleading messages so that the intruder is never in a position to break the key or get to know the message.

We give some super special codes which are stegnanographic super special codes.

***Example 4.4:*** Let C(S) be a super special super code. We call this to be a steganographic super special code.

For instance we have a group which works with some name which starts in T and the messages sent to one another is highly confidential for it involves huge amount of money transactions or military secrets.

So if any intruder breaks open the message the company may run a very big loss or the nation's security would be in danger.

Now the super special code C(S) is given by

$$C(S)=\left[\begin{array}{c|c|c|c|c|c|c|c|c} C & C & C & C & C & C & C & C & C \\ \hline C & C & C & C & C & C & C & C & C \\ \hline C & C & C & C & C & C & C & C & C \\ \hline C & C & C & C & C & C & C & C & C \\ \hline C & C & C & C & C & C & C & C & C \\ \hline C & C & C & C & C & C & C & C & C \\ \hline C & C & C & C & C & C & C & C & C \end{array}\right]$$

where C is the code used. The hidden message is that the receiver is advised to use the truth table or a steganographic image of T in C(S) i.e.,

| | | | | | | | | |
|---|---|---|---|---|---|---|---|---|
| × | × | × | × | × | × | × | × | × |
| | | | | × | | | | |
| | | | | × | | | | |
| | | | | × | | | | |
| | | | | × | | | | |
| | | | | × | | | | |
| | | | | × | | | | |

i.e., the group is advised to read only the messages present in the first row and the 5th column of C(S) all other messages do not carry any sense to them for they are only misleading messages or the messages can take place as the first alphabet of every member of the group.

For instance K is the first alphabet of some member of the group then we need to decode the message from

| | | | | | |
|---|---|---|---|---|---|
| × | | | | | × |
| × | | | | × | |
| × | | | × | | |
| × | | × | | | |
| × | × | | | | |
| × | | × | | | |
| × | | | × | | |
| × | | | | × | |
| × | | | | | × |

i.e., first column and one of the opposite diagonals of the two 5 $\times$ 5 matrices so that letter K is formed.

This will be the key he has to use to decode the received message

$$y(S) = \left[\begin{array}{c|c|c|c|c|c}
y_1 & y_2 & y_3 & y_4 & y_5 & y_6 \\ \hline
y_7 & y_8 & y_9 & y_{10} & y_{11} & y_{12} \\ \hline
y_{13} & y_{14} & y_{15} & y_{16} & y_{17} & y_{18} \\ \hline
y_{19} & y_{20} & y_{21} & y_{22} & y_{23} & y_{24} \\ \hline
y_{25} & y_{26} & y_{27} & y_{28} & y_{29} & y_{30} \\ \hline
y_{31} & y_{32} & y_{33} & y_{34} & y_{35} & y_{36} \\ \hline
y_{37} & y_{38} & y_{39} & y_{40} & y_{41} & y_{42} \\ \hline
y_{43} & y_{44} & y_{45} & y_{46} & y_{47} & y_{48} \\ \hline
y_{49} & y_{50} & y_{51} & y_{52} & y_{53} & y_{54}
\end{array}\right].$$

He needs to decode only the messages $y_1$, $y_7$, $y_{13}$, $y_{19}$, $y_{25}$, $y_{31}$, $y_{37}$, $y_{43}$, $y_{49}$, $y_{26}$, $y_{21}$, $y_{16}$, $y_{11}$, $y_6$, $y_{33}$, $y_{40}$, $y_{47}$ and $y_{54}$ which is easily seen to form the letter K.

It can also be at times symbols like ‘cross’ or asterisk or star.

For instance if y(S) is the received message given by

$$y(S) = \left[\begin{array}{c|c|c|c|c|c|c}
y_1 & y_2 & y_3 & y_4 & y_5 & y_6 & y_7 \\ \hline
y_8 & y_9 & y_{10} & y_{11} & y_{12} & y_{13} & y_{14} \\ \hline
y_{15} & y_{16} & y_{17} & y_{18} & y_{19} & y_{20} & y_{21} \\ \hline
y_{22} & y_{23} & y_{24} & y_{25} & y_{26} & y_{27} & y_{28} \\ \hline
y_{29} & y_{30} & y_{31} & y_{32} & y_{33} & y_{34} & y_{35} \\ \hline
y_{36} & y_{37} & y_{38} & y_{39} & y_{40} & y_{41} & y_{42} \\ \hline
y_{43} & y_{44} & y_{45} & y_{46} & y_{47} & y_{48} & y_{49}
\end{array}\right]$$

each $y_i$ is a code word from the code C, $1 \leq i \leq 49$. The receiver should and need to decode only $y_4$ $y_{11}$ $y_{18}$ $y_{25}$ $y_{32}$ $y_{38}$ $y_{46}$ $y_{15}$ $y_{16}$ $y_{17}$ $y_{19}$ $y_{15}$ $y_{20}$ and $y_{21}$. This forms the cross. Thus these can also be given as finite series or a finite arithmetic progression for instance arithmetic progression with first term 4 and difference 7 last term 46 and arithmetic progression with first term 15 common difference 1 and last term 21. Any nice mathematical technique can be used to denote the codes which is essential for the member to be decoded to get the message.

**DEFINITION 4.1:** *Let C(S) be a super special super code if each of the codes $C_i$ in C(S) is cyclic then we call C(S) to be super special super cyclic code.*

We illustrate this by the following example.

***Example 4.5:*** Let C(S) be a super cyclic code. Let

$$C(S) = \left[\begin{array}{c|c} C_1^1 & C_2^1 \\ \hline C_1^2 & C_2^2 \end{array}\right]$$

where each $C_j^i$ is a cyclic code $1 \le i \le 2$ and $1 \le j \le 2$.

$$H(S) = \left[\begin{array}{c|c} H_1^1 & H_2^1 \\ \hline H_1^2 & H_2^2 \end{array}\right]$$

where

$$H_1^1 = \begin{bmatrix} 0 & 0 & 1 & 0 & 1 & 1 & 1 \\ 0 & 1 & 0 & 1 & 1 & 1 & 0 \\ 1 & 0 & 1 & 1 & 1 & 0 & 0 \end{bmatrix} = H_1^2$$

and

$$H_2^1 = \begin{bmatrix} 0 & 0 & 1 & 0 & 0 & 1 \\ 0 & 1 & 0 & 0 & 1 & 0 \\ 1 & 0 & 0 & 1 & 0 & 0 \end{bmatrix} = H_2^2$$

i.e.,

$$H(S) = \left[\begin{array}{ccccccc|cccccc} 0 & 0 & 1 & 0 & 1 & 1 & 1 & 0 & 0 & 1 & 0 & 0 & 1 \\ 0 & 1 & 0 & 1 & 1 & 1 & 0 & 0 & 1 & 0 & 0 & 1 & 0 \\ 1 & 0 & 1 & 1 & 1 & 0 & 0 & 1 & 0 & 0 & 1 & 0 & 0 \\ \hline 0 & 0 & 1 & 0 & 1 & 1 & 1 & 0 & 0 & 1 & 0 & 0 & 1 \\ 0 & 1 & 0 & 1 & 1 & 1 & 0 & 0 & 1 & 0 & 0 & 1 & 0 \\ 1 & 0 & 1 & 1 & 1 & 0 & 0 & 1 & 0 & 0 & 1 & 0 & 0 \end{array}\right]$$

is a super matrix given.

Given a super special cyclic code

$$x(S) = \left[\begin{array}{c|c} x_1^1 & x_2^1 \\ \hline x_1^2 & x_2^2 \end{array}\right]$$

where each $x_j^i$ is a cyclic code $1 \le i, j \le 2$. Any super code word x(S) in C(S) is of the form

$$\left[\begin{array}{c|c} x_1^1 & x_1^2 \\ \hline x_2^1 & x_2^2 \end{array}\right]$$

$$= \left[\begin{array}{ccccccc|cccccc} 1 & 0 & 0 & 0 & 1 & 1 & 0 & 1 & 1 & 0 & 1 & 1 & 0 \\ \hline 1 & 1 & 1 & 1 & 1 & 1 & 1 & 1 & 1 & 1 & 1 & 1 & 1 \end{array}\right]. \, H(S)\left[x(S)\right]^t$$

$$= \left[\begin{array}{c|c} H_1^1 & H_2^1 \\ \hline H_1^2 & H_2^2 \end{array}\right]\left[\begin{array}{c|c} x_1^1 & x_2^1 \\ \hline x_1^2 & x_2^2 \end{array}\right]^T$$

$$= \left[\begin{array}{c|c} H_1^1 & H_2^1 \\ \hline H_1^2 & H_2^2 \end{array}\right]\left[\begin{array}{c|c} \left(x_1^1\right)^t & \left(x_2^1\right)^t \\ \hline \left(x_1^2\right)^t & \left(x_2^2\right)^t \end{array}\right]$$

$$= \left[\begin{array}{c|c} H_1^1\left(x_1^1\right)^t & H_2^1\left(x_2^1\right)^t \\ \hline H_1^2\left(x_1^2\right)^t & H_2^2\left(x_2^2\right)^t \end{array}\right]$$

$$= \left[\begin{array}{ccc|ccc} 0 & 0 & 0 & 0 & 0 & 0 \\ \hline 0 & 0 & 0 & 0 & 0 & 0 \end{array}\right]$$

so $x(S) \in C(S)$.

Thus we see in general if C(S) is a super special super code then any super code $x(S) \in C(S)$ will be of the form

$$\left[\begin{array}{c|c|c|c} x_1^1 & x_1^2 & \cdots & x_1^n \\ \hline x_2^1 & x_2^2 & \cdots & x_2^n \\ \hline \vdots & \vdots & \cdots & \vdots \\ \hline x_m^1 & x_m^2 & \cdots & x_m^n \end{array}\right],$$

where each $x_i^j$ is a row vector $1 \leq i \leq n$ and $1 \leq j \leq m$. This type of super matrix from now on wards will be known as super row cell matrix i.e., each cell in x(S) is a row vector.

Now we say two super row cell matrices x(S) and y(S) are equal if and only if each $x_i^j = y_i^j$. We see two super row cell matrices are of same order if and only if number of row cells in x(S) is equal to number of row cells in y(S) and number of columns cells in x(S) is equal to number of column cells in y(S) and further the number of elements in the row cell $x_j^i$ is the same as the number of elements in the row cell $y_j^i$ ; $1 \leq i \leq n$ and $1 \leq j \leq m$ where

$$y(S) = \left[\begin{array}{c|c|c|c} y_1^1 & y_1^2 & \cdots & y_1^n \\ \hline y_2^1 & y_2^2 & \cdots & y_2^n \\ \hline \vdots & \vdots & \cdots & \vdots \\ \hline y_m^1 & y_m^2 & \cdots & y_m^n \end{array}\right].$$

We will illustrate this situation by examples.

***Example 4.6:*** Let

$$x(S) = \left[\begin{array}{c|c|c} x_1^1 & x_1^2 & x_1^3 \\ \hline x_2^1 & x_2^2 & x_2^3 \\ \hline x_3^1 & x_3^2 & x_3^3 \end{array}\right]$$

$$= \left[\begin{array}{c|c|c} 1\,1\,0\,0\,1 & 0\,1\,1\,1 & 1\,1\,1\,0\,0\,0\,1\,1 \\ \hline 1\,1\,0\,1\,0 & 1\,1\,1\,1 & 1\,1\,0\,0\,1\,1\,1\,1 \\ \hline 1\,1\,1\,1\,0 & 0\,1\,0\,1 & 0\,1\,0\,1\,0\,1\,1\,1 \end{array}\right]$$

be a super cell matrix.

$$y(S) = \left[\begin{array}{c|c|c} y_1^1 & y_1^2 & y_1^3 \\ \hline y_2^1 & y_2^2 & y_2^3 \\ \hline y_3^1 & y_3^2 & y_3^3 \end{array}\right]$$

be a super cell row matrix i.e., let

$$y(S) = \left[\begin{array}{c|c|c} 0\ 1\ 1\ 0\ 0 & 1\ 1\ 1\ 1 & 0\ 0\ 1\ 1\ 1\ 1\ 0\ 0 \\ \hline 0\ 1\ 0\ 0\ 1 & 1\ 0\ 1\ 0 & 0\ 1\ 0\ 1\ 0\ 1\ 1\ 1 \\ \hline 1\ 1\ 1\ 1\ 0 & 0\ 1\ 0\ 1 & 1\ 0\ 1\ 0\ 1\ 0\ 1\ 1 \end{array}\right]$$

we say y(S) and x(S) are of same order of same type. However y(S) ≠ x(S). Suppose

$$x(S) = \left[\begin{array}{c|c|c} 1\ 1\ 1\ 1\ 1 & 0\ 0\ 0\ 1 & 1\ 1\ 0\ 0\ 1\ 1\ 0\ 0 \\ \hline 0\ 0\ 0\ 1\ 1 & 1\ 0\ 0\ 0 & 1\ 0\ 0\ 0\ 0\ 0\ 0\ 1 \\ \hline 1\ 0\ 0\ 0\ 1 & 1\ 1\ 1\ 1 & 0\ 0\ 0\ 0\ 0\ 1\ 1\ 1 \end{array}\right]$$

and

$$y(S) = \left[\begin{array}{c|c|c} 1\ 1\ 1\ 1\ 1 & 0\ 0\ 0\ 1 & 1\ 1\ 0\ 0\ 1\ 1\ 0\ 0 \\ \hline 0\ 0\ 0\ 1\ 1 & 1\ 0\ 0\ 0 & 1\ 0\ 0\ 0\ 0\ 0\ 0\ 1 \\ \hline 1\ 0\ 0\ 0\ 1 & 1\ 1\ 1\ 1 & 0\ 0\ 0\ 0\ 0\ 1\ 1\ 1 \end{array}\right]$$

x(S) = y(S) if each row cell in them are identical. Thus x(S) = y(S). Now x(S), y(S) can also be called as super row cell vectors. Now we define the dot product of two super cell vectors if and only if they are of same order i.e., each row vector in the corresponding cells of x(S) and y(S) are of same length. Let

$$x(S) = \left[\begin{array}{c|c|c|c} x_1^1 & x_1^2 & \cdots & x_1^m \\ \hline x_2^1 & x_2^2 & \cdots & x_2^m \\ \hline \vdots & \vdots & \cdots & \vdots \\ \hline x_n^1 & x_n^2 & \cdots & x_n^m \end{array}\right]$$

and

$$y(S) = \left[\begin{array}{c|c|c|c} y_1^1 & y_1^2 & \cdots & y_1^m \\ \hline y_2^1 & y_2^2 & \cdots & y_2^m \\ \hline \vdots & \vdots & \cdots & \vdots \\ \hline y_n^1 & y_n^2 & \cdots & y_n^m \end{array}\right].$$

Then we define the dot product of x(S) with y(S) as

$$(x(S), y(S)) = \left[\begin{array}{c|c|c|c} (x_1^1, y_1^1) & (x_1^2, y_1^2) & \cdots & (x_1^m, y_1^m) \\ \hline (x_2^1, y_2^1) & (x_2^2, y_2^2) & \cdots & (x_2^m, y_2^m) \\ \hline \vdots & \vdots & \cdots & \vdots \\ \hline (x_n^1, y_n^1) & (x_n^2, y_n^2) & \cdots & (x_n^m, y_n^m) \end{array}\right]$$

where $(x_j^i, y_j^i)$ is the usual inner (dot) product of two vectors $1 \le i \le m$, $1 \le j \le n$. Basically all these row vectors are from the subspaces of the vector spaces.

Now we illustrate this situation by the following example.

***Example 4.7:*** Let

$$x(S) = \left[\begin{array}{c|c} 110 & 111101 \\ \hline 111 & 011101 \\ \hline 001 & 100010 \\ \hline 010 & 011001 \end{array}\right]$$

and

$$y(S) = \left[\begin{array}{c|c} 010 & 110001 \\ \hline 101 & 100011 \\ \hline 011 & 101010 \\ \hline 110 & 010101 \end{array}\right]$$

be any two super cell row vectors of same order. Now how does the dot or inner product of x(S) with y(S) look like

$$\left(x(S), y(S)\right) = \left[\begin{array}{c|c} \left((1\,1\,0),(0\,1\,0)\right) & (1\,1\,1\,1\,0\,1),(1\,1\,0\,0\,0\,1) \\ \hline \left((1\,1\,1),(1\,0\,1)\right) & \left((0\,1\,1\,1\,0\,1),(1\,0\,0\,0\,1\,1)\right) \\ \hline \left((0\,0\,1),(0\,1\,1)\right) & \left((1\,0\,0\,0\,1\,0),(1\,0\,1\,0\,1\,0)\right) \\ \hline \left((0\,1\,0),(1\,1\,0)\right) & \left((0\,1\,1\,0\,0\,1),(0\,1\,0\,1\,0\,1)\right) \end{array}\right]$$

$$= \left[\begin{array}{c|c} 1 & 1 \\ \hline 0 & 1 \\ \hline 1 & 0 \\ \hline 1 & 0 \end{array}\right]$$

is just a super matrix which can be called as a super cell matrix as usual $4 \times 2$ matrix is divided as cells.

Now having defined the notion of inner product of super row cell vectors or matrices, we proceed to define super orthogonal super code.

**DEFINITION 4.2:** *Let C(S) be a super special super binary code i.e., the code words are from the field of characteristics two i.e., from $Z_2 = \{0, 1\}$. Let*

$$x(S) = \left[\begin{array}{c|c|c|c} x_1^1 & x_1^2 & \cdots & x_1^m \\ \hline x_2^1 & x_2^2 & \cdots & x_2^m \\ \hline \vdots & \vdots & \cdots & \vdots \\ \hline x_n^1 & x_n^2 & \cdots & x_n^m \end{array}\right]$$

*be a super code word in C(S). Suppose there exists a super code word y(S) in C(S) where*

$$y(S) = \left[\begin{array}{c|c|c|c} y_1^1 & y_1^2 & \cdots & y_1^m \\ \hline y_2^1 & y_2^2 & \cdots & y_2^m \\ \hline \vdots & \vdots & \cdots & \vdots \\ \hline y_n^1 & y_n^2 & \cdots & y_n^m \end{array}\right]$$

*is such that*

$$\left(x(S), y(S)\right) = \left[\begin{array}{c|c|c|c} 0 & 0 & \cdots & 0 \\ \hline 0 & 0 & \cdots & 0 \\ \hline \vdots & \vdots & \cdots & \vdots \\ \hline 0 & 0 & \cdots & 0 \end{array}\right],$$

*then we say the super code words are orthogonal to each other. Now suppose C(S) is a super special super code then the set of all super codes {w(S) / (w(S), (x(S)) = (0)(S) for all x(S) $\in$ C(S)} is defined to be the super special orthogonal code of the super special code C(S).*

*Here w(S) $\subseteq$ {The collection of all n $\times$ m super cell row vectors of same order as x(S) with entries from {0, 1}} =V(S). In fact this collection can be realised as super special vector space over $Z_2$ = {0, 1}.We denote this super special orthogonal super code by $(C(S))^{\perp}$. We see 0(S) is the only code word in C(S), then the whole space V(S) is orthogonal to 0(S).*

We will illustrate this situation by the following example.

***Example 4.8:*** Let

$$\left\{\left[\begin{array}{c|c} 0\,0\,0\,0 & 0\,0\,0\,0 \\ \hline 0\,0\,0\,0 & 0\,0\,0\,0 \end{array}\right], \left[\begin{array}{c|c} 0\,0\,0\,0 & 0\,0\,0\,0 \\ \hline 0\,0\,0\,0 & 1\,0\,1\,1 \end{array}\right] \left[\begin{array}{c|c} 0\,0\,0\,0 & 0\,0\,0\,0 \\ \hline 1\,0\,1\,1 & 0\,0\,0\,0 \end{array}\right],\right.$$

$$\left[\begin{array}{c|c} 1\,0\,1\,1 & 0\,0\,0\,0 \\ \hline 0\,0\,0\,0 & 0\,0\,0\,0 \end{array}\right]; \left[\begin{array}{c|c} 0\,0\,0\,0 & 1\,0\,1\,1 \\ \hline 0\,0\,0\,0 & 0\,0\,0\,0 \end{array}\right],$$

$$\left[\begin{array}{c|c} 0\,0\,0\,0 & 0\,1\,0\,1 \\ \hline 0\,0\,0\,0 & 0\,0\,0\,0 \end{array}\right], \left[\begin{array}{c|c} 0\,0\,0\,0 & 0\,0\,0\,0 \\ \hline 0\,0\,0\,0 & 0\,1\,0\,1 \end{array}\right],$$

$$\left.\left[\begin{array}{c|c} 0\,0\,0\,0 & 0\,0\,0\,0 \\ \hline 0\,0\,0\,0 & 1\,1\,1\,0 \end{array}\right], \left[\begin{array}{c|c} 0\,0\,0\,0 & 0\,0\,0\,0 \\ \hline 1\,1\,1\,0 & 0\,0\,0\,0 \end{array}\right], \left[\begin{array}{c|c} 1\,1\,1\,0 & 0\,1\,0\,1 \\ \hline 1\,0\,1\,1 & 0\,0\,0\,0 \end{array}\right]\right\}$$

Thus

$$C(S) = \left[\begin{array}{c|c} C & C \\ \hline C & C \end{array}\right]$$

where C is the linear binary (4, 2) code given by {(0 0 0 0 ), (1 0 1 1), (0 1 0 1), (1 1 1 0)}. $C^{\perp}$ = {(0 0 0 0), (1 1 0 1), (0 1 1 1), (1 0 1 0)} is the orthogonal code of C.

The super special orthogonal code of C(S) given by $C(S)^{\perp}$ is given by

$$(C(S)^{\perp}) = \left[\begin{array}{c|c} C^{\perp} & C^{\perp} \\ \hline C^{\perp} & C^{\perp} \end{array}\right].$$

Clearly

$$(C(S)), (C(S)^{\perp}) = 0(S).$$

We give a general method of finding $(C(S)^{\perp})$ given C(S).
Let

$$(C(S)) = \left[\begin{array}{c|c|c|c} C_1^1 & C_2^1 & \cdots & C_n^m \\ \hline C_1^2 & C_2^2 & \cdots & C_2^m \\ \hline \vdots & \vdots & \cdots & \vdots \\ \hline C_n^1 & C_n^2 & \cdots & C_n^m \end{array}\right];$$

$C_j^i$ are linear codes $1 \le i \le m$ and $1 \le j \le n$.

$$\left(C(S)\right)^{\perp} = \left[\begin{array}{c|c|c|c} \left(C_1^1\right)^{\perp} & \left(C_2^1\right)^{\perp} & \cdots & \left(C_1^m\right)^{\perp} \\ \hline \left(C_1^2\right)^{\perp} & \left(C_2^2\right)^{\perp} & \cdots & \left(C_2^m\right)^{\perp} \\ \hline \vdots & \vdots & \cdots & \vdots \\ \hline \left(C_n^1\right)^{\perp} & \left(C_n^2\right)^{\perp} & \cdots & \left(C_n^m\right)^{\perp} \end{array}\right]$$

$$(C(S)), (C(S)^{\perp}) = 0(S).$$

Now having defined cyclic orthogonal super special super codes C(S) we just indicate how error is detected and corrected.

We just illustrate this by the following example.

***Example 4.9:*** Let

$$\left(C(S)\right)=\left[\begin{array}{c|c|c|c} C_1^1 & C_2^1 & \cdots & C_1^m \\ \hline C_1^2 & C_2^2 & \cdots & C_2^m \\ \hline \vdots & \vdots & \cdots & \vdots \\ \hline C_n^1 & C_n^2 & \cdots & C_n^m \end{array}\right]$$

be a super special super code.

Let the super parity check matrix H(S) be associated with it, that is

$$H(S)=\left[\begin{array}{c|c|c|c} H_1^1 & H_2^1 & \cdots & H_1^m \\ \hline H_1^2 & H_2^2 & \cdots & H_2^m \\ \hline \vdots & \vdots & \cdots & \vdots \\ \hline H_n^1 & H_n^2 & \cdots & H_n^m \end{array}\right].$$

Any super special code word x(S) $\in$ C(S) is of the form

$$x(S)=\left[\begin{array}{c|c|c|c} x_1^1 & x_1^2 & \cdots & x_1^m \\ \hline x_2^1 & x_2^2 & \cdots & x_2^m \\ \hline \vdots & \vdots & \cdots & \vdots \\ \hline x_n^1 & x_n^2 & \cdots & x_n^m \end{array}\right].$$

If $H(S)(x(S))^t = 0(S)$ then we assume x(S) $\in$ C(S).

If $H(S)(x(S))^t \neq 0(S)$ for a received code word x(S) then we assume the received message has an error.

Thus we have

$$H(S)\left[\left(x(S)\right)\right]^t$$

$$= \left[\begin{array}{c|c|c|c} H_1^1(x_1^1)^t & H_1^2(x_1^2)^t & \cdots & H_1^m(x_1^m)^t \\ \hline H_2^1(x_2^1)^t & H_2^2(x_2^2)^t & \cdots & H_2^m(x_2^m)^t \\ \hline \vdots & \vdots & \cdots & \vdots \\ \hline H_n^1(x_n^1)^t & H_n^2(x_n^2)^t & \cdots & H_n^m(x_n^m)^t \end{array}\right] = (0(S));$$

$x(S) \in C(S)$ otherwise we use the usual error correcting techniques to each cell i.e., to each $H_j^i(x_j^i)^t$ to obtain the corrected code word; $1 \leq i \leq n$; $1 \leq j \leq m$.

Thus we have indicated how these super special codes can be used when ARQ is impossible. Also these codes will be very useful in cryptology. Further use of these codes can spare both time and economy. The applications of super special row codes and super special column codes can be obtained from super special codes with appropriate modifications.

# FURTHER READING

# INDEX

**G**

**H**

**I**

**L**

**M**

**O**

**P**

**R**

**S**

# ABOUT THE AUTHORS

**Dr.W.B.Vasantha Kandasamy** is an Associate Professor in the Department of Mathematics, Indian Institute of Technology Madras, Chennai. In the past decade she has guided 13 Ph.D. scholars in the different fields of non-associative algebras, algebraic coding theory, transportation theory, fuzzy groups, and applications of fuzzy theory of the problems faced in chemical industries and cement industries.

She has to her credit 646 research papers. She has guided over 68 M.Sc. and M.Tech. projects. She has worked in collaboration projects with the Indian Space Research Organization and with the Tamil Nadu State AIDS Control Society. She is presently working on a research project funded by the Board of Research in Nuclear Sciences, Government of India. This is her 47th book.

On India's 60th Independence Day, Dr.Vasantha was conferred the Kalpana Chawla Award for Courage and Daring Enterprise by the State Government of Tamil Nadu in recognition of her sustained fight for social justice in the Indian Institute of Technology (IIT) Madras and for her contribution to mathematics. The award, instituted in the memory of Indian-American astronaut Kalpana Chawla who died aboard Space Shuttle Columbia, carried a cash prize of five lakh rupees (the highest prize-money for any Indian award) and a gold medal.
She can be contacted at vasanthakandasamy@gmail.com
Web Site: http://mat.iitm.ac.in/home/wbv/public_html/

---

**Dr. Florentin Smarandache** is a Professor of Mathematics at the University of New Mexico in USA. He published over 75 books and 150 articles and notes in mathematics, physics, philosophy, psychology, rebus, literature.

In mathematics his research is in number theory, non-Euclidean geometry, synthetic geometry, algebraic structures, statistics, neutrosophic logic and set (generalizations of fuzzy logic and set respectively), neutrosophic probability (generalization of classical and imprecise probability). Also, small contributions to nuclear and particle physics, information fusion, neutrosophy (a generalization of dialectics), law of sensations and stimuli, etc. He can be contacted at smarand@unm.edu

---

**K. Ilanthenral** is the editor of The Maths Tiger, Quarterly Journal of Maths. She can be contacted at ilanthenral@gmail.com